\documentclass{amsart}
\begin{document}
\title[Arithmetic Lie theory]{Arithmetic differential equations  on $GL_n$, II:\\
arithmetic Lie theory}
\author{Alexandru Buium and Taylor Dupuy}
\def \cB{{\mathcal B}}
\def \dim{\text{dim}}
\def \Rp{R_p}
\def \Rpi{R_{\pi}}
\def \dpi{\d_{\pi}}
\def \bT{{\bf T}}
\def \cI{{\mathcal I}}
\def \cH{{\mathcal H}}
\def \cJ{{\mathcal J}}
\def \ZN{\bZ[1/N,\zeta_N]}
\def \tA{\tilde{A}}
\def \o{\omega}
\def \tB{\tilde{B}}
\def \tC{\tilde{C}}
\def \alph{A}
\def \bet{B}
\def \bsigma{\bar{\sigma}}
\def \y{^{\infty}}
\def \Ra{\Rightarrow}
\def \uBS{\overline{BS}}
\def \lBS{\underline{BS}}
\def \lB{\underline{B}}
\def \<{\langle}
\def \>{\rangle}
\def \hL{\hat{L}}
\def \cU{\mathcal U}
\def \cF{\mathcal F}
\def \S{\Sigma}
\def \st{\stackrel}
\def \sd{Spec_{\d}\ }
\def \pd{Proj_{\d}\ }
\def \s{\sigma_2}
\def \i{\sigma_1}
\def \bs{\bigskip}
\def \cD{\mathcal D}
\def \cC{\mathcal C}
\def \cT{\mathcal T}
\def \cK{\mathcal K}
\def \cX{\mathcal X}
\def \sX{X_{set}}
\def \cY{\mathcal Y}
\def \cS{X}
\def \cR{\mathcal R}
\def \cE{\mathcal E}
\def \tcE{\tilde{\mathcal E}}
\def \cP{\mathcal P}
\def \cA{\mathcal A}
\def \cV{\mathcal V}
\def \cM{\mathcal M}
\def \cL{\mathcal L}
\def \cN{\mathcal N}
\def \tcM{\tilde{\mathcal M}}
\def \caS{\mathcal S}
\def \cG{\mathcal G}
\def \cB{\mathcal B}
\def \tG{\tilde{G}}
\def \cF{\mathcal F}
\def \h{\hat{\ }}
\def \hp{\hat{\ }}
\def \tS{\tilde{S}}
\def \tP{\tilde{P}}
\def \tA{\tilde{A}}
\def \tX{\tilde{X}}
\def \tcS{\tilde{X}}
\def \tT{\tilde{T}}
\def \tE{\tilde{E}}
\def \tV{\tilde{V}}
\def \tC{\tilde{C}}
\def \tI{\tilde{I}}
\def \tU{\tilde{U}}
\def \tG{\tilde{G}}
\def \tu{\tilde{u}}
\def \chu{\check{u}}
\def \tx{\tilde{x}}
\def \tL{\tilde{L}}
\def \tY{\tilde{Y}}
\def \d{\delta}
\def \e{\chi}
\def \bW{\mathbb W}
\def \bV{{\mathbb V}}
\def \bF{{\bf F}}
\def \bE{{\bf E}}
\def \bC{{\bf C}}
\def \bO{{\bf O}}
\def \bR{{\bf R}}
\def \bA{{\bf A}}
\def \bB{{\bf B}}
\def \cO{\mathcal O}
\def \ra{\rightarrow}
\def \bx{{\bf x}}
\def \f{{\bf f}}
\def \bX{{\bf X}}
\def \bH{{\bf H}}
\def \bS{{\bf S}}
\def \bF{{\bf F}}
\def \bN{{\bf N}}
\def \bK{{\bf K}}
\def \bE{{\bf E}}
\def \bB{{\bf B}}
\def \bQ{{\bf Q}}
\def \bd{{\bf d}}
\def \bY{{\bf Y}}
\def \bU{{\bf U}}
\def \bL{{\bf L}}
\def \bQ{{\bf Q}}
\def \bP{{\bf P}}
\def \bR{{\bf R}}
\def \bC{{\bf C}}
\def \bD{{\bf D}}
\def \bM{{\bf M}}
\def \bZ{{\mathbb Z}}
\def \xtoleqr{x^{(\leq r)}}
\def \hU{\hat{U}}
\def \k{\kappa}
\def \ee{\overline{p^{\k}}}

\newtheorem{THM}{{\!}}[section]
\newtheorem{THMX}{{\!}}
\renewcommand{\theTHMX}{}
\newtheorem{theorem}{Theorem}[section]
\newtheorem{corollary}[theorem]{Corollary}
\newtheorem{lemma}[theorem]{Lemma}
\newtheorem{proposition}[theorem]{Proposition}
\theoremstyle{definition}
\newtheorem{definition}[theorem]{Definition}
\theoremstyle{remark}
\newtheorem{remark}[theorem]{Remark}
\newtheorem{example}[theorem]{\bf Example}
\numberwithin{equation}{section}
\address{University of New Mexico \\ Albuquerque, NM 87131}
\email{buium@math.unm.edu, taylor.dupuy@gmail.com} 
\maketitle

\begin{abstract}

Motivated by the search of a concept of linearity in the theory of arithmetic differential equations \cite{book} we introduce here an arithmetic analogue of Lie algebras, of Chern connections, and of Maurer-Cartan connections. 
Our arithmetic analogues of Chern connections are  certain  remarkable lifts of Frobenius on the $p$-adic completion of $GL_n$ which are uniquely determined by certain compatibilities  with the ``outer" involutions defined by symmetric (respectively antisymmetric) matrices. The Christoffel symbols of our arithmetic Chern connections will involve a matrix analogue of the Legendre symbol. The analogues of Maurer-Cartan connections can then be viewed as families of ``linear" flows  attached to each of our Chern connections.
We will also investigate the compatibility of lifts of Frobenius with the inner automorphisms  of $GL_n$; in particular we will prove the existence and uniqueness of certain arithmetic analogues of ``isospectral flows" on the space of matrices. 
\end{abstract}

\section{Main concepts and results}

\subsection{Introduction} This paper is, in principle,  part of a series of papers \cite{adel1, adel3} but can be read independently of the other papers in the series. 
This series of papers is motivated, in particular, by the problem of finding an arithmetic analogue of linear differential equations.
 A general theory of arithmetic differential equations was introduced in \cite{char, book}; we will not assume here any familiarity with \cite{char, book} but let us recall that the idea in loc. cit.  was to replace derivation operators with Fermat quotient operators. However, the concept of linearity in this theory seems elusive and turns out to be subtle. A naive attempt, modeled on Kolchin's logarithmic derivative \cite{kolchin} (equivalently on the Maurer-Cartan connection \cite{sharpe}),  would be to consider {\it classical differential cocycles} of $GL_n$ into its (usual) Lie algebra ${\mathfrak g}{\mathfrak l}_n$, i.e. cocycles (in the usual sense) given by ``$\d$-functions" (in the sense of \cite{char,book}); but it turns out \cite{adel1} that no  such cocycles exist (except for those that  ``factor through the determinant").  To remedy the situation we will introduce here an arithmetic analogue of Lie algebras, arithmetic analogues of the Maurer-Cartan connections, and a notion of linearity for arithmetic differential equations relative to a given lift of Frobenius on the $p$-adic completion of $GL_n$. We will then construct some remarkable such lifts; they will be uniquely determined by certain compatibilities with the involutions defining  the (various forms of  the) split classical groups $SL_n, SO_n, Sp_{n}$. In their turn these compatibilities are analogous to those defining  Chern connections in differential geometry. On the other hand  the Christoffel symbols giving our arithmetic Chern connections will involve a   matrix analogue of the Legendre symbol.   The ``linear" arithmetic differential equations attached to these lifts  admit remarkable ``prime integrals" and ``symmetries" which will be used later in \cite{adel3} when we discuss the Galois side of the theory. The above story pertains to the  ``outer" involutions of $GL_n$ and  their corresponding classical symmetric spaces; we will also investigate, in the present paper, the behavior of lifts of Frobenius with respect to inner involutions (or more generally inner automorphisms) of $GL_n$ and we will study  lifts of Frobenius on the corresponding conjugacy classes. In particular we will prove the existence and uniqueness of certain arithmetic analogues of ``isospectral flows" on the space of matrices.

\subsection{Acknowledgement}
The authors are indebted to P. Cartier, C. Boyer, and D. Vassiliev  for inspiring discussions. Also the first author would like to acknowledge partial support from the Hausdorff Institute of Mathematics in Bonn, from the NSF through grant DMS 0852591,  from the Simons Foundation
(award 311773), 
and from the Romanian National Authority
for Scientific Research, CNCS - UEFISCDI, project number
PN-II-ID-PCE-2012-4-0201.

\subsection{Arithmetic differential equations}
 To explain the main results of this paper let us quickly recall the   $\delta$-arithmetic setting of \cite{char, book}.
 We denote by $R$  the unique complete discrete valuation ring with maximal ideal generated by an odd prime $p$ and  residue field $k=R/pR$ equal to the algebraic closure ${\mathbb F}_p^a$ of ${\mathbb F}_p$; then we denote by $\d:R\ra R$ the map
 $$\d x=\frac{\phi(x)-x^p}{p}$$ where $\phi:R\ra R$ is the unique ring homomorphism lifting the $p$-power Frobenius on the residue field $k$. We refer to $\d$ as a {\it $p$-derivation}; it is, of course, not a derivation. We denote by 
$R^{\d}$ the monoid of constants $\{\lambda\in R;\d \lambda=0\}$; so $R^{\d}$ consists of $0$ and all roots of unity in $R$. Also we denote by $K$ the fraction field of $R$.
 In this context we will consider (smooth) schemes of finite type $X$ over $R$ or, more generally,  (smooth) {\it $p$-formal schemes} of finite type, by which we mean formal schemes locally isomorphic to $p$-adic completions of (smooth) schemes of finite type; we denote by $X(R)$ the set of $R$-points of $X$; if there is no danger of confusion we often simply write $X$ in place of $X(R)$. Groups in the category of smooth $p$-formal schemes will be called {\it smooth group $p$-formal schemes}.
 For any ring $A$ we denote by $\widehat{A}$ the $p$-adic completion of $A$. In particular we consider the ring $R[x_0,...,x_n]\h$, $p$-adic completion of the ring of polynomials; its elements are the restricted power series with coefficients in $R$.
A map $f:R^N\ra R^M$
will be called a {\it $\d$-map} of order $n$ if there exists an $M$-vector  $F=(F_j)$,  $F_j\in R[x_0,...,x_n]\h$ of restricted power series with coefficients in $R$ in $N$-tuples of variables $x_0,...,x_n$, such that 
\begin{equation}
\label{singe}
f(a)=F(a,\d a,...,\d^n a),\ \ a\in R^N.\end{equation}
 One can then consider affine smooth schemes $X,Y$ and say that a map $f:X(R)\ra Y(R)$ is a {\it $\d$-map} of order $n$ if there exist  embeddings $X\subset {\mathbb A}^N$, $Y\subset {\mathbb A}^M$ such that $f$ is induced by a $\d$-map $R^N\ra R^M$ of order $n$; we simply write $f:X\ra Y$. More generally if $X$ is any smooth scheme and $Y$ is an affine scheme a set theoretic map $X\ra Y$ is called a {\it $\d$-map of order $n$} if there exists an affine cover $X=\bigcup X_i$ such that all the induced maps $X_i\ra Y$ are $\d$-maps of order $n$.
 If $G,H$ are smooth group schemes a {\it $\d$-homomorphism} $G\ra H$ is, by definition,  a $\d$-map which is also a homomorphism. An {\it arithmetic differential equation} (or simply a $\d$-{\it equation}) is an equation of the form $f(u)=0$ where $f:X\ra {\mathbb A}^N$ is a $\d$-map and $u\in X(R)$. All of the above concepts can be generalized, in an obvious way, to the case when schemes are replaced by $p$-formal schemes.
 We shall review these concepts, from a slightly different, but equivalent perspective,  in section 2 following 
\cite{char,book}.

\subsection{Arithmetic analogue of Lie algebra}
 In the present paper we will introduce  arithmetic analogues of some basic concepts  of Lie theory.  The theory will be  introduced for arbitrary linear algebraic groups
but here is how it looks in the case of $GL_n$ and its subgroups. We start by discussing the arithmetic analogue of the Lie algebra of a linear algebraic group.

Let $GL_n=Spec\ R[x,\det(x)^{-1}]$, where $x$ is a matrix of indeterminates.
The {\it $\d$-Lie algebra} of $GL_n$  is defined as the set ${\mathfrak g}{\mathfrak l}_n$  of $n\times n$ matrices over $R$ equipped with the non-commutative group law $+_{\d}:{\mathfrak g}{\mathfrak l}_n\times {\mathfrak g}{\mathfrak l}_n\ra {\mathfrak g}{\mathfrak l}_n$ given by
$$a+_{\d}b:=a+b+pab,$$
where the addition and multiplication in the right hand side are those of ${\mathfrak g}{\mathfrak l}_n$, viewed as an associative $R$-algebra. (This makes ${\mathfrak g}{\mathfrak l}_n$ the group of $R$-points of a natural group  $p$-formal scheme.) We denote by $-_{\d}$ the inverse with respect to $+_{\d}$.
There is a natural {\it $\d$-adjoint action} $\star_{\d}$ of $GL_n$ on ${\mathfrak g}{\mathfrak l}_n$ given by
$$a\star_{\d} b:=\phi(a) \cdot b \cdot \phi(a)^{-1}.$$
Here $\phi(a)$ is the matrix obtained from $a$ by applying $\phi:R\ra R$ to the coefficients.
Let now $G$ be a smooth subgroup scheme of $GL_n$.
Assume the ideal of $G$ in $\cO(GL_n)$ is generated by polynomials $f_i(x)\in R[x]$.
Then recall that the Lie algebra $L(G)$ identifies, as an additive group, to the  subgroup of the usual additive group $({\mathfrak g}{\mathfrak l}_n,+)$ consisting of all matrices $a$ satisfying
 $$``\epsilon^{-1}"f_i(1+\epsilon a)=0,$$
where $\epsilon^2=0$. Let
 $f^{(\phi)}_i$ be the polynomials obtained from $f_i$ by applying  $\phi$ to the coefficients.
Then 
we  define the {\it $\d$-Lie algebra} $L_{\d}(G)$ as the  subgroup of $({\mathfrak g}{\mathfrak l}_n,+_{\d})$ consisting of all the matrices $a\in {\mathfrak g}{\mathfrak l}_n$ satisfying
$$p^{-1}f_i^{(\phi)}(1+pa)=0.$$
In this formulation the factor $p^{-1}$ is, of course, not necessary;
but it will be important later to consider the above equations, with $p^{-1}$ present, in order to define a group $p$-formal scheme whose group of $R$-points is $L_{\d}(G)$.

 Note that $L_{\d}(G)$ and  $L(G)$ are not equal as subsets of $L(GL_n)={\mathfrak g}{\mathfrak l}_n$. For instance if $G=SL_2\subset GL_2$ then $L(SL_2)={\mathfrak sl}_2$ consists of all $a\in {\mathfrak g}{\mathfrak l}_2$ with $\text{tr}(a)=0$ whereas $L_{\d}(SL_2)$ consists of all $a\in {\mathfrak g}{\mathfrak l}_2$ such that $\text{tr}(a)+p\cdot \det(a)=0.$
 
 There is an arithmetic analogue of the Lie bracket  which can be naturally introduced at this point. This will play a role in subsequent work but will only play a minor role in the present paper; so the reader may choose to skip the definition below. For the definition of the arithmetic Lie bracket we need to first introduce a higher order version of the above construction as follows. For $r\geq 1$ an integer the {\it order $r$  $\d$-Lie algebra} of $GL_n$
 is the group $L_{\d}^{r}(GL_n)$ whose underlying set is, again, ${\mathfrak g}{\mathfrak l}_n$, and whose group law $+_{\d,r}$ is given by 
 $$a+_{\d,r}b=a+b+p^{r}ab.$$
 If $r$ is clear from the context we write $+_{\d,*}$ instead of $+_{\d,r}$.
 If $G$ is a subgroup scheme of $GL_n$ with ideal generated by $f_i(x)$   we define the  {\it order $r$ $\d$-Lie algebra} $L_{\d}^{r}(G)$ as the  subgroup of $({\mathfrak g}{\mathfrak l}_n,+_{\d,r})$ consisting of all the matrices $a\in {\mathfrak g}{\mathfrak l}_n$ satisfying
$$p^{-r}f_i^{(\phi^{r})}(1+p^{r}a)=0.$$
So $L_{\d}(G)=L_{\d}^1(G)$. Now, for $r,s\geq 1$, define the {\it $\d$-Lie brackets}
$$[\ ,\ ]_{\d}:L_{\d}^{r}(GL_n)\times L_{\d}^{s}(GL_n)\ra L_{\d}^{r+s}(GL_n)$$
by the formula
\begin{equation}
\label{snuggles}
[\alpha,\beta]_{\d}  := \frac{(1+p^{r}\phi^{s}(\alpha))(1+p^{s}\phi^{r}(\beta))(1+p^{r}\phi^{s}(\alpha))^{-1}(1+p^{s}\phi^{r}(\beta))^{-1}-1}{p^{r+s}}.
\end{equation}
This map is a $\d$-map of order $\max\{r,s\}$ and, for any subgroup scheme $G\subset GL_n$, it induces a $\d$-map of $p$-formal schemes,
$$[\ ,\ ]_{\d}:L_{\d}^{r}(G)\times L_{\d}^{s}(G)\ra L_{\d}^{r+s}(G).$$
These brackets are easily seen to satisfy the following congruences:
$$[\alpha,\beta]_{\d} \equiv \alpha^{(p^{s})} \beta^{(p^{r})}- \beta^{(p^{r})}\alpha^{(p^{s})}\ \ \text{mod}\ \ p,$$
for $\alpha\in L_{\d}^r(G)$, $\beta\in L_{\d}^s(G)$, 
where the $(p^{r})$ superscript means ``raising all entries of the corresponding matrix to the $p^{r}$-th power". As a consequence one has linearity, antisymmetry and the Jacobi identity mod $p$:

\medskip

$1) \ \ [\alpha_1+_{\d,*}\alpha_2,\beta]_{\d}\equiv [\alpha_1,\beta]_{\d}+_{\d,*}[\alpha_2,\beta]_{\d}\ \ \ \text{mod}\ \ \ p$

\medskip

$2) \ \ [\alpha,\beta]_{\d}+_{\d,*}[\beta,\alpha]_{\d} \equiv 0\ \ \ \text{mod}\ \ \ p$

\medskip

$3)\ \  [[\alpha,\beta]_{\d},\gamma]_{\d}+_{\d,*}[[\beta,\gamma]_{\d},\alpha]_{\d}+_{\d,*}
[[\beta,\gamma]_{\d},\alpha]_{\d}\equiv 0,\ \ \ \text{mod}\ \ \ p$

\medskip

\noindent where $\gamma\in L_{\d}^t(G)$. Also we will prove that these brackets are functorial in $G$ with respect to group scheme homomorphisms. Finally if one defines the ``exponential" maps between sets of points
\begin{equation}
\label{exex}
ex^r:L_{\d}^r(G)\ra G,\ \ ex^r(\alpha)=1+p^r\phi^{-r}(\alpha),\end{equation}
and $ex:=ex^1$
then these maps  are group homomorphisms and we have the following obvious formulae
\begin{equation}
\label{cald}
ex^{r+s}([\alpha,\beta]_{\d})=[ex^r(\alpha),ex^s(\beta)],\ \ \alpha\in L_{\d}^r(G),\ \  \beta\in L_{\d}^s(G),\end{equation}
\begin{equation}
\label{rece}
ex^2([\alpha,\beta]_{\d})=ex((ex(\alpha)\star_{\d} \beta)-_{\d}\beta),\ \ \alpha\in L_{\d}(G),\ \  \beta\in L_{\d}(G),\end{equation}
where $[\ ,\ ]:G\times G\ra G$ is the usual commutator $[a,b]=aba^{-1}b^{-1}$. Of course, the maps $ex^r$ are {\it not} $\d$-maps. On the other hand the maps $ex^r$ are obviously injective and the equation \ref{cald} uniquely determines the brackets $[\ ,\ ]_{\d}$ for all $r,s$. Also equation \ref{rece} uniquely determines $[\ ,\ ]_{\d}$ on $L_{\d}(G)$.

\subsection{Arithmetic analogue of Maurer-Cartan connection}
In what follows, for a matrix $a=(a_{ij})$ with entries in $R$ we set $a^{(p)}=(a_{ij}^p)$, $\phi(a)=(\phi(a_{ij}))$, $\d a=(\d a_{ij})$; so $\phi(a)=a^{(p)}+\d a$.  
Consider the ring $\cO(GL_n)\h=R[x,\det(x)^{-1}]\h$ where $x=(x_{ij})$ is a matrix of indeterminates, and 
$\h$ denotes, as before, $p$-adic completion. 
Assume also that one is given a ring endomorphism $\phi_{GL_n}$ of $\cO(GL_n)\h$ lifting Frobenius, i.e. a ring endomorphism whose reduction mod $p$ is the $p$-power Frobenius on 
$\cO(GL_n)\h/(p)=k[x,\det(x)^{-1}]$; we still denote by $\phi_{GL_n}:\widehat{GL_n}\ra \widehat{GL_n}$ the induced morphism of $p$-formal schemes and we refer to it as a lift of Frobenius on $\widehat{GL_n}$.
Consider the matrices $\Phi(x)=(\phi_{GL_n}(x_{ij}))$ and $x^{(p)}=(x_{ij}^p)$
with entries in $\cO(GL_n)\h$; then
 $\Phi(x)=x^{(p)}+p\Delta(x)$ where $\Delta(x)$ is some matrix with entries in $\cO(GL_n)\h$ which, due to an analogy to be discussed later,  we refer to as the {\it Christoffel symbol}. Note that $\phi_{GL_n}:\widehat{GL_n}\ra \widehat{GL_n}$ is not a morphism over $R$; nevertheless it induces naturally a map (still denoted by $\phi_{GL_n}:GL_n=GL_n(R)\ra GL_n=GL_n(R)$) between the corresponding sets of points, via the formula $\phi_{GL_n}(a)=\phi^{-1}(\Phi(a))$, $a\in GL_n=GL_n(R)$. (The latter map is, of course, not a $\d$-map.) 
Furthermore given a lift of Frobenius $\phi_{GL_n}$ as above one defines an operator
$$\d_{GL_n}:\cO(GL_n)\h\ra \cO(GL_n)\h,\ \ \d_{GL_n}(f)=\frac{\phi_{GL_n}(f)-f^p}{p},$$
referred to as the {\it $p$-derivation} associated to $\phi_{GL_n}$.

 A smooth closed subgroup scheme $G\subset GL_n$ will be called {\it $\phi_{GL_n}$-horizontal} if $\phi_{GL_n}$ sends the ideal of $G$ into itself; in this case we have a lift of Frobenius endomorphism $\phi_G$ on $\widehat{G}$, equivalently on $\cO(G)\h$. 
Furthermore $\d_{GL_n}$ induces an operator $\d_G:\cO(G)\h\ra \cO(G)\h$ referred to as the {\it $p$-derivation} associated to $\phi_G$.

For a fixed lift of Frobenius $\phi_{GL_n}$ on $\widehat{GL_n}$ (and hence for fixed data $\d_{GL_n}$, $\Phi$, $\Delta$) 
we define the {\it arithmetic logarithmic derivative} $l\d:GL_n\ra {\mathfrak g}{\mathfrak l}_n$  by 
\begin{equation}
\label{lda}
l\d a:=\frac{1}{p}\left(\phi(a)\Phi(a)^{-1}-1\right)=(\d a-\Delta(a))(a^{(p)}+p\Delta(a))^{-1}.\end{equation}
This is a $\d$-map of order $1$ and is viewed as an arithmetic analogue of the Maurer-Cartan connection (or of the Kolchin logarithmic derivative).
Note that $l\d$ does not satisfy the cocycle condition but rather a condition which we shall call the {\it skew cocycle condition}; cf. the body of the paper. In particular $l\d$ will satisfy the cocycle condition ``modulo a term of order $0$". 
For $G$ a $\phi_{GL_n}$-horizontal subgroup $l\d$ above induces a $\d$-map of order $1$, 
$l\d:G\ra L_{\d}(G).$
Note  that if $ex:L_{\d}(G)\ra G$ is the ``exponential" defined in \ref{exex} and $\phi_G:G\ra G$ is the map on points induced  by $\phi_{GL_n}:GL_n\ra GL_n$ (i.e. $\phi_G(a)=\phi^{-1}(\Phi(a))$) then we have the following obvious formula in $G$ for $a\in G$:
\begin{equation}
ex(l\d a)=a\cdot \phi_G(a)^{-1}.\end{equation}

Now  any $\alpha\in L_{\d}(G)$  defines  an equation of the form
\begin{equation}
\label{truee}
l\d u=\alpha,
\end{equation}
with unknown $u\in G$; such an equation will be referred to as a {\it $\d_G$-linear} equation or  a {\it $\Delta$-linear} equation (or simply a {\it $\d$-linear} equation if $\Delta$, and hence $\d_G$, has been fixed and is clear from the context).  
Setting $$\Delta^{\alpha}(x) =  \alpha \cdot \Phi(x) +\Delta(x),$$
 $$\Phi^{\alpha}(x)  =  x^{(p)}+p\Delta^{\alpha}(x)  =  \epsilon \cdot \Phi(x),$$
  where $\epsilon=1+p\alpha$, we see that Equation \ref{truee} 
 is equivalent to the  equation:
\begin{equation}\label{clack} 
\d u = \Delta^{\alpha}(u)\end{equation}
and also to the equation
\begin{equation}
\label{clack2}
\phi(u)  =   \Phi^{\alpha}(u).
\end{equation}
 Note that  equation \ref{clack2} is {\it not} a  linear difference equation in the sense of \cite{SVdP}; indeed linear difference equations for $\phi$ have the form $\phi(u)=\alpha \cdot u$. Equations of the form  \ref{truee} will be  studied  in some detail in  \cite{adel3}. Note that for any $\alpha$ as above one can attach a lift of Frobenius $\phi^{\alpha}_{GL_n}$ on $\widehat{GL_n}$ by the formula $\phi^{\alpha}_{GL_n}(x)=\Phi^{\alpha}(x)$; then its associated $p$-derivation $\d^{\alpha}_{GL_n}$ satisfies 
 \begin{equation}
 \label{deltaalpha}
 \d^{\alpha}_{GL_n}(x)=\Delta^{\alpha}(x).\end{equation}

\subsection{Compatibilities with translations and involutions}
In order to implement the formalism above and obtain an interesting theory  we need to construct lifts of Frobenius that satisfy  certain compatibilities with standard group theoretic concepts. We start by discussing  compatibility with translation action of subgroups.
  Let $G$ be, as before, a smooth closed subgroup scheme of $GL_n$; let $H$  be a smooth closed subgroup scheme of $G$. We say that $\phi_G$ is 
 {\it left  compatible} with $H$ if $H$ is $\phi_{GL_n}$-horizontal and
  the following diagram is commutative:
 \begin{equation}
 \label{vigil1} 
  \begin{array}{rcl}
\widehat{H}\times \widehat{G} & \ra & \widehat{G}\\
\phi_H\times \phi_G \downarrow & \  & \downarrow \phi_G\\
\widehat{H}\times \widehat{G} & \ra & \widehat{G}\end{array}\end{equation}
where $\phi_H:\widehat{H}\ra \widehat{H}$ is induced by $\phi_G$ and the horizontal maps are given by multiplication.
Similarly we say that $\phi_G$ is 
 {\it right  compatible} with $H$ if $H$ is $\phi_{GL_n}$-horizontal and
  the following diagram is commutative:
\begin{equation}
\label{vigil2}
\begin{array}{rcl}
\widehat{G}\times \widehat{H} & \ra & \widehat{G}\\
\phi_G\times \phi_H \downarrow & \  & \downarrow \phi_G\\
\widehat{G}\times \widehat{H} & \ra & \widehat{G}\end{array}
\end{equation}
 We say $\phi_G$ is {\it  bicompatible} with $H$ if it is both left and right compatible with $H$.

  Next we discuss compatibility of lifts of Frobenius with involutions or more generally with automorphisms.
  By an automorphism of $G$ we understand  an automorphism  $\tau:G\ra G$ of $G$ viewed as a group scheme over $R$; so $(xy)^{\tau}=x^{\tau}y^{\tau}$ on points. 
  We say $\tau$ is an {\it involution} if $x^{\tau \tau}=x$ on points. (Note that we allow the identity to be an involution; 
  classically, in the Cartan theory of symmetric spaces, one  rules out the case when $\tau$ is the identity but it is convenient here to consider this case as well.)
  Let $\tau$ be an automorphism (not necessarily an involution).
  One can attach to $\tau$ a closed subgroup scheme $G^+$ defined as   the fixed closed subscheme of the automorphism $\tau$; on points,
  \begin{equation}
  \label{Gtaujos}
  G^+=\{x\in G;x^{\tau}=x\}.\end{equation}
 The identity component $(G^+)^{\circ}$ will be called the group defined by $\tau$. Similarly one defines the closed subscheme $G^-$ by
  \begin{equation}
  \label{Gtausus}
  G^-=\{x\in G;x^{\tau}=x^{-1}\}.\end{equation}
  note that $G^-$ is not a subgroup of $G$ in general; but it is stable under $x\mapsto x^{\nu}$ for all $\nu\in \bZ$. For any $x\in G$ write $x^{-\tau}:=(x^{-1})^{\tau}=(x^{\tau})^{-1}$.
 One can attach to $\tau$  a morphism of schemes, $\cH$, which we shall refer to as the {\it quadratic map} attached to $\tau$,
     \begin{equation}
    \label{hash}
\cH:G\ra G,\ \ \     \cH(x)=x^{-\tau}x.\end{equation}
More generally, for any $g\in G$, we may consider the left translation $L_g:G\ra G$ by $g$, $L_g(x)=gx$, and the composition $\cH_g:=L_g\circ \cH$ so
\begin{equation}
\label{gilmore}
\cH_g:G\ra G,\ \ \ \cH_g(x)=gx^{-\tau}x.
\end{equation}
In particular, for $g=1$, we have $\cH_1=\cH$.
  One can also define the morphism of schemes, ${\mathcal B}$, which we shall refer to as the {\it bilinear map} attached to $\tau$,
 \begin{equation}
 \label{scorpion4}
 {\mathcal B}:G\times G\ra G, \ \  {\mathcal B}(x,y)=x^{-\tau}y.
 \end{equation}
 More generally for $g\in G$ we may define 
   a morphism
  \begin{equation}
 \label{scorpion42}
 {\mathcal B}_g:G\times G\ra G, \ \  {\mathcal B}(x,y)=gx^{-\tau}y.
 \end{equation}
  Note that $G^+=\cH^{-1}(1)$. Also, if $\tau$ is an involution then $\cH$ factors through a map $\cH:G\ra G^-$. The pair $(G,G^+)$ is then classically called a {\it symmetric pair} and the orbit space $G^+\backslash G$ is called a {\it symmetric space}.
    
  Assume $\tau$ is any automorphism and let $\phi_G$ and $\phi_{G,0}$ be two  lifts of Frobenius   on $\widehat{G}$.
Then $\phi_{G}$ is said to be {\it $\cH_g$-horizontal} with respect to  $\phi_{G,0}$ if the following diagram is  commutative:
  \begin{equation}
  \label{sleepless} \begin{array}{rcl}
\widehat{G} & \stackrel{\phi_{G}}{\longrightarrow} & \widehat{G}\\
\cH_g  \downarrow & \  & \downarrow \cH_g \\
\widehat{G} & \stackrel{\phi_{G,0}}{\longrightarrow} & \widehat{G}\\
\end{array}\end{equation}
We say $\phi_{G}$ is {\it ${\mathcal B}_g$-symmetric} with respect to $\phi_{G,0}$ if the following diagram is  commutative:
\begin{equation}\label{seatle}
\begin{array}{rcl}
\widehat{G} & \stackrel{\phi_{G,0} \times \phi_{G}}{\longrightarrow} & \widehat{G}\times \widehat{G}\\
\phi_{G} \times \phi_{G,0} \downarrow & \  & \downarrow {\mathcal B}_g\\
 \widehat{G}\times \widehat{G} & \stackrel{{\mathcal B}_g}{\longrightarrow} & \widehat{G}\end{array}\end{equation}
 Clearly ${\mathcal B}_g$-symmetry is equivalent to ${\mathcal B}$-symmetry; also, if $L_g$ commutes with $\phi_{G,0}$ then $\cH_g$-symmetry is equivalent to $\cH$-symmetry.
Note that if $\phi_{G,0}(1)=1$ and if
$\phi_G$ is ${\mathcal H}$-horizontal with respect to $\phi_{G,0}$  then  the group $S:=(G^+)^{\circ}$ (identity component of $G^+$)  is $\phi_{G}$-horizontal and $\phi_G(1)\in S$; we refer to $S$ as the group defined by $\tau$ and  note that, under the assumptions above,  there is an induced lift of Frobenius $\phi_S$ on $\widehat{S}$.
Also note that in the special case when $G=GL_n$ and $\phi_{GL_n,0}(x)=x^{(p)}$, viewing $\cH_g$ as a matrix $\cH_g(x)$ with entries in $R[x,\det(x)^{-1}]\h$, we have that the commutativity of \ref{sleepless} is equivalent to the condition that
$\d_{GL_n}\cH_g=0,$
which can be interpreted as saying that ${\mathcal H}_g$ is a {\it prime integral} for the {\it $\d$-flow} $\d_{GL_n}$.

We will prove the following uniqueness result (cf. Corollary \ref{unique}):  if $\phi_{G,1},\phi_{G,2}$ are two lifts of Frobenius on $\widehat{G}$ which  are both $\cH_g$-horizontal and ${\mathcal B}_g$-symmetric with respect to $\phi_{G,0}$ then $\phi_{G,1}$ and $\phi_{G,2}$   must coincide.

Lifts of Frobenius that are  both $\cH_g$-horizontal and ${\mathcal B}_g$-symmetric with respect to a lift of Frobenius $\phi_{G,0}$ can be viewed as an analogue of  the  {\it hermitian Chern connections} in \cite{GH}, p. 73,  and also as analogues of the {\it Duistermaat connections} in \cite{dui}; cf. the discussion in section \ref{clicking}. 

\subsection{Outer involutions on $GL_n$}
Set $G=GL_n$. We also let  $T\subset GL_n$ be the maximal torus of diagonal matrices, we let $N$ be the normalizer of $T$ in $GL_n$, and we let
 $W\subset GL_n$ be the subgroup of $GL_n$ of all permutation matrices;
 note that $T,W,N$ are smooth group schemes and $N=TW=WT$, in particular the {\it Weyl group} $N/T$ is isomorphic to $W$.
Throughout the paper we will let  $\phi_{GL_n,0}(x)$ be the lift of Frobenius on $\widehat{GL_n}$ defined by $\phi_{GL_n,0}(x):=x^{(p)}$; we will prove that this
 $\phi_{GL_n,0}(x)$
  is the unique lift of Frobenius on $\widehat{GL_n}$ that is bicompatible with $N$  and extends to a lift of Frobenius
 on $\widehat{{\mathfrak g}{\mathfrak l}_n}$ (where we view $\widehat{GL_n}$ as an open set of 
 $\widehat{{\mathfrak g}{\mathfrak l}_n}$). Let $x^t$ denote the transpose of $x$. Let us consider  an involution of $G$ of the form $x^{\tau}=q^{-1}(x^t)^{-1}q$
where $q\in GL_n$, $q^t=\pm q$. 
Then the group $SO(q):=(G^+)^{\circ}$ defined by $\tau$ is smooth.
Also if we let 
$$G_{\pm}=\{a\in G;a^t=\pm a\}$$ 
be the locus of symmetric/antisymmetric matrices in $G$ then 
$$G^-=q^{-1}G_{\pm},$$
so, in particular, since $G_{\pm}$ is smooth (it is an open set in an affine space), $G^-$ is also smooth.
We stress here the fact that $G_{\pm}$  depends only on the sign in $q^t=\pm q$ but {\it not}  on $q$ itself,  whereas $G^+$ and $G^-$  effectively depend on $q$. Note that in the situation above the quadratic map $\cH$ and the bilinear map ${\mathcal B}$ attached to $\tau$ (also referred to as {\it attached to} $q$) are given by 
$$\cH(x)=q^{-1}x^tqx,\ \ \ {\mathcal B}(x,y)=q^{-1}x^tqy.$$
Also a special role will be played by the maps 
$$\cH_q(x)=x^tqx,\ \ {\mathcal B}_q(x,y)=x^tqy$$
which we shall refer to as the {\it normalized}
 quadratic map, respectively the {\it normalized} bilinear map, attached to $q$.
The involutions $\tau$ discussed above will be referred to as {\it outer} involutions. The case of {\it inner} involutions will be discussed later.

It is worth noting that, since $2$ is invertible in $R$ and  $R^{\times}$
is ``closed under the square root operation",   Sylvester's  theorem implies  that for any $q\in GL_n(R)$ with $q^t=q$ there exists $u\in GL_n(R)$ such that $q=u^tu$. In particular for any $q_1,q_2\in GL_n(R)$ with $q_1^t=q_1$ and $q_2^t=q_2$ there exists $u\in GL_n(R)$ such that  $SO(q_2)=u^{-1}\cdot SO(q_2)\cdot u$. Nevertheless, since conjugation  by $u$ 
does not commute in general with the operation $x\mapsto x^{(p)}$ on $GL_n$ the groups $SO(q_1)$ and $SO(q_2)$ behave, in general, quite differently from the viewpoint of the theory we are developing here.

\subsection{Split classical groups}
The  split classical groups $GL_n,SO_n,Sp_n$ are defined by involutions on $G=GL_n$ as follows. We start with
 $GL_n$ itself which is defined by  the identity $x^{\tau}=x$.  We refer to the identity $\tau$ as the canonical involution defining $GL_n$. 
 By a {\it split} matrix we understand a matrix $q\in GL_n$ of one of the following forms
 \begin{equation}
 \label{scorpion3}
 \left(\begin{array}{cl} 0 & 1_r\\-1_r & 0\end{array}\right),\ \ 
\left( 
\begin{array}{ll} 0 & 1_r\\1_r & 0\end{array}\right),\ \ 
\left( \begin{array}{lll} 1 & 0 & 0\\
0 & 0 & 1_r\\
0 & 1_r & 0\end{array}\right),\ \ 
\end{equation}
 where $n=2r, 2r, 2r+1$ respectively, and $1_r$ is the $r\times r$ identity matrix. 
 We denote by $Sp_{2r}, SO_{2r}, SO_{2r+1}$ the groups  defined by the  involutions on $G=GL_n$ given by 
$x^{\tau}=q^{-1}(x^t)^{-1}q$
where $q$ is equal to one of the split matrices \ref{scorpion3} respectively.
 We call this $\tau$ the {\it canonical} involution defining $Sp_{2r}, SO_{2r}, SO_{2r+1}$ respectively. All these groups are smooth over $R$. In our previously introduced notation these groups are, of course,  equal to $(G^+)^{\circ}$.  One also easily checks that $G^-$ is $\phi_{GL_n,0}$-horizontal.
  
\subsection{Special linear group}  For $SL_n$ the picture is slightly different. Indeed, in this case, we need to assume that
 $p\not| n$ and we take
 \begin{equation}
 \label{scorpion1}
 G:=GL_n':=\{\left( \begin{array}{cc}x & 0 \\ 0 & t\end{array}\right) \in GL_{n+1};\ x\in GL_n,\ t\in {\mathbb G}_m,\ \det(x)^2=t^n\};\end{equation}
   note that projection $\pi:G\ra GL_n$ onto the $x$ component is an \'{e}tale homomorphism. In particular any lift of Frobenius on $\widehat{GL_n}$ can be extended uniquely to a lift of Frobenius on $\widehat{G}$. Consider the involution $\tau$ on  $G$ defined by
  \begin{equation}
  \label{scorpion2}\left(\begin{array}{cc} x & 0\\0 & t\end{array}\right)^{\tau}=\left(\begin{array}{ll} t^{-1}x & 0\\0 & t^{-1}\end{array}\right)\end{equation}
 So the group $(G^+)^{\circ}$ defined by $\tau$ projects isomorphically onto $SL_n$ via $\pi:G\ra GL_n$. We will say, by abuse of terminology, that $SL_n$ is defined by $\tau$ and we call $\tau$ the canonical involution defining $SL_n$. 
 Note that 
 $$\cH\left( \begin{array}{cc}x & 0 \\ 0 & t\end{array}\right)=
 \left( \begin{array}{cc}t\cdot 1 & 0 \\ 0 & t^2\end{array}\right);$$
 so $(\pi\circ \cH)^n=\det^2\circ \pi$.
  Furthermore let us equip $\widehat{G}$ with the unique lift of Frobenius $\phi_{G,0}$ extending the lift of Frobenius $\phi_{GL_n,0}$ on $\widehat{GL_n}$.
A lift of Frobenius $\phi_{GL_n}$ on $\widehat{GL_n}$ will be said to be ${\mathcal H}$-horizontal with respect to $\phi_{GL_n,0}$  if the unique  lift of Frobenius $\phi_G$ on $\widehat{G}$ induced by $\phi_{GL_n}$ is ${\mathcal H}$-horizontal with respect to $\phi_{G,0}$; a similar definition is given for ${\mathcal B}$-symmetry. 

Now  if $S$ is any of the split classical groups $GL_n, SL_n, SO_n, Sp_n$
 we define $N_S:=N\cap S$, $T_S:=T\cap S$
 (scheme theoretic intersections). Then $N_S$ and $T_S$ are smooth over $R$, $T_S$ is a maximal torus of $S$, and $N_S$ is the normalizer of $T_S$ in $S$.
 (However the Weyl group $W_S:=N_S/T_S$ is not isomorphic to $W\cap S$ in general!) Finally we recall that  one defines the set of {\it roots} of $S$ as a subset of the characters of the maximal torus $T_S$ of $S$; for any root of $S$ there is a well defined additive {\it root subgroup}
   $U_{\chi}\subset S$ of $S$ and a natural isomorphism ${\mathbb G}_a\simeq U_{\chi}$, ${\mathbb G}_a=Spec\ R[z]\h$, $z$ a primitive element of the Hopf algebra $R[z]$.
   A lift of Frobenius $\phi_{GL_n}$ will be called {\it compatible} with a root $\chi$ of $S$
   if $U_{\chi}$ is $\phi_{GL_n}$-horizontal and the corresponding induced lift of Frobenius on $\widehat{{\mathbb G}_a}$ is given by $z\mapsto z^p$.
   We will need, in the statement below, one more piece of terminology: a root $\chi$ of a classical group $S$ will be called {\it abnormal} if $S=SO_n$, $n$ is odd, and $\chi$ is a ``shortest root" of that group (cf. the body of the paper). 
      
      \subsection{Unitary group}
      Let $q_0\in GL_{2r}$ be the split matrix with $q_0^t=-q_0$ i.e. 
      $$q_0=\left( \begin{array}{cl} 0 & 1_r\\ - 1_r &0\end{array}\right).$$
      Then denote by $GL_r^c:=C(q_0)$ the centralizer of $q_0$ in $GL_{2r}$; it is a  closed (smooth) subgroup scheme of $GL_{2r}$ and its points are exactly the points of $GL_{2r}$ of the form
      \begin{equation}
      \label{z}
      z:=\left(\begin{array}{rr} a & b\\ - b  & a\end{array}\right).\end{equation}
      If we fix a square root $\sqrt{-1}\in R$ of $-1$ we have a natural isomorphism over $R$,
      \begin{equation}
      \label{papers}
      GL_r^c\ra GL_r\times GL_r,\ \ z\mapsto (z^c,(z^t)^c),
      \end{equation}
       where  $z$ is as in \ref{z},  $z^c:=a+\sqrt{- 1}\cdot b$, and hence $(z^t)^c=a-\sqrt{-1}\cdot b$.
      Also $GL_r^c$ is stable under transposition. We denote by $x\mapsto x^*$ the restriction to $GL_r^c$ of the transposition $GL_{2r}\ra GL_{2r}$, $x\ra x^t$  and consider the involution on $GL_r\times GL_r$ given by  $(u,v)^*=(v^t,u^t)$. Then the isomorphism \ref{papers} commutes with $*$.
      Let now $q\in GL_r^c$ be such that $q^*=q$, i.e., 
      $$q=\left(\begin{array}{rr} q_1 & q_2\\ - q_2  & q_1\end{array}\right),\ \ q_1^t=q_1,\ \ q_2^t=- q_2;$$
      we refer to such a $q$ as a {\it hermitian} matrix. 
      We denote by $U^c(q)$ the closed subgroup scheme of $GL_r^c$ defined by the equations $x^*qx=q$; we refer to $U^c(q)$ as the {\it unitary} group of $q$.  
      Note that, for $q=1$, the induced map 
      $$U^c_r:=U^c(q)\ra GL_r,\ \ z\mapsto z^c=a+\sqrt{- 1} \cdot b,$$
with $z$ as in \ref{z},      is an isomorphism; its inverse is given by
$g\mapsto z$ with $z$ as in \ref{z} and 
$$a=\frac{1}{2}(g+(g^t)^{-1}),\ \ b=\frac{1}{2\sqrt{-1}}(g-(g^t)^{-1}).$$ 
In particular $U_r^c$ is connected. Going back to an arbitrary $q$, 
          clearly, $GL_r^c$ is $\phi_{GL_{2r},0}$-horizontal so 
     $\phi_{GL_{2r},0}$ induces a lift of Frobenius $\phi_{GL_r^c,0}$ on $\widehat{GL_r^c}$. On the other hand we have an involution  $x\mapsto q^{-1}(x^*)^{-1}q$ on $GL_r^c$ and the group defined by this involution is the identity component $U^c(q)^{\circ}$. Also we have induced maps ${\mathcal H}_q:\widehat{GL_r^c}\ra \widehat{GL_r^c}$ and ${\mathcal B}_q:\widehat{GL_r^c}\times \widehat{GL_r^c}\ra \widehat{GL_r^c}$. Finally consider the homomorphism 
     $$c:GL_r^c\ra GL_r,\ \ c(z)=z^c$$
      obtained by composing the isomorphism \ref{papers} with the first projection. Let us say that a lift of Frobenius $\phi_{GL_r^c}$ on $\widehat{GL_r^c}$ is $c$-horizontal with respect to a lift of Frobenius $\phi_{GL_r}$ on $\widehat{GL_r}$ if the following diagram is commutative:
     \begin{equation}
     \begin{array}{rcl}
     \widehat{GL_r^c} & \stackrel{\phi_{GL_r^c}}{\longrightarrow} & 
     \widehat{GL_r^c}\\
     c \downarrow & \  & \downarrow c\\
     \widehat{GL_r} & \stackrel{\phi_{GL_r}}{\longrightarrow} & 
     \widehat{GL_r}\end{array}
     \end{equation}
     The terms {\it hermitian} and {\it unitary} are justified by the (obvious) link
     (via Weil restriction) with the classical hermitian and unitary matrices; the letter $c$ was used to invoke {\it complexification}.

   \subsection{Main results on outer involutions}
 With the notation above we have the following result for arbitrary outer involutions; cf. Propositions \ref{existences} and \ref{miramar}.

 \begin{theorem}
 \label{pretzel}
Let $q\in GL_n$ be any matrix with $q^t=\pm q$ and let $\tau$ be the involution on $GL_n$ defined by $x^{\tau}=q^{-1}(x^t)^{-1}q$.
Let  $\cH_q$ and ${\mathcal B}_q$ be the normalized quadratic map
and the normalized  bilinear map attached to $q$. 
Then:

1)  (Lifts on $\widehat{GL_n}$) There is a unique lift of Frobenius $\phi_{GL_n}$ on $\widehat{GL_n}$ that is $\cH_q$-horizontal and ${\mathcal B}_q$-symmetric with respect to $\phi_{GL_n,0}$. 

2) (Lifts on $\widehat{GL_r^c}$) If, in addition,  $n=2r$ and $q^t=q\in GL_r^c$ then $GL_r^c$ is $\phi_{GL_{2r}}$-horizontal; in particular the lift of Frobenius $\phi_{GL_r^c}$ on $\widehat{GL_r^c}$ induced by $\phi_{GL_{2r}}$ is $\cH_q$-horizontal and ${\mathcal B}_q$-symmetric with respect to $\phi_{GL_r^c,0}$ (and is the unique lift of Frobenius with these properties). 

3) (Non-existence of lifts on $\widehat{GL_r}$) Assume the situation in 2) with $q=1=1_n$ the identity, and consider the lift of Frobenius $\phi_{GL_r^c}$ on $\widehat{GL_r^c}$ defined in  2). Then there is no lift of Frobenius $\phi_{GL_r}$ on $\widehat{GL_r}$ such that $\widehat{GL_r^c}$
is $c$-horizontal with respect to $\phi_{GL_r}$.
 \end{theorem}

\begin{remark}
Assertions 1 and 2 in the above theorem can be viewed as an analogue of  results on the existence/uniqueness  for Chern connections in differential geometry. On the other hand assertion 3 should be viewed as contrasting with the situation in differential geometry; cf. section \ref{clicking}.\end{remark}

\begin{remark}
Under  the hypothesis of assertion 3 of Theorem \ref{pretzel}, $U^c(q)$ is $\phi_{GL_r^c}$-horizontal and hence $\phi_{GL_r^c}$ induces on it a lift of Frobenius $\phi_{U^c(q)}$. On the other hand the restriction of $c:GL_r^c\ra GL_r$ to $U^c(q)$ is an isomorphism $c_{U^c(q)}:U^c(q)\ra GL_r$ and hence $\phi_{U^c(q)}$ induces, via $c_{U^c(q)}$, a lift of Frobenius $\phi_{GL_r}$ on $GL_r$; by assertion 3, $\phi_{GL_r^c}$ is, of course {\it not} $c$-horizontal with respect to this $\phi_{GL_r}$.
\end{remark}

\begin{remark}
Let $\phi_{GL_n}$ be the lift of Frobenius in assertion 1 of Theorem \ref{pretzel} and consider the matrix 
$\phi_{GL_n}(x)=\Phi(x)$ with entries in $R[x\det(x)^{-1}]\h$. Then Proposition \ref{existences}  will provide the following formula for the value of $\Phi(x)$ at the identity matrix $1$:
\begin{equation}
\label{valueatone} 
\Phi(1)=(1+p(q^{(p)})^{-1}\d q))^{-1/2}:=1+\sum_{i=1}^{\infty}\left(\begin{array}{c}-1/2\\i\end{array}\right) p^i ((q^{(p)})^{-1}\d q)^i,
\end{equation} 
where $\left(\begin{array}{c}-1/2\\i\end{array}\right)$ are the binomial coefficients. In particular, the value $\Delta(1)$ of the Christoffel symbol at $x=1$ satisfies
\begin{equation}
\label{Christoffel}
\Delta(1)\equiv -\frac{1}{2} (q^{(p)})^{-1}\d q\ \ \ \text{mod}\ \ \ p.
\end{equation}
For $n=1$ and $q\in \bZ$ \ref{valueatone} simplifies to the following formula (cf. \cite{difmod}, Introduction):
\begin{equation}
\label{legendre}
\Phi(1)=q^{(p-1)/2}\cdot \left(\frac{q}{p}\right),
\end{equation}
where $\left(\frac{q}{p}\right)$ denotes the Legendre symbol. (Indeed $q^{(p-1)/2}\cdot \Phi(1)$ has square $1$ and is $\equiv q^{(p-1)/2}$ mod $p$.) This makes \ref{valueatone}, ``up to a polynomial function of the entries of in $q$", a matrix analogue of the Legendre symbol and hence our Christoffel symbol involves a matrix analogue of the Legendre symbol..
\end{remark}

 For the involutions defining the split classical groups we have the following result; 
  cf. Propositions \ref{hepatitis}, \ref{night}, \ref{existences}.

   \begin{theorem}\label{laugh}
 Let $S$ be any  of the groups $GL_n, SL_n, SO_n, Sp_n$, let
 $\tau$ be the canonical involution defining $S$, and let ${\mathcal H}$ and ${\mathcal B}$ be the quadratic and bilinear maps attached to $\tau$.
 Then the following hold.
 
 1) (Compatibility with outer involutions). There exists a  unique lift of Frobenius $\phi_{GL_n}$ on $\widehat{GL_n}$ that is ${\mathcal H}$-horizontal and ${\mathcal B}$-symmetric with respect to $\phi_{GL_n,0}$.
  
  2) (Compatibility with normalizer of maximal torus). The lift of Frobenius  $\phi_{GL_n}$ in assertion 1   is right compatible with $N$  and  also   left compatible (and hence bicompatible) with $N_S$.

  3) (Compatibility  with roots).  The lift of Frobenius
  $\phi_{GL_n}$ in assertion 1 is compatible with a root $\chi$ of $S$ if and only if $\chi$ is not abnormal.
  \end{theorem}

 \begin{remark}
 \label{integrals}
 i) The ${\mathcal H}$-horizontality in assertion 1 of  Theorem \ref{laugh} follows directly from Theorem \ref{pretzel}. It will be supplemented by the following.
 Assume $S=GL_n, SO_n, Sp_n$. Then 
for any $\alpha\in L_{\d}(S)$,  we have $\d_{GL_n}^{\alpha}(\cH)=0,$
where $\d_{GL_n}^{\alpha}$
is as in \ref{deltaalpha};
equivalently $\phi_{GL_n}^{\alpha}$ is ${\mathcal H}$-horizontal with respect to $\phi_{GL_n,0}$.
 The above will be seen to imply  that for any $\alpha\in L_{\d}(S)$, the components of $\cH$ are {\it prime integrals} of the equation $\d u=\Delta^{\alpha}(u)$ in the sense that for any solution $u$ of this equation we have 
 $\d({\mathcal H}(u))=0$. 
 In case $S=SL_n$ we have similar statements. Indeed set ${\mathcal H}^*(x)=\det(x)$. Then for any $\alpha\in L_{\d}(S)$,  we have $\d_{GL_n}^{\alpha}(\cH^*)=0$. Moreover, for any $\alpha\in L_{\d}(S)$, $\cH^*$ is a  {\it prime integral} of the equation $\d u=\Delta^{\alpha}(u)$ in the sense that for any solution $u$ of this equation we have 
 $\d({\mathcal H}^*(u))=0$. 
For more on this see
 the discussion around Equation \ref{primeintegral} and Remark \ref{zgomott} below.  
 
ii) The ${\mathcal B}$-symmetry in assertion 1 of  Theorem \ref{laugh} will   imply the following. 
 Assume we are in case $S=SO_n,Sp_n$. For $b\in {\mathfrak g}{\mathfrak l}_n$ set $b^{\tau}=-_{\d}(q^{-1}b^tq)$ and let
 ${\mathfrak g}{\mathfrak l}_n^{+}$ (respectively ${\mathfrak g}{\mathfrak l}_n^{-}$) be the set of all $b\in {\mathfrak g}{\mathfrak l}_n$ such that $b^{\tau}=b$ (respectively $b^{\tau}=-_{\d}b$). It is easy to show that for any element $b\in  {\mathfrak g}{\mathfrak l}_n$ one has a unique decomposition
 $b=b^++_{\d}b^-$, with $b^+\in {\mathfrak g}{\mathfrak l}_n^{+}$ and 
 $b^-\in {\mathfrak g}{\mathfrak l}_n^{-}$ which will be referred to as the {\it Cartan decomposition} of $b$ with respect to $\tau$.
 Let $l\d:GL_n\ra {\mathfrak g}{\mathfrak l}_n$ be the arithmetic logarithmic derivative associated to $\phi_{GL_n}$ and let $l\d_0:GL_n\ra {\mathfrak g}{\mathfrak l}_n$ be the arithmetic logarithmic derivative associated to $\phi_{GL_n,0}(x)=x^{(p)}$, i.e. $l\d_0 (a)=\d a \cdot (a^{(p)})^{-1}$. Then
 it will follow that for any $a\in G$ one has
 $$l\d (a)\in l\d_0(a)+_{\d} {\mathfrak g}{\mathfrak l}_n^{-}.$$
 In particular if $a\in S$ and 
 $l\d_0 (a)=(l\d_0(a))^++_{\d}(l\d_0(a))^-$
 is the Cartan decomposition of $l\d_0(a)$ then $$l\d(a)=(l\d_0(a))^+.$$
 The latter condition pins down completely the restriction $\phi_S$ of $\phi_{GL_n}$ to $S$.
 
iii) Assertion 2
 of  Theorem \ref{laugh} is reminiscent of conditions satisfied by Cartan connections in principal bundles encountered in classical differential geometry; it is the starting point in \cite{adel3} of the analysis of symmetries of $\d$-linear equations. Indeed
  if $l\d:GL_n\ra {\mathfrak g}{\mathfrak l}_n$ is the arithmetic logarithmic derivative associated to $\phi_{GL_n}$ then assertion 2 above will imply that for all $a\in  N_S$  and $b\in GL_n$ (alternatively for all $a\in GL_n$ and $b\in N$) we have
  $$
\label{cris}l\d(ab)=a \star_{\d} l\d (b) +_{\d} l\d(a).$$

 iv) Assertion 3 of Theorem \ref{laugh} will actually follow from a more general result saying that $\phi_{GL_n}$ coincides with $\phi_{GL_n,0}$ on the set of $R$-points $S\cap \phi_{GL_n,0}^{-1}(S)$. 
 \end{remark}

\begin{remark}
 Theorem \ref{laugh} for $SO_n,Sp_n$ deals with the involutions of $GL_n$ of the form 
 $x^{\tau}=q^{-1}(x^t)^{-1}q$; these can be referred to as {\it outer} involutions.
 The involution defining $SL_n$ can also be referred to as an {\it outer} involution. On the other hand we can consider {\it inner} involutions, or more generally inner automorphisms  $x^{\tau}=q^{-1}xq$, where $q\in GL_n$ is fixed, and ask for their place in our theory. The trivial case $q=1$ (or, more generally, $q$ scalar) is again, covered by Theorem \ref{laugh}. But for  $q$  non-scalar the resulting picture is rather different from that in Theorem \ref{laugh}, as we shall explain in what follows.
 \end{remark}
 
 \subsection{Compatibilities with conjugation} Let $G$ be a smooth affine group scheme over $R$, $H$ a closed smooth subscheme, and  $G^*\subset G$ an open set which is invariant under the action of $G$ on $G$ by conjugation; also let $H^*=H\cap G^*$ and let 
  ${\mathcal C}:G\times G\ra G$ be the conjugation map ${\mathcal C}(a,b)=b^{-1}\star a:=b^{-1}ab$, inducing a map ${\mathcal C}:H^*\times G\ra G^*$. Let $\phi_{G,0}$ be a lift of Frobenius on $\widehat{G}$ and $\phi_{G^*}$ a lift of Frobenius on $\widehat{G^*}$. We say that $\phi_{G,0}$ is {\it ${\mathcal C}$-horizontal} with respect to $\phi_{G^*}$
  if $H$ is $\phi_{G,0}$-horizontal and, upon denoting by $\phi_{H^*,0}$ the  lift of Frobenius on $\widehat{H^*}$ induced by $\phi_{G,0}$, the following diagram is commutative:
  \begin{equation}
  \label{love}
  \begin{array}{rcl}
  \widehat{H^*}\times \widehat{G} & \stackrel{\phi_{H^*,0}\times \phi_{G,0}}{\longrightarrow} & \widehat{H^*}\times \widehat{G}\\
  {\mathcal C}\downarrow & \  & \downarrow {\mathcal C}\\
  \widehat{G^*} & \stackrel{\phi_{G^*}}{\longrightarrow} & \widehat{G^*}\end{array}
  \end{equation}
  
 Going back to the case $G=GL_n$ we take, in the above discussion,  $H=T$, the diagonal torus. Consider the characteristic polynomial
 \begin{equation}
 \label{characteristic}
 \det(s\cdot 1_n-x)=\sum_{i=0}^n(-1)^i{\mathcal P}_i(x)s^{n-i}\in R[x][s]
 \end{equation}
 and its discriminant $D^*(x)\in R[x]$.
 We let $G^*=GL_n^*$ be the open set of {\it regular matrices}, i.e., matrices for which $D^*(x)$ is invertible; so $T^*$ consists of the regular diagonal matrices, i.e. diagonal matrices with diagonal entries distinct mod $p$. Moreover we consider the lift of Frobenius $\phi_{G,0}(x)=x^{(p)}$ on $\widehat{G}$. 
Finally 
 consider the affine space ${\mathbb A}^n=Spec\ R[z_1,...,z_n]$, the {\it characteristic polynomial map} 
 \begin{equation}
 \label{mapp}
{\mathcal P}:G\ra {\mathbb A}^n,\ \ {\mathcal P}(a)=({\mathcal P}_1(a),...,{\mathcal P}_n(a)),\end{equation}
 and the lift of Frobenius $\phi_{{\mathbb A}^n,0}$ on $\widehat{{\mathbb A}^n}$ defined by $\phi_{{\mathbb A}^n,0}(z_i)=z_i^{p}$. We say that a lift of Frobenius $\phi_{G^{**}}$ on an open set $\widehat{G^{**}}$ of $G$ is {\it ${\mathcal P}$-horizontal} with respect to $\phi_{{\mathbb A}^n,0}$ if the following diagram is commutative:
  \begin{equation}
  \label{loving}
  \begin{array}{rcl}
  \widehat{G^{**}} & \stackrel{\phi_{G^{**}}}{\longrightarrow} & \widehat{G^{**}}\\
  {\mathcal P}\downarrow & \  & \downarrow {\mathcal P}\\
  \widehat{{\mathbb A}^n} & \stackrel{\phi_{{\mathbb A}^n,0}}{\longrightarrow} & \widehat{{\mathbb A}^n}\end{array}
  \end{equation}
  On the other hand we will say that $\phi_{G^{**}}$ is a {\it $T$-deformation} of $\phi_{G,0}$ if there is a commutative diagram
   \begin{equation}
  \label{hating}
  \begin{array}{ccl}
  \widehat{G^{**}} & \stackrel{\phi_{G^{**},0}\times \phi_{G^{**}}}{\longrightarrow} & \widehat{G^{**}}\times \widehat{G^{**}}\\
 \downarrow & \  & \downarrow {\mathcal Q}\\
  \widehat{T} & \subset & \widehat{G}\end{array}
  \end{equation}
  where ${\mathcal Q}(x,y)=x^{-1}y$ and $\phi_{G^{**},0}$ is induced by $\phi_{G,0}$; the left vertical map in \ref{hating} is then necessarily unique. Assuming, in addition, that $G^{**}$ is invariant under conjugation by $N$ we say that $\phi_{G^{**}}$ is ${\mathcal C}$-compatible with $N$ if the following diagram is commutative:
   \begin{equation}
  \label{drinking}
  \begin{array}{rcl}
  \widehat{G^{**}} \times \widehat{N}& \stackrel{{\mathcal C}}{\longrightarrow} & \widehat{G^{**}}\\
\phi_{G^{**}}\times \phi_{N,0}  \downarrow & \  & \downarrow \phi_{G^{**}}\\
  \widehat{G^{**}} \times \widehat{N} & \stackrel{{\mathcal C}}{\longrightarrow} & \widehat{G^{**}}\end{array}
  \end{equation}
  where $\phi_{N,0}$ is induced by $\phi_{G,0}(x)=x^{(p)}$.
  
 \subsection{Main result on inner automorphisms} With this notation we have the following result; cf. Propositions  \ref{bade}, \ref{douche}, \ref{sick}, \ref{sickness}, \ref{leak}, \ref{crazy}. Let, in the next Theorem,  $G=GL_n$.

  \begin{theorem}
  \label{cavalier} \ 
  
  1) (Compatibility with conjugation) There exists a unique lift of Frobenius $\phi_{G^*}$ on $\widehat{G^*}$ such that $\phi_{G,0}$  is ${\mathcal C}$-horizontal with respect to $\phi_{G^*}$.
  
  2) (Singularity along the discriminant) For $n\geq 3$ the lift of Frobenius $\phi_{G^*}$ 
  in assertion 1 
  does not extend to a lift of Frobenius on the whole of $\widehat{G}$.
  
  3) (Compatibility with characteristic polynomial) There exists an open set $G^{**}$ of $G$  and  a  lift of Frobenius $\phi_{G^{**}}$ on $\widehat{G^{**}}$ such that  $G^{**}$ is invariant under conjugation by $N$, $G^{**}\cap T=T^*$,
  $\phi_{G^{**}}$
   is ${\mathcal P}$-horizontal with respect to $\phi_{{\mathbb A}^n,0}$, and   $\phi_{G^{**}}$ is a $T$-deformation of $\phi_{G,0}$.
   Moreover, for any such  $G^{**}$, $\phi_{G^{**}}$ is unique with the above properties and  is ${\mathcal C}$-compatible with $N$.
  
  Furthermore let $q\in G$ be such that $q^{(p)}=q$ and consider the involution $x^{\tau}=q^{-1}xq$ on $G$. Then the following hold.
 
  4) (Inner involutions) If $q\in T$, $q^2=1$, $q$ non-scalar,
then
there is no lift of Frobenius $\phi_{G}$  on $\widehat{G}$ that is ${\mathcal H}$-horizontal  with respect to $\phi_{G,0}$. 

5) (Regular inner automorphisms) If $q\in T^*$  then there exists a  lift of Frobenius $\phi_{G}$ on $\widehat{G}$ such that $\phi_{G,0}$  is ${\mathcal H}$-horizontal  with respect to $\phi_{G}$. 

6)  (Non-regular inner automorphisms) If $q\in T\backslash T^*$  and $q$ is  non-scalar then there is no  lift of Frobenius $\phi_{G}$ on $\widehat{G}$ such that $\phi_{G,0}$  is ${\mathcal H}$-horizontal  with respect to $\phi_{G}$. 
\end{theorem}

\begin{remark}
\label{seculca} 
i) We stress the fact that in Theorem \ref{laugh} as well as in assertion 4 of Theorem \ref{cavalier} the condition under consideration is that
 $\phi_{G}$  be ${\mathcal H}$-horizontal  with respect to $\phi_{G,0}$.
On the contrary,
in assertions 5 and 6 of Theorem \ref{cavalier}  the condition under consideration
is that  $\phi_{G,0}$  be ${\mathcal H}$-horizontal  with respect to $\phi_{G}$; hence the roles of $\phi_{G}$ and $\phi_{G,0}$ have been switched.
So the point of assertion 4 in Theorem \ref{cavalier} is that the paradigm of Theorem \ref{laugh}, which was appropriate for outer involutions, is not appropriate for inner involutions. An appropriate paradigm for regular inner automorphisms (i.e. inner automorphisms defined by regular diagonal matrices) is obtained by switching 
the roles of $\phi_{G}$ and $\phi_{G,0}$ as shown by assertion 5 in Theorem \ref{cavalier}. But for non-regular inner automorphisms  even the above switching of roles between the lifts of Frobenius will not fix the problem; cf. assertion 6 of Theorem \ref{cavalier}. 

ii) The maps ${\mathcal C}$ and ${\mathcal P}$ fit into the following commutative diagram
$$
\begin{array}{rcl}
T\times G & \stackrel{{\mathcal C}}{\longrightarrow} & G\\
pr_1\downarrow & \  & \downarrow{\mathcal P}\\
T & \stackrel{{\mathcal S}}{\longrightarrow} & {\mathbb A}^n
\end{array}
$$
where $pr_1$ is the first projection and the components ${\mathcal S}_1,{\mathcal S}_2,...$ of ${\mathcal S}$ are the fundamental symmetric polynomials  ${\mathcal S}_1(t_1,...,t_n)=\sum_i t_i$, ${\mathcal S}_2(t_1,...,t_n)=\sum_{i<j}t_it_j$, etc.
Note however that the ${\mathcal C}$-horizontality in assertion 1 and the ${\mathcal P}$-horizontality in assertion 3 of Theorem \ref{cavalier} do not imply that 
$\phi_{G,0}$ is ${\mathcal P}\circ{\mathcal C}$-horizontal with respect to $\phi_{{\mathbb A}^n,0}$ (in the obvious sense); indeed $\phi_{G^*}$ and $\phi_{G^{**}}$ do not coincide; and actually the above ${\mathcal P}\circ{\mathcal C}$-horizontality
is trivially seen to fail. It is also a fact 
 that $\phi_{T,0}$ is {\it not} ${\mathcal S}$-horizontal with respect to $\phi_{{\mathbb A}^n,0}$ (in the obvious sense).
 
 iii) Assertion 3 in Theorem \ref{cavalier} will be complemented as follows. 
 Let
 $$\phi_{G^{**}}(x)=:x^{(p)}+p\Delta^{**}(x).$$
 More generally,
 if $\alpha(x)$ is any matrix with coefficients in $\cO(\widehat{G^{**}})$ and
 $\epsilon(x)=1+p\alpha(x)$,
 consider the lift of Frobenius $\phi_{G^{**}}^{(\alpha)}$ on $\widehat{G^{**}}$ defined by
 $$\phi_{G^{**}}^{(\alpha)}(x):=\epsilon(x)\cdot 
 (x^{(p)}+p\Delta^{**}(x))
 \cdot \epsilon(x)^{-1}=:x^{(p)}+p\Delta^{**(\alpha)}(x).$$
 Then 
 $\phi_{G^{**}}^{(\alpha)}$ is ${\mathcal P}$-horizontal with respect to $\phi_{{\mathbb A}^n,0}$. Consequently
 the components ${\mathcal P}_i$ of ${\mathcal P}$ are {\it prime integrals} of the equation 
 \begin{equation}
 \label{picture}
 \d u=\Delta^{**(\alpha)}(u)\end{equation}
  in the sense that for any solution $u$ of this equation and for any $i=1,...,n$, we have $\d({\mathcal P}_i(u))=0$; cf. the discussion around Equation \ref{primeintegral} below.  Note however that, assuming the eigenvalues of $u$ are in $R$,  the conditions $\d({\mathcal P}_i(u))=0$ {\it do not} imply (and are {\it not} implied by) the condition that  $\d$ annihilate the eigenvalues of $u$.  These two conditions are inequivalent (although, in special cases, related) ways to express the idea that the spectrum of $u$ does not ``vary" arithmetically.
  
   By the way the equation \ref{picture}
  is trivially seen to be equivalent to the equation
  \begin{equation}
  \label{movies}
  l\d u=-_{\d}(u\star_{\d} \alpha(u))+_{\d}\alpha(u),\end{equation}
  where $l\d:G^{**}\ra L_{\d}(G)$ is the {\it arithmetic logarithmic derivative} attached to $\phi_{G^{**}}$, defined on $G^{**}$ (rather than on the whole of $G$) by the same formula as the usual arithmetic logarithmic derivative \ref{lda}.

iv) The lift $\phi_{G}$ in assertion 5 of Theorem \ref{cavalier} is not unique.

v) If $q\in T$ has the form 
$q=\text{diag}(q_1\cdot   1_{r_1},...,q_s\cdot 1_{r_s})$, $\sum_{i=1}^s r_i=n$,
with $q_1,...,q_s\in R^{\times}$ distinct
then the fixed group $G^+$ of $x^{\tau}=q^{-1}xq$ in $G=GL_n$ has the form
$$G^+=GL_{r_1}\times...\times GL_{r_s},$$
  diagonally embedded into $G$, while the image of $\cH:G^+\backslash G\ra G$
  is the left translate by $q^{-1}$ of the adjoint orbit of $q$. Note that $G^+$ is always, in this case, $\phi_{GL_n,0}$-horizontal; so the non-existence of $\phi_{GL_n}$ in assertion 6 of Theorem \ref{cavalier} may come as a surprise. 
  
 vi) Assertion 5 in Theorem \ref{cavalier} will follow from assertion 1 and assertion 2 will follow from assertion 6. So ${\mathcal C}$-horizontality and ${\mathcal H}$-horizontality are closely intertwined.
 \end{remark}
  
 \subsection{Comparison with classical differential equations}
   \label{clicking}
It is interesting to see what the theory in this paper corresponds to in the classical case of usual differential equations. So assume ({\it in this subsection  only}!), that we are in the framework of differential algebra \cite{kolchin}; this framework approximates, in the algebraic setting, the situation classically encountered in the theory of Lie and Cartan.  So we assume, in this subsection,  that $R$ is a (commutative, unital) ${\mathbb Q}$-algebra  equipped with a derivation $\d$. For any scheme $X$ over $R$ we denote by $X(R)$ the set of its $R$-points. (The typical example we have in mind is that of the ring $R$ of smooth or analytic, real or complex functions on ${\mathbb R}^n$ or ${\mathbb C}^m$; $\d$ is then the usual partial derivative with respect to one of the coordinates.)

\subsubsection{Maurer-Cartan connection}
The {\it Kolchin logarithmic derivative} 
\begin{equation}
\label{kolld}
l\d:GL_n(R)\ra {\mathfrak g}{\mathfrak l}_n(R)\end{equation}
 is the map $l\d(a)=\d a \cdot a^{-1}$ where if $a=(a_{ij})$ then $\d a:=(\d a_{ij})$. This map is an algebraic incarnation of the Maurer-Cartan connection and our map \ref{lda} is an arithmetic analogue of it. (Note, however, that unlike the Kolchin logarithmic derivative \ref{kolld}, our map \ref{lda} is not intrinsically associated to $GL_n$ but also depends on an extra datum which is a lift of Frobenius on $GL_n$. Pinning down the lift of Frobenius, and hence \ref{lda},  in terms of compatibilities with outer involutions defining the classical groups  was the main purpose of Theorems \ref{pretzel} and \ref{laugh}.) If $G\subset GL_n$
is a smooth subgroup scheme defined by equations with coefficients in the ring of constants $R^{\d}=\{c\in R;\d c=0\}$, then one gets that
$l\d$ induces a map $l\d:G(R)\ra {\mathfrak g}(R)$ where ${\mathfrak g}=L(G)\subset {\mathfrak g}{\mathfrak l}_n$ is the Lie algebra of $G$. For any $\alpha\in {\mathfrak g}{\mathfrak l}_n(R)$ one can consider the {\it linear differential equation} $\d u=\alpha u$ with unknown $u\in GL_n(R)$; this equation is, of course, equivalent to the equation $l\d u=\alpha$ of which \ref{truee} is an arithmetic analogue.

\subsubsection{Connections in principal bundles}
Let us consider now $G=GL_n$ and the ring $\cO(G)=R[x,\det(x)^{-1}]$. Consider the unique derivation
  $\d_{G,0}:\cO(G)\ra \cO(G)$ extending $\d:R\ra R$ such that $\d_{G,0}x=0$. Also, for each 
$\alpha\in {\mathfrak g}{\mathfrak l}_n={\mathfrak g}{\mathfrak l}_n(R)$,
consider the unique derivation
  $\d_G:=\d_{G,0}^{\alpha}:\cO(G)\ra \cO(G)$ extending $\d:R\ra R$ such that  \begin{equation}
  \label{deltaa} \d_{G} x=\alpha x.\end{equation}

  \subsubsection{Compatibility with outer involutions}
  Consider  an involution  $x\mapsto x^{\tau}$ on $G$, let $\cH:\cO(G)\ra \cO(G)$ be the $R$-algebra map induced  by  the map $G\ra G$, $x\mapsto x^{-\tau}x$, 
  consider, for $g\in G$, the $R$-algebra map $\cH_g:\cO(G)\ra \cO(G)$ induced  by  the map $G\ra G$, $x\mapsto gx^{-\tau}x$,
  and let
 ${\mathcal B}_g: \cO(G)\ra \cO(G)\otimes_R \cO(G)$ be the $R$-algebra map defined by the map $G\times G\ra G$, $(x_1,x_2)\mapsto gx_1^{-\tau}x_2$.
 Let us say that $\d_{G}$ is $\cH_g$-horizontal with respect to $\d_{G,0}$ if the following diagram is commutative:
  \begin{equation}
 \label{got}
 \begin{array}{ccc}
 \cO(G) & \stackrel{\d_{G}}{\longleftarrow} & \cO(G)\\
 \cH_g \uparrow &\ &\uparrow \cH_g\\
 \cO(G) & \stackrel{\d_{G,0}}{\longleftarrow} & \cO(G)\end{array}
 \end{equation}
 Similarly let us say that $\d_{G}$ is ${\mathcal B}_g$-symmetric with respect to $\d_{G,0}$ if the following diagram is commutative:
 \begin{equation}
 \label{mail}
 \begin{array}{rcl}
 \cO(G) & \stackrel{\d_{G}\otimes 1+1\otimes \d_{G,0}}{\longleftarrow} & \cO(G)\otimes_R \cO(G)\\
 \d_{G,0}\otimes 1+1\otimes \d_{G} \uparrow & \  & \uparrow {\mathcal B}_g\\
 \cO(G)\otimes_R \cO(G) & \stackrel{{\mathcal B}_g}{\longleftarrow} & \cO(G)\end{array}
 \end{equation}
 These diagrams can be viewed as differential algebraic analogues of the diagrams \ref{sleepless} and \ref{seatle} respectively.

 Theorem \ref{pretzel} should be viewed as an arithmetic analogue of the following
 facts. Let  $x^{\tau}=q^{-1}(x^t)^{-1}q$ where $q\in G(R)=GL_n(R)$, $q^t=\pm q$. Hence
${\mathcal H}(x)=q^{-1}x^tqx$ and $\cH_q(x)=x^tqx$. Then 
 $\d_{G}$ is $\cH_q$-horizontal with respect to $\d_{G,0}$
 if and only if 
\begin{equation}
 \label{apa}
\d q+ \alpha^tq+q\alpha=0.\end{equation}
On the other hand 
$\d_{G}$ is ${\mathcal B}_q$-symmetric with respect to $\d_{G,0}$
 if and only if
 \begin{equation}
 \label{apu}
 \alpha^tq-q\alpha=0.\end{equation}
Note that the system consisting of the equations \ref{apa} and \ref{apu} 
(where $q$ is viewed as given and $\alpha$ is viewed as unknown) 
has a unique solution $\alpha\in {\mathfrak g}{\mathfrak l}_n$ equal to 
\begin{equation}
\label{slime}
\alpha=-\frac{1}{2}q^{-1}\d q.
\end{equation}
In other words there is a unique $\alpha \in {\mathfrak g}{\mathfrak l}_n$ such that 
$\d_G:=\d_{G,0}^{\alpha}$ is $\cH_q$-horizontal and ${\mathcal B}_q$-symmetric with respect to $\d_{G,0}$ and this unique $\alpha$ is given by equation \ref{slime}. Let us refer to $\d_G$ as the {\it Chern connection} attached to $q$; this is an analogue of  hermitian Chern connection on a hermitian bundle on a complex manifold, cf.  \cite{GH}, p. 73,  and also of the Duistermaat connections in \cite{dui}.
The analogy is not an entirely direct one;  see  the discussion in subsection \ref{herm} below.  By the way our formula \ref{Christoffel}  should be viewed as an arithmetic analogue of formula \ref{slime} above.

\subsubsection{Chern connections versus hermitian Chern connections}
\label{herm} 
Assume we have an involution on $R$ i.e., a ring automorphism $R\ra R$, $a\mapsto\overline{a}$, whose square is the identity: $\overline{\overline{a}}=a$; we do not require that this automorphism be different from the identity. Also define the derivation $\overline{\d}$ on $R$ by the formula $\overline{\d}(a):=\overline{\d(\overline{a})}$; we further assume that $\d$ and $\overline{\d}$ commute on $R$.
(The typical example we have in mind is that where $R$ is the ring of complex valued smooth functions on a domain in ${\mathbb C}$ and  $\d=\frac{\partial}{\partial z}$ where $z$ is a complex coordinate on ${\mathbb C}$.)
Let 
$$q^t=q=\left(\begin{array}{rr}q_1 & q_2\\ -q_2 & q_1\end{array}\right) \in GL^c_r(R)\subset GL_{2r}(R)$$
 and consider the attached matrix $q^c:=q_1+\sqrt{-1}\cdot q_2\in GL_r(R)$.
 Let $x=\left(\begin{array}{cc}a & b\\ c & d\end{array}\right)$ be a $2r\times 2r$ matrix of variables and $v,w$ be matrices of $r\times r$ variables.
  Consider now 
 the natural isomorphism 
 \begin{equation}
 \label{foot}
 GL_r^c\ra GL_r\times GL_r\end{equation}
  given by the $R$-algebra map
 from 
 $$\cO(GL_r\times GL_r)=R[v,\det(v)^{-1}]\otimes_R R[w,\det(w)^{-1}]$$
  to 
 \begin{equation}
\label{jose}
\cO(GL_r^c)=R[x,\det(x)^{-1}]/(a-d,b+c)=R[a,b,\det\left(\begin{array}{rr}a & b\\ -b & a\end{array}\right)^{-1}]\end{equation}
sending  $v$ and $w$ into the classes of $x^c:=a+\sqrt{-1}\cdot b$ and $(x^t)^c:=a-\sqrt{-1}\cdot b$ respectively.
 Denote by $\d_{GL_r}$ the derivation on $R[v,\det(v)^{-1}]$   which is $\d$ on $R$ and satisfies $\d_{GL_r} v=\alpha^c v$, where $\alpha:=-\frac{1}{2}q^{-1}\d q$, hence
 \begin{equation}
 \label{maradona}
 \alpha^c=-\frac{1}{2}(q^c)^{-1}(\d q)^c=-\frac{1}{2} (q^c)^{-1}\d (q^c).\end{equation}
 We may refer to $\d_{GL_r}$ as the {\it hermitian Chern connection} on 
 $$GL_r=Spec\ R[v,\det(v)^{-1}]$$
  attached to the matrix $q^c$.
   On the other hand  let
 $\d_{GL_{2r}}$ be the Chern connection attached to $q\in GL_{2r}(R)$.  So if
 $\d_{GL_{2r},0}$ is
  the derivation 
 on $R[x,\det(x)^{-1}]$ that lifts
 $\d$ on $R$ and vanishes on $x$ then $\d_{GL_{2r}}$ is 
 the unique  derivation 
 on $R[x,\det(x)^{-1}]$ that lifts
 $\d$ on $R$
and is ${\mathcal H}_q$-horizontal and ${\mathcal B}_q$-symmetric with respect  $\d_{GL_{2r},0}$.
By \ref{slime},  $\d_{GL_{2r}} x=\alpha x$. Then  the link between the hermitian Chern connection $\d_{GL_r}$ attached to $q^c$ and the  Chern connection $\d_{GL_{2r}}$ attached to $q$ is given as follows. 
The derivation $\d_{GL_{2r}}$ induces a derivation, which we  denote by $\d_{GL^c_{r}}$, on the ring \ref{jose} (which we could refer to as the {\it hermitian Chern connection} on $GL_r^c$ attached to $q$); then we claim that the derivation 
$\d_{GL^c_{r}}$
 sends $R[v,\det(v)^{-1}]$ into itself and 
its restriction to $R[v,\det(v)^{-1}]$ equals $\d_{GL_r}$. In other words if $c:GL_r^c\ra GL_r$ is the composition of the isomorphism \ref{foot} with the second projection and we still denote by $c$ the induced algebra map between the ring of regular functions of the two groups then we claim that the following diagram is commutative:
\begin{equation}
\label{cartofi}
\begin{array}{rcl}
\cO(GL_r^c) & \stackrel{\d_{GL_r^c}}{\longleftarrow} &
\cO(GL_r^c)\\
c \uparrow & \  & \uparrow c\\
\cO(GL_r) & \stackrel{\d_{GL_r}}{\longleftarrow} &
\cO(GL_r)\end{array}
\end{equation}
To check the claim  set $\alpha=\left(\begin{array}{rr} \alpha_1 & \alpha_2\\
-\alpha_2&\alpha_1\end{array}\right)$; then, since $\d_{GL_{2r}}x=\alpha x$,
we have
$$\begin{array}{rcl}
\d_{GL_{2r}}v & = & \d_{GL_{2r}}(a+\sqrt{-1}\cdot b)\\
\  & \  & \  \\
\  & = &
\alpha_1 a-\alpha_2 b+\sqrt{-1}(\alpha_1 b+\alpha_2 a)\\
\  & \  & \  \\
\  & = &
 \alpha^c\cdot v,\end{array}$$
which proves our claim. 

The existence of a derivation $\d_{GL_r}$ making the diagram \ref{cartofi} commute is in stark contrast with assertion 3 of our Theorem \ref{pretzel}. In other words the hermitian Chern connection $\d_{GL_r^c}$ above has a direct arithmetic analogue, $\phi_{GL_r^c}$,  while the hermitian Chern connection $\d_{GL_r}$ does not have a direct arithmetic analogue: the natural candidate $\phi_{GL_r}$ does not exist.  

A fuller picture of the analogy with hermitian geometry is obtained as follows.
Denote by $\overline{\d}_{GL_r}$ the unique derivation on $R[v,\det(v)^{-1}]$ that sends $v$ into $0$ and equals the derivation $\overline{\d}$ on $R$. Then $\d_{GL_{2r,0}}$ and $\overline{\d}_{GL_r}$ are related as follows. The derivation $\d_{GL_{2r},0}$ induces a derivation, which we  denote by $\d_{GL^c_{r},0}$, on the ring \ref{jose}; this induced derivation  sends $R[v,\det(v)^{-1}]$ into itself and if we denote by $\d_{GL_r,0}:R[v,\det(v)^{-1}]\ra R[v,\det(v)^{-1}]$ the restriction of $\d_{GL^c_{r},0}$, we have the  equality
\begin{equation}
\label{cartofiori}
\overline{\d}_{GL_r}=\overline{(\ \ )}\circ \d_{GL_r,0}\circ \overline{(\ \ )},\end{equation}
where $\overline{(\ \ )}$ is the automorphism of $R[v\det(v)^{-1}]$ that fixes $v$ and lifts $a\mapsto \overline{a}$ on $R$.
The pair of derivations $(\d_{GL_r},\overline{\d}_{GL_r})$ on $R[v,\det(v)^{-1}]$  is an analogue of the hermitian ``Chern" connection on a hermitian vector bundle on a complex manifold \cite{GH}, p. 73; cf. also the discussion in subsection \ref{vectorbundle}.
\medskip

From now on we assume, for simplicity that $R$ is an algebraically closed field of characteristic zero; varieties over $R$ will be identified with their sets of $R$-points.

\subsubsection{Connections on vector bundles}\label{vectorbundle}
We include, in what follows,  a short digression on the ``vector bundle" version of the above ``principal bundle" formalism.
Start with an $n$-dimensional vector space $V$ over $R$ and, again a derivation $\d$ on $R$. By a {\it connection} on $V$ we understand here an additive group homomorphism $\nabla:V\ra V$ such that 
$$\nabla(a v)=(\d a)v+a \nabla v$$
for all $v\in V$, $a\in R$. Assume we are given a non-degenerate $R$-bilinear map
$$B:V\times V\ra V$$
which is either symmetric or antisymmetric (which we view as an analogue of a Riemannian metric or a  $2$-form respectively). We say that $\nabla$ is $B$-{\it horizontal} if 
\begin{equation}
 \label{eggroll} 
 \d (B(u,v))=B(\nabla u,v)+B(u,\nabla v)
 \end{equation}
 for all $u,v\in V$. (If this is the case we view $\nabla$ as an analogue of a connection that is compatible with a metric or a $2$-form, respectively.) Moreover, given an $R$-linear linear endomorphism $\Lambda$ of $V$ we say that $\Lambda$ is $B$-{\it symmetric}  if 
 \begin{equation}
 \label{tvtv}
 B(\Lambda u, v) =B(u,\Lambda v)\end{equation}
 for all $u,v\in V$. It is easy to see that for any given  non-degenerate bilinear map $B$ and any connection $\nabla_0$  there is exactly one connection $\nabla$ such that $\nabla$ is $B$-horizontal and $\nabla-\nabla_0$ is $B$-symmetric.
  The uniqueness is clear. The existence is proved exactly in the same way one proves the existence of Chern connections. Indeed, let $V^*$ be the $R$-linear dual of $V$, equipped with the {\it dual connection}
  $$\nabla_0^*:V^*\ra V^*$$
  defined as the unique connection with the property that
  $$\d\langle u^*,v\rangle=\langle \nabla_0^* u^*,v\rangle+\langle u^*,\nabla_0 v\rangle,$$
  for $u^*\in V^*$ and $v\in V$; here $\langle\ ,\ \rangle:V^*\times V\ra R$ is the duality pairing. Let $B^*:V\ra V^*$ be the linear map defined by
  $$\langle B^*(u),v\rangle:=B(u,v).$$
  Then define the connection $\nabla:V\ra V$  by
  \begin{equation}
  \label{theone}
  \nabla:=\frac{1}{2}\nabla_0+\frac{1}{2}\cdot (B^*)^{-1}\circ \nabla_0^*\circ B^*.
  \end{equation}
It  is easy to see that $\nabla$ is $B$-horizontal and $\nabla-\nabla_0$ is $B$-symmetric.
 
 The concepts of $B$-horizontality and $B$-symmetry relate to our concepts of $\cH_q$-horizontality and ${\mathcal B}_q$-symmetry as follows. Start with a matrix $\alpha\in {\mathfrak g}{\mathfrak l}_n$ and consider the $R$-linear space $V=R^n$ whose elements we view as column vectors. Consider the connection
$$\nabla:V\ra V,\ \ \ \nabla v=\d v-\alpha v,$$
for $v\in V$,
where if $v=(v_i)$ has components $v_i$ then $\d v:=(\d v_i)$.  Also for $q\in GL_n$ with $q^t=\pm q$ we consider the non-degenerate bilinear (symmetric, respectively antisymmetric) map
$$B:V\times V\ra R,\ \ \ B(u,v):=u^tqv.$$
  It is trivial then to see that $\nabla$ is $B$-horizontal if and only the equation \ref{apa} holds, hence if and only if $\d_G:=\d_{G,0}^{\alpha}$ is $\cH_q$-horizontal with respect to $\d_{G,0}$.
  This makes the commutativity of the diagram \ref{sleepless} an analogue of the condition for a connection to be compatible with a metric or a  $2$-form, respectively. On the other hand we may consider the connection
$$\nabla_0:V\ra V,\ \ \ \nabla_0 v=\d v;$$
this is the unique connection that kills the standard basis of $R^n$.
  Then $\nabla-\nabla_0$ is $B$-symmetric  if and only if the equation \ref{apu} holds, hence if and only if $\d_G=\d_{G,0}^{\alpha}$ is ${\mathcal B}_q$-symmetric with respect to $\d_{G,0}$.

  Note that, in our present context, there is no analogue of the concept of {\it torsion} in Riemannian geometry (because we only have one distinguished derivation $\d$ on the base). In particular, our $\nabla$ is {\it not} an analogue of the Levi-Civita connection of a Riemannian metric. 
  
  Let us also record the hermitian paradigm. Assume the field $R$ has the form $R_0\otimes_{\bZ}\bZ[\sqrt{-1}]$ where $R_0$ is some field and assume $a\mapsto \overline{a}$ is the $R_0$-automorphism of $R$ sending $\sqrt{-1}\mapsto -\sqrt{-1}$.
  Let
  $V$ is an $n$-dimensional vector space over $R$ and let $\d$ and $\overline{\d}$ be two commuting derivations on $R$ such that $\overline{\d} (\overline{a})=\overline{\d a}$ for all $a\in R$. Also let $H:V\times V\ra R$ be a hermitian form (with respect to the involution on $R$). Define a {\it $\d$-connection} on $V$ to be
  an additive operator $\nabla_{\d}:V\ra V$ such that
  $\nabla_{\d}(a v)=\d a \cdot v+ a\cdot \nabla_{\d} v$ for $v\in V$ and $a\in R$;
  similarly define a {\it $\overline{\d}$-connection} on $V$ to be
  an additive operator $\nabla_{\overline{\d}}:V\ra V$ such that
  $\nabla_{\overline{\d}}(a v)=\overline{\d} a \cdot v+ a\cdot \nabla_{\overline{\d}} v$. Define a {\it hermitian connection} to be 
   a pair $\nabla=(\nabla_{\d},\nabla_{\overline{\d}})$
   consisting of a $\d$-connection and a $\overline{\d}$-connection on $V$.
    Say that $\nabla$ is {\it compatible} with $H$ if one of the following two equivalent conditions is satisfied:
  $$\d H(v,w)=H(\nabla_{\d} v,w)+H(v,\nabla_{\overline{\d}}w),\ \ v,w\in V$$
  $$\overline{\d} H(v,w)=H(\nabla_{\overline{\d}} v,w)+H(v,\nabla_{\d}w),\ \ v,w\in V.$$
  Assume we are given a hermitian form $H$ and 
  a $\overline{\d}$-connection $\nabla_0$ on $V$ such that the kernel of $\nabla_0$ spans the $R$-vector space $V$.
  Then there exists  a unique hermitian connection $\nabla=(\nabla_{\d},\nabla_{\overline{\d}})$ on $V$  compatible with $H$ and such that $\nabla_{\overline{\d}}=\nabla_0$. 
  Indeed  take an  $R$-basis $(e_i)$ of $V$ killed by $\nabla_0$ and  consider the matrix $h=(h_{ij})$, where $h_{ij}=H(e_i,e_j)$. Then take $\nabla_{\overline{\d}}=\nabla_0$ and  define $\nabla_{\d}$  by $\nabla_{\d}e_i=\sum_j \overline{\beta}_{ji}e_j$, where 
  the matrix $\beta=(\beta_{ij})$ is defined by $\beta=h^{-1}\d h$. Cf. the computation of the hermitian ``Chern" connection in \cite{GH}, p.73. This and the equation \ref{maradona} justifies the terminology of {\it hermitian Chern connection} used in relation to the equation \ref{maradona}.
   
\subsubsection{Prime integrals}
Going back to our general discussion of the involution $x^{\tau}=q^{-1}(x^t)^{-1}q$ on $G=GL_n$, let $SO(q)$ be the identity component of the fixed group $G^+=\{x\in G;x^{\tau}=x\}$ and let ${\mathfrak s}{\mathfrak o}(q)$ be its Lie algebra.
If $\alpha\in {\mathfrak s}{\mathfrak o}(q)$ and $\d_G:=\d_{G,0}^{\alpha}$ then, trivially, 
 \begin{equation}
 \label{low}
 \d_{G}({\mathcal H}(x))=0;\end{equation}
 so 
 ${\mathcal H}$ is a {\it prime integral} of the equation $\d u=\alpha u$ in the sense that for any solution $u\in G$ of this equation we have $\d({\mathcal H}(u))=0$. 
 There is a  similar ``classical"  analogue of the $SL_n$ case of our Theorem \ref{laugh}. In particular, for ${\mathcal H}^*(x):=\det(x)$ and any $\alpha\in {\mathfrak s}{\mathfrak l}_n$, we have
 \begin{equation}
 \label{slow}
 \d_{G}({\mathcal H}^*(x))=0,\end{equation}
 hence
 ${\mathcal H}^*$ is a {\it prime integral} of the equation $\d u=\alpha u$ in the sense that for any solution $u\in G$ of this equation  we have $\d(\det(u))=0$.
 The equality \ref{slow} follows from the general fact that for any ring equipped with a derivation $\d$ and for any invertible matrix $z$ with coefficients in that ring we have
 \begin{equation}
 \label{general}
 \d(\det(z))=\text{tr}(\d z \cdot z^{-1})\det(z).\end{equation}
  The statement i) of 
 Remark \ref{integrals} is an arithmetic analogue of the  statements \ref{low} and \ref{slow}. This suggests that the family of $p$-derivations $\d_{G}^{\alpha}$ in \ref{deltaalpha} should be viewed as an arithmetic analogue of the ``classical" family of derivations  $\d_{G,0}^{\alpha}$ in \ref{deltaa}. The latter can be viewed as {\it linear flows} on $G$; this justifies viewing the $p$-derivations in \ref{deltaalpha} as arithmetic analogues of linear flows. The main difference between the ``classical" and the arithmetic case is that the derivation $\d_{G,0}$ in the classical case does not depend on $\tau$ (it simply maps $x$ into $0$) whereas the $p$-derivations $\d_{G}=\d_{G}^0$ in Theorem \ref{laugh} depend (and they do so in a rather interesting way) on $\tau$.
 
\subsubsection{Compatibility with inner involutions}
 If, instead of the above ``outer" involutions, we consider  ``inner" involutions $x^{\tau}=q^{-1}xq$ with $q^2=1$ then the following hold. First note that ${\mathcal H}(x)=q^{-1}x^{-1}qx$ and the fixed group of $\tau$ is the centralizer of $q$,  $C(q)=\{x\in G;xq=qx\}$. The Lie algebra of the latter is, of course, ${\mathfrak c}(q)=\{\alpha\in {\mathfrak g}{\mathfrak l}_n;\alpha q=q\alpha\}$. Now if $\alpha\in {\mathfrak c}(q)$ then, for $\d_G:=\d_{G,0}^{\alpha}$,
$$ \d_{G}({\mathcal H}(x))=-q^{-1}x^{-1}(\alpha x)x^{-1} qx
+q^{-1}x^{-1}q(\alpha x)=
0$$
 so 
 ${\mathcal H}$ is a {\it prime integral} of the equation $\d u=\alpha u$ in the sense that for any solution $u\in G$ of this equation we have $\d({\mathcal H}(u))=0$. This fact does not seem to have a direct analogue in the arithmetic setting; cf. assertion 4 in Theorem \ref{cavalier}.
 
 In seeking  a classical fact of which  assertion 3  in Theorem \ref{cavalier} and part  iii) in Remark \ref{seculca} are arithmetic analogues  consider again an arbitrary derivation $\d_{G}:\cO(G)\ra \cO(G)$, $\d_{G}x=\Delta(x)$ and consider the diagram
 \begin{equation}
 \label{cough}
\begin{array}{rcl}
\cO(G) & \stackrel{\d_{G}}{\longleftarrow} & \cO(G)\\
P \uparrow & \  & \uparrow P\\
R[z] &\stackrel{\d_0}{\longleftarrow} & R[z]\end{array}
 \end{equation}
  where $R[z]=R[z_1,...,z_n]$, $P(z_i)={\mathcal P}_i(x)$, $\det(s\cdot 1-x)=\sum_{i=0}^n(-1)^i{\mathcal P}_i(x)s^{n-i}$, $\d_0z_i=0$. Then diagram \ref{cough} can be viewed as an analogue of the diagram \ref{loving}. 
  Assertion 3 in Theorem \ref{cavalier} is analogous to the easily checked fact that if $$\Delta(x)=\alpha x$$ with $\alpha=\alpha(x)$ a diagonal matrix with entries in $\cO(G)$ then \ref{cough} is commutative if and only if $\alpha=0$.
  Also part iii) in Remark \ref{seculca} is an analogue of the fact that if
  $$\Delta(x)=[\alpha,x]:=\alpha x-x\alpha$$ where $\alpha=\alpha(x)$ is any matrix with coefficients in $\cO(G)$ then  the diagram \ref{cough} is commutative.
  Indeed to check the latter is equivalent to checking that $\d_{G}({\mathcal P}_i(x))=0$ for all $i$. Consider the ring $\cO(G)[s,\det(s\cdot 1-x)^{-1}]$ and the unique derivation $\d_{G}$ on this ring that extends $\d_{G}:\cO(G)\ra \cO(G)$ and sends $s\mapsto 0$. We want to show that 
  $\d_{G}(\det(s\cdot 1-x))$ vanishes. But
$$\d_{G}(\det(s\cdot 1-x))  =  \text{tr}(\d_{G}(s\cdot 1-x)\cdot (s\cdot 1-x)^{-1})\det(s\cdot 1-x);$$
so it is enough to show that the above trace vanishes. Now
$$
\begin{array}{rcl}
\text{tr}(\d_{G}(s\cdot 1-x)\cdot (s\cdot 1-x)^{-1}) & = & \text{tr}((x\alpha-\alpha x) \cdot (s\cdot 1-x)^{-1})\\
\  & \  & \  \\
\  & = & \text{tr}(x\alpha(s\cdot 1-x)^{-1})-\text{tr}(\alpha x (s\cdot 1-x)^{-1})\\
\  & \  & \  \\
\  & = & \text{tr}(\alpha(s\cdot 1-x)^{-1}x)-\text{tr}(\alpha x (s\cdot 1-x)^{-1})\\
\  & \  & \  \\
\  & = & 0,
\end{array}
$$
because $x$ and $(s\cdot 1-x)^{-1}$ commute, and our claim is proved.
In particular one gets that the polynomials ${\mathcal P}_i(x)$ are {\it prime integrals} of the equation 
\begin{equation}
\label{sun}
\d u=[\alpha(u),u],\end{equation}
equivalently of the equation
\begin{equation}
\label{moon}
l\d u=-(u\star \alpha(u))+\alpha(u),
\end{equation}
in the sense that for any solution $u$ of this equation we have $\d({\mathcal P}_i(u))=0$.
Here $\star$ is the adjoint action, $u\star v:=uvu^{-1}$.
This is analogous to part iii) in Remark  \ref{seculca}; especially Equation \ref{moon}
is analogous to Equation \ref{movies}. Cf. Remark \ref{death} in the body of the paper. Equation \ref{sun} can be viewed as an ``isospectral flow" on the space of matrices; so the lift of Frobenius $\phi_{G^{**}}$ in assertion 3 of Theorem 
\ref{cavalier} can be viewed as an arithmetic analogue of such an ``isospectral flow".
  
\subsection{Organization of the paper}
The  paper is organized as follows. Section 2 reviews some of the basic concepts 
in \cite{char,book} and adds some complements to them.  In section 3 we introduce arithmetic Lie theory in the ``abstract case". In section 4 we specialize the general theory to the case of $GL_n$. In section 5 we consider outer involutions and their classical symmetric spaces and we prove Theorems \ref{pretzel} and  \ref{laugh}. In section 6 we examine inner automorphisms and conjugacy classes and we prove Theorem \ref{cavalier}.

\section{Review of $p$-jets}

This section is devoted to reviewing some of the material in \cite{char, book, BYM}.

Throughout the paper we fix an odd  prime $p\in \bZ$. Unless otherwise specified rings will always be assumed commutative with unit. For any ring $A$ we denote by $\widehat{A}$ the $p$-adic completion of $A$. Similarly for any scheme $X$ of finite type over $R$ we denote by $\widehat{X}$ the $p$-adic completion of $X$. A $p$-formal scheme (respectively a $p$-formal scheme of finite type over a base ring) is a formal scheme locally isomorphic to the $p$-adic completion of a Noetherian scheme (respectively of a scheme of finite type over the base ring). If $u:A\ra B$ is a ring homomorphism we usually still denote by $u:Spec\ B\ra Spec\ A$ the induced map; and conversely if $f:X\ra Y$ is a morphism of schemes we still denote by $f:\cO_Y\ra f_*\cO_X$ the induced morphism on functions.

A {\it $p$-derivation} on a ring $A$ is a map of sets $\d:A\ra A$ such that $\d(1)=0$ and, for all $a,b\in A$,
$$\begin{array}{rcl}
\d(a+b) &= &\d(a)+\d(b)+C_p(a,b)\\
\d(ab)&=&a^p\d(b)+b^p\d(a)+p\d(a)\d(b),\end{array}$$
where $C_p(x,y)\in \bZ[x,y]$ is the polynomial $C_p(x,y)=p^{-1}(x^p+y^p-(x+y)^p)$.

A {\it lift of Frobenius} on a ring $A$ will mean a ring homomorphism $\phi=\phi_A:A\ra A$ whose reduction mod $p$ is the $p$-power Frobenius $A/pA\ra A/pA$. 
If $\d:A\ra A$ is a $p$-derivation then $\phi:A\ra A$ defined by $\phi(a)=a^p+p\d a$ is a lift of Frobenius. Conversely
if $A$ is $p$-torsion free and $\phi:A\ra A$ is a lift of Frobenius then the operator $$\d=\d_A:A\ra A,\ \ \d a=\frac{\phi(a)-a^p}{p}$$ is a  $p$-derivation referred to as the $p$-derivation attached to $\phi$. 

A {\it lift of Frobenius} on a scheme (or $p$-formal scheme) $X$ will mean a morphism of schemes (or $p$-formal schemes respectively)  $\phi=\phi_X:X\ra X$ whose reduction mod $p$ is the $p$-power Frobenius;  if $X$ is a $p$-formal scheme and $\cO_X$ is $p$-torsion free the above construction globalizes to yield an operator
 $\d=\d_X:\cO_X \ra \cO_X$ referred to as the $p$-{\it derivation} attached to $\phi_X$. If $f\in \cO(X)$ is a global function then we say $f$ is {\it $\d_X$-constant}
 if $\d_X f=0$. 
  If $X$ is a scheme and $\phi_{\widehat{X}}$ is a lift of Frobenius on $\widehat{X}$ we usually denote $\phi_{\widehat{X}}$ and $\d_{\widehat{X}}$ simply by $\phi_X$ and $\d_X$ respectively.
  
If $\pi:X \ra Y$ is a morphism of schemes (or $p$-formal schemes) and $X,Y$ are given lifts of Frobenius $\phi_X,\phi_Y$ we say that $\pi$ is {\it horizontal} with respect to $\phi_X,\phi_Y$ if the following diagram commutes
$$\begin{array}{rcl}
X & \stackrel{\phi_X}{\longrightarrow} & X\\
\pi \downarrow & \  & \downarrow \pi\\
Y & \stackrel{\phi_Y}{\longrightarrow} & Y\end{array}$$
We also say that $\phi_X$ is {\it $\pi$-horizontal} with respect to $\phi_Y$.
If $Z \subset X$ is a closed subscheme (or closed $p$-formal subscheme) and we are given a lift of Frobenius $\phi_X$ we say $Z$ is {\it $\phi_X$-horizontal} if there exists a (necessarily unique) lift of Frobenius $\phi_Z$ such that $Z\subset X$ is horizontal.
If all the objects above are $p$-torsion free then we have at our disposal attached $p$-derivations and we sometimes say {\it $\d_X$-horizontal} instead of {\it $\phi_X$-horizontal}.

Recall from the Introduction that $R$ denotes, for us, the complete discrete valuation ring with maximal ideal generated by  $p$ and algebraically closed residue field $k={\mathbb F}_p^a$. Then $R$ comes equipped with a unique lift of Frobenius $\phi=\phi_R:R\ra R$. 
If $X$ is a scheme (or a $p$-formal scheme) over $R$ then lifts of Frobenius $\phi_X:X\ra X$
will always be assumed to be such that $X\ra Spec\ R$ is horizontal with respect to $\phi_X,\phi_R$; in particular $\phi_X$ in {\it not} a morphism of $R$-schemes. Nevertheless  there is a naturally induced map between $R$-points $\phi_{X,*}:X(R)\ra X(R)$
that sends a point $P:Spec\ R\ra X$ into the $R$-point 
$$\phi_{X,*}(P):=\phi_X \circ P\circ \phi_R^{-1}:Spec\ R\ra Spec\ R \ra X\ra X;$$ the  map $\phi_{X,*}:X(R)\ra X(R)$ is, of course,  {\it not} regular in general  and it is {\it not} generally a $\d$-map in the sense of \cite{char, book}. Indeed, in case, say, $X={\mathbb A}^1=Spec\ R[x]$, if $\phi_X(x)=\Phi(x)\in R[x]\h$
and $P\in {\mathbb A}^1(R)\simeq R$ corresponds to $x\mapsto a\in R$ then $\phi_{X,*}(P)\in {\mathbb A}^1(R)\simeq R$
corresponds to $x\mapsto \phi_R^{-1}(\Phi(a))\in R$. We will usually simply write $\phi_X$ instead of $\phi_{X,*}$ in what follows.
More generally if $X$, $Y$ are schemes (or $p$-formal schemes) over $R$ then a morphism of schemes over $\bZ$ (or $p$-formal schemes over $\bZ_p$) $f:X\ra Y$
 will be called $\phi^r$-linear if the following diagram is commutative
$$
\begin{array}{ccc}
X & \stackrel{f}{\longrightarrow} & Y\\
\downarrow & \  & \downarrow \\
Spec\ R & \stackrel{\phi^r}{\longrightarrow} & Spec\ R\end{array}
$$
Any such $f$ induces a map of sets $f_*:X(R)\ra Y(R)$, $f_*(P)=f\circ P\circ \phi_R^{-r}$. Again, for simplicity, we will write $f$ instead of $f_*$.

Recall from the Introduction that for any  smooth scheme (or $p$-formal scheme) $X$ over $R$  we usually continue to denote by $X$ the set  $X(R)$ of $R$-points of $X$. If $Y\subset X$ is a closed smooth subscheme (respectively $p$-formal scheme) then it is trivial to check that $Y$ is $\phi_X$-horizontal if and only if for all $R$-points $P\in Y$ one has $\phi_X(P)\in Y$. (The latter is an easy exercise using Nullstellensatz over $k$ and the surjectivity of $Y(R)\ra Y(k)$.) Also two lifts of Frobenius on a smooth $\widehat{X}$ coincide if and only if the corresponding maps on points coincide. We will later need the following:

\begin{lemma}\label{pat}
Let $Z\subset Y \subset X$ be closed embeddings of smooth schemes and let $\phi_{X,0}$ and $\phi_{X,1}$ be two lifts of Frobenius on $\widehat{X}$. Assume that $Z$ is $\phi_{X,0}$-horizontal and that $\phi_{X,0}$ and $\phi_{X,1}$ coincide on the set $Y\cap\phi^{-1}_{X,0}(Y)$. Then $Z$ is also $\phi_{X,1}$-horizontal and the restrictions of $\phi_{X,0}$ and $\phi_{X,1}$ to $Z$ coincide.
\end{lemma}

{\it Proof}.
It is enough to show that for any point $P\in Z$ we have $\phi_{X,1}(P)\in Z$ and $\phi_{X,1}(P)=\phi_{X,0}(P)$. Since $Z$ is $\phi_{X,0}$-horizontal we have $\phi_{X,0}(P)\in Z$. So $P\in Y\cap \phi_{X,0}^{-1}(Y)$. Hence $\phi_{X,1}(P)=\phi_{X,0}(P)\in Z$ and we are done.
\qed

\bigskip

The following is also useful and trivial to check:

\begin{lemma}
\label{liftingg}
Let $U$ be an affine open subscheme of the affine space ${\mathbb A}^n$ over $R$ and let $X$ be a closed $p$-formal subscheme of $\widehat{U}$. Then any lift of Frobenius on $X$ can be prolonged to a lift of Frobenius on $\widehat{U}$.
\end{lemma}

Next we review $p$-jet spaces by  closely following the discussion in \cite{adel1}. 
For $x$ a tuple of indeterminates over $R$ and tuples of indeterminates $x',...,x^{(n)},...$
we consider the ring of polynomials in infinitely many variables $R\{x\}=R[x,x',x'',...]$ and we still denote by $\phi:R\{x\}\ra R\{x\}$  the unique lift of Frobenius
extending $\phi$ on $R$ such that $\phi(x)=x^p+px'$, $\phi(x')=x^p+px''$, etc. Then we consider the $p$-derivation 
 $\d:R\{x\}\ra R\{x\}$  defined by $$\d f=\frac{\phi(f)-f^p}{p}.$$  
 Now
for any affine scheme of finite type 
$$X=Spec\ \frac{R[x]}{(f)}$$ over $R$, where $f$ is a tuple of polynomials, we define the {\it $p$-jet spaces} of $X$ as being the $p$-formal schemes
\begin{equation}
\label{fort}
J^n(X)=Spf\ \frac{R[x,x',...,x^{(n)}]\h}{(f,\d f,...,\d^n f)}.\end{equation}
So $J^0(X)=\widehat{X}$.
For $X$ of finite type but not necessarily affine we define $J^n(X)=\bigcup J^n(X_i)$ where $X=\bigcup X_i$ is an affine cover and the gluing is an obvious one. The spaces $J^n(X)$ have an obvious universality property for which we refer to \cite{char,book} and can be defined for $X$ a $p$-formal scheme  of finite type as well. 
There are natural $R$-morphisms
\begin{equation}
\label{defofpi}\pi:J^n(X)\ra J^{n-1}(X)\end{equation}
induced by the inclusions $R[x,x',...,x^{(n-1)}]\h\subset R[x,x',...,x^{(n)}]\h$ and also
$\phi$-linear morphisms
\begin{equation}
\label{defofpiphi}
\pi_{\phi}:J^n(X)\ra J^{n-1}(X)\end{equation}
induced by the morphisms $\phi:R[x,x',...,x^{(n-1)}]\h\ra R[x,x',...,x^{(n)}]\h$.  
If $X/R$ is smooth then $J^n(X)$ is
a smooth $p$-formal scheme.
The universality property yields natural maps on sets of $R$-points $\nabla^n:X(R)\ra J^n(X)(R)$; for $X={\mathbb A}^m$ (the affine space),
$J^n(X)={\mathbb A}^{m(n+1)}$ and $\nabla^n(a)=(a,\d a,...,\d^na)$.
Let $X,Y$ be schemes of finite type over $R$; by a {\it $\d$-map} of order $n$, $f:X\ra Y$, we understand a map of $p$-formal schemes $J^n(X)\ra J^0(Y)=\widehat{Y}$. Two $\d$-maps 
$X\ra Y$ and $Y\ra Z$ of orders $n$ and $m$ respectively can be composed (using the universality property) to yield a $\d$-map of order $n+m$. Any $\d$-map $f:X\ra Y$ induces a set theoretic map  $f_*:X(R)\ra Y(R)$ defined by $f_*(P)=f(\nabla^n(P))$; if $X,Y$ are smooth the map $f_*$ uniquely determines the map $f$ and, in this case, we simply write $f$ instead of $f_*$ (and $X,Y$ instead of $X(R),Y(R)$).
A $\d$-map $X\ra Y$ of order zero 
is nothing but a map of $p$-formal schemes $\widehat{X}\ra \widehat{Y}$.
The functors $J^n$ commute with products and send groups into groups. By a {\it $\d$-homomorphism}
$f:G\ra H$ between two group schemes (or group $p$-formal schemes) we understand a group homomorphism $J^n(G)\ra J^0(H)=\widehat{H}$. 

We end  by discussing the formalism of flows and prime integrals, following \cite{BYM}. 
Let $X$ be a an affine smooth  scheme over $R$.
A system of {\it arithmetic differential equations} of order $r$ is simply a subset ${\mathcal E}$ of $\cO(J^r(X))$. 
By a {\it solution} of ${\mathcal E}$ we mean an $R$-point $P\in X$ such that $f(P)=0$ for all $f\in {\mathcal E}$. We denote by $Sol({\mathcal E})\subset X$ the set of solutions of ${\mathcal E}$. By a {\it prime integral} of ${\mathcal E}$ we mean a function $\cH\in \cO(X)\h$ such that 
\begin{equation}
\label{primeintegral}
\d(\cH(P))=0,\ \ \text{for all}\ \ P\in Sol({\mathcal E}).\end{equation}
 Intuitively $\cH$ is ``constant along the solutions of ${\mathcal E}$".
 
  Let now $\d_X$ be a $p$-derivation on $\cO(X)\h$. Then one can define a system of arithmetic differential equations of order $1$, ${\mathcal E}(\d_X)\subset \cO(J^1(X))$, which we can refer to as the $\d$-{\it flow} associated to $\d_X$; by definition we take 
 ${\mathcal E}(\d_X)$ to be the ideal in $\cO(J^1(X))$ generated by all elements of the form $\d f-\d_X f$ where $f\in \cO(X)\h$. (Here $\d f$ is the image of $f$ under the ``universal" $p$-derivation $\d:\cO(X)\h \ra \cO(J^1(X))$  while $\d_X f$ is the image of $f$ under the $p$-derivation $\d_X:\cO(X)\h\ra \cO(X)\h\subset \cO(J^1(X))$.) By a $\d$-flow on $X$ we will understand an ideal in $\cO(J^1(X))$  of the form ${\mathcal E}(\d_X)$ where
 $\d_X$ is some  $p$-derivation on $\cO(X)\h$. If we denote by ${\mathcal F}(X)$ the set of $\d$-flows on $X$ then we have natural bijections
 $$\begin{array}{rcl}
 {\mathcal F}(X) & := & \{\text{$\d$-flows on $X$}\}\\
 \  & \simeq & \{\text{$p$-derivations on $\cO(X)\h$}\}\\
 \  & \simeq & \{\text{lifts of Frobenius on $\widehat{X}$}\}\\
 \  &  \simeq & \{\text{sections of $J^1(X)\ra \widehat{X}$}\}\end{array}$$
 If $\d_X$ is a $p$-derivation on $\cO(X)\h$ and $\sigma:\widehat{X}\ra J^1(X)$ is the corresponding section then the ideal ${\mathcal E}(\d_X)$ equals the ideal
 of the image of $\sigma$.
 
 \begin{lemma}\label{olema}
Let $\cO(X)=(R[x_1,...,x_n]/I)_g$ for indeterminates $x_1,...,x_n$ and some $g\in R[x_1,...,x_n]$ and let us continue to denote by  $x_i\in \cO(X)\h$  the images of $x_i$. Then the elements
 \begin{equation}
 \label{elements}
 \d x_i-\d_X x_i,\ \ i=1,...,n\end{equation}
 generate the ideal ${\mathcal E}(\d_X)$. In particular 
 $$Sol({\mathcal E}(\d_X))=Sol(\d x_1-\d_X x_1,...,\d x_n-\d_X x_n).$$\end{lemma}
 
 {\it Proof}.
  Denote by ${\mathcal I}$ the ideal of $\cO(J^1(X))$ generated by the elements \ref{elements} and consider the set
 $${\mathcal A}=\{f\in \cO(X)\h;\d f-\d_X f\in {\mathcal I}\}.$$
 Then ${\mathcal A}$ is easily checked to be a subring of $\cO(X)\h$ and it clearly contains $R$ and $x_i$ for all $i$. It also contains the image of $1/g$ because
 $$\d\left( \frac{1}{g}\right) -\d_X \left( \frac{1}{g} \right)=
-\frac{\d g-\d_X g}{(g^p+p\d g)(g^p+p \d_X g)}.$$
By Noetherianity ${\mathcal I}$ is $p$-adically  closed in $\cO(J^1(X))$ so ${\mathcal A}$ is $p$-adically closed in $\cO(X)$; so it must coincide with $\cO(X)\h$.\qed

\begin{lemma}
\label{incaolema} 
If $\cH\in \cO(X)\h$ is such that $\d_X \cH=0$ then $\cH$ is a prime integral for ${\mathcal E}(\d_X)$.
\end{lemma}

{\it Proof}.
Indeed if $P\in Sol({\mathcal E}(\d_X))$ then 
$$0=(\d \cH-\d_X \cH)(P)=(\d \cH)(P)-(\d_X \cH)(P)=(\d \cH)(P).$$
\qed

\begin{remark}\label{zgomott} 
In the terminology above, any $\Delta$-linear equation $\d u=\Delta^{\alpha}(u)$, cf.  \ref{clack}, may
 be identified with  the  system of arithmetic differential equations on $GL_n$,
\begin{equation}
\label{robinet}
\d x_{ij}-\Delta^{\alpha}_{ij}(x),\ \ i,j=1,...,n,\end{equation}
where $(\Delta^{\alpha}_{ij}(x))=\Delta^{\alpha}(x)$, and hence, by Lemma \ref{olema}, it has the same solution set as the    $\d$-flow ${\mathcal E}(\d^{\alpha}_{GL_n})$ on $GL_n$, where $\d^{\alpha}_{GL_n}(x)=\Delta^{\alpha}(x)$.
Also note that, by Lemma \ref{incaolema}, for any  $\cH\in \cO(GL_n)\h$ with $\d_{GL_n}^{\alpha}\cH=0$ we have that $\cH$ is a prime integral of the system \ref{robinet}; i.e. for all solutions  $u$ of $\d u=\Delta^{\alpha}(u)$ we have $\d({\mathcal H}(u))=0$. \end{remark}

\section{Abstract arithmetic Lie theory}

In this section we introduce 
arithmetic analogues of some basic Lie theoretic concepts in the case of
``abstract" group schemes. The ``concrete" case of classical groups will be analyzed in the next sections.

\begin{definition}
Let $G$ be a smooth group scheme over $R$ (or more generally a group $p$-formal scheme).
 The {\it $\d$-Lie algebra} of $G$  is the group  $p$-formal scheme:
$$L_{\d}(G):=Ker(J^1(G)\ra J^0(G)).$$\end{definition}

\begin{remark}
This construction is obviously functorial in the following sense: if $f:G\ra G'$ is a homomorphism of group schemes (or group $p$-formal schemes) then we have an induced homomorphism of group $p$-formal schemes $L_{\d}(f):L_{\d}(G)\ra L_{\d}(G')$. \end{remark}

\begin{remark}
If $L(G)$ is the Lie algebra of $G$ (viewed as a vector group scheme)  then we have a (non-canonical!) isomorphism
of $p$-formal schemes
$$L_{\d}(G)\simeq \widehat{L(G)};$$
cf. \cite{char}.
This is not an isomorphism of groups in general; indeed $L(G)$ is always commutative whereas $L_{\d}(G)$ is commutative if and only if $G$ is commutative. 
\end{remark}

In what follows we freely use our convention that if $X$ is a scheme (or a $p$-formal scheme) of finite type we denote its set $X(R)$ of $R$-points simply by $X$; in particular we write $a\in X$ instead of $a\in X(R)$. We also identify, as usual, the $\d$-maps $X\ra Y$ with the maps between the sets of $R$-points. Note that if $X$ is a smooth $R$-scheme then $\widehat{X}(R)\simeq X(R)$; so, according to our conventions, both these sets will sometimes be denoted simply by $X$.

In what follows we define the $\d$-{\it adjoint action}. Consider  the morphism  $G\times G \ra G$ which on $R$-points acts as $(g,h)\mapsto ghg^{-1}$, $g,h\in G$. We get an obvious induced morphism of $p$-formal schemes (i.e. $\d$-map of order $0$)
$J^1(G)\times J^1(G)\ra J^1(G)$ which, by restriction, yields a $\d$-map of order $0$ 
\begin{equation}
\label{hominid}
\star:J^1(G)\times L_{\d}(G)\ra L_{\d}(G),\end{equation}
which is, of course, an action that can be referred to as the adjoint action.
 Now the identity $id:J^1(G)\ra J^1(G)$ defines a $\d$-homomorphism of order $1$, 
 $\nabla^1=id_*$,
 $$\nabla^1:G\ra J^1(G).$$
  We get an induced $\d$-map of order $1$, 
  \begin{equation}
  \label{monkey}
  \nabla^1\times id:G\times L_{\d}(G)\ra J^1(G)\times L_{\d}(G).\end{equation}
  
  \begin{definition}
  The  $\d$-map of order $1$
\begin{equation}\label{dino}
\star_{\d}:G \times L_{\d}(G)\ra L_{\d}(G)\end{equation}
   obtained by composing \ref{monkey}  with  \ref{hominid} is called the  $\d$-{\it adjoint action}. \end{definition}
   
   This map is indeed an action. By the way one has a similar construction in the  $\delta$-algebraic setting; in this case the resulting action   has order $0$ and is the usual adjoint action; cf. \cite{aim}.

Here are two  other basic actions. Denote by ${\mathcal F}(G)$ the set of $\d$-flows on $G$; we identify $\d$-flows on $G$ with sections $\sigma:\widehat{G}\ra J^1(G)$ in the category of $p$-formal schemes of the projection $J^1(G)\ra \widehat{G}$. (N.B. $\sigma$ are not group homomorphisms in general!) Then ${\mathcal F}(G)$ is acted upon, via left multiplication, by the group $L_{\d}(G)(G)$ of $p$-formal scheme maps $\widehat{G}\ra L_{\d}(G)$; of course ${\mathcal F}(G)$ is a principal homogeneous space for this action so it is either empty or in bijection with $L_{\d}(G)(G)$. On the other hand the group $L_{\d}(G)=L_{\d}(G)(R)$ acts via left multiplication on ${\mathcal F}(G)$; if $\alpha\in L_{\d}(G)$ and $\sigma\in {\mathcal F}(G)$, $\sigma:\widehat{G}\ra J^1(G)$, then the action is denoted by $(\alpha,\sigma)\mapsto \sigma^{\alpha}$ and is defined by $\sigma^{\alpha}(v)=\alpha \cdot \sigma(v)$. So if one fixes, once and for all, a section $\sigma_G\in {\mathcal F}(G)$ then we have a natural identification 
$$L_{\d}(G)\backslash L_{\d}(G)(G)\ \simeq\ L_{\d}(G)\backslash {\mathcal F}(G).$$
The orbits of $L_{\d}(G)$ on $L_{\d}(G)(G)$ and ${\mathcal F}(G)$ respectively will  be referred to as {\it linear equivalence classes}. There is also a right action of $G=G(R)$ on ${\mathcal F}(G)$ defined as follows: for any $u\in G$ and any $\sigma\in {\mathcal F}(G)$, $\sigma:\widehat{G}\ra J^1(G)$, the action is denoted by $(u,\sigma)\mapsto \sigma^{\star u}$ and  is given by $\sigma^{\star u}(v)=\nabla^1(u)^{-1}\cdot \sigma(uv)$.

Assume now in addition that $G$ is affine. Then, by smoothness, there exists a lift of Frobenius $\phi_G:\cO(G)\h\ra \cO(G)\h$ extending $\phi:R\ra R$ (and not necessarily compatible with the comultiplication) hence a lift of Frobenius $\phi_G$ on $\widehat{G}$; 
one can attach to $\phi_G$ the $p$-derivation $\d_G:\cO(G)\h\ra \cO(G)\h$ defined by $$\d_G(f)=p^{-1}(\phi_G(f)-f^p).$$ 
The latter induces (and is actually  equivalent to giving) a section $\sigma_G:\widehat{G}\ra J^1(G)$ of the projection $J^1(G)\ra \widehat{G}$ in the category of $p$-formal schemes. Fix from now on $\phi_G$ (equivalently $\d_G$ or $\sigma_G$.)
We may then consider 
the $\d$-map of order $0$, 
$$\{\ ,\ \}_G:\widehat{G} \times \widehat{G}\ra L_{\d}(G)$$ 
 which on $R$-points is given by:
 $$\{\ ,\ \}_G:G \times G\ra L_{\d}(G),\ \ \ \{a,b\}_G = \sigma_G(a)\sigma_G(b)\sigma_G(ab)^{-1}.$$
Recall that by a {\it $\d_G$-horizontal} (or {\it $\phi_G$-horizontal}) subgroup scheme of $G$ we understand a smooth closed subgroup  scheme $H\subset G$ such that
 the ideal $I_H$ of $\widehat{H}$ in $\cO(G)\h$ satisfies $\d_G(I_H)\subset I_H$; the latter is equivalent to $\phi_G(I_H)\subset I_H$ (due to the flatness of $H$).
 The indices $G$ in $\phi_G, \d_G, \sigma_G, \{\ ,\ \}_G$ are not meant to indicate that these objects are attached in any natural way to $G$; rather the index is meant to suggest that these objects are to be considered as part of the data whenever $G$ is being considered. Dropping the index $G$ from $\phi_G,\d_G$ would  create confusion even if the reference to $G$ was clear, because $\phi$ and $\d$ have other meanings;
 but we will often drop the index $G$ from $\sigma_G, \{\ ,\ \}_G$ whenever the reference to $G$ is clear.
 
\begin{definition}\label{high}
A $\d$-map $f:G\ra L_{\d}(G)$ is called a  {\it skew $\d$-cocycle} if it satisfies 
$$f(ab)=(a\star_{\d} f(b)) \cdot f(a)\cdot  \{a,b\}$$
for all $a,b \in G$.
We say $f$ is {\it  $\d$-coherent} if for any $\d_G$-horizontal  closed subgroup scheme $H\subset G$ we have that $f(H)\subset L_{\d}(H)$.
\end{definition}

It would be interesting to have a description of all ($\d$-coherent) skew $\d$-cocycles for $G=GL_n$, say, which is analogous to the description of all {\it classical $\d$-coherent $\d$-cocycles} for $GL_n$ in \cite{adel1}.
In any case examples of $\d$-coherent skew $\d$-cocycles will be given below and one might expect that these examples are essentially the only ones.

If $f$ is a fixed skew $\d$-cocycle then 
we set
$$G^f=\{b\in G;f(b)=1\}.$$
This is not a subgroup of $G$ in general.
For all $a\in G$ and all $b\in  G^f$ with $\{a,b\}=1$ we have $f(ab)=f(a)$. 

\begin{definition}
The map
$l\d:G\ra L_{\d}(G)$,
given  on $R$-points  by 
$$l\d(a)=\nabla^1(a)\sigma(a)^{-1},\ a\in G$$
 (operations performed in $J^1(G)$) is called 
 the {\it arithmetic logarithmic derivative} attached to $\phi_G$.
\end{definition}

 The map $l\d$ above is a $\d$-map of order $1$ and is easily seen to be a skew $\d$-coherent $\d$-cocycle; we view it  as an analogue of Kolchin's logarithmic derivative.
  In view of the Remark above
  for any $b\in G^{l\d}$ with $\{a,b\}=1$ we have $l\d(ab)=l\d(a)$. 
  
By the way, the map of $p$-formal schemes 
$$l\d^1:J^1(G)\ra L_{\d}(G)$$ inducing the $\d$-map 
$l\d:G\ra L_{\d}(G)$ is given by $l\d^1(g)=g\cdot \sigma(\pi(g))^{-1}$; indeed, clearly, $l\d^1(\nabla^1(a))=l\d (a)$ i.e. $(l\d^1)_*=l\d$.
   The map $l\d^1$ should be viewed as an analogue of the  Maurer-Cartan connection \cite{sharpe}, p. 98. Indeed the Maurer-Cartan connection is a natural morphism of schemes $T(G)\ra L(G)$ from the tangent bundle $T(G)$ of $G$ to the Lie algebra of $G$ which induces linear isomorphisms between the fibers of $T(G)\ra G$ and $L(G)$.   On the other hand $J^1(G)$ is an arithmetic analogue of $T(G)$ and $L_{\d}(G)$ is an arithmetic analogue of $L(G)$;
   also $l\d^1$ induces  left $L_{\d}(G)$-equivariant bijections between the fibers of $J^1(G)\ra \widehat{G}$ and $L_{\d}(G)$.
   
    \begin{definition}
  Let $G$ be a smooth affine group scheme equipped with a lift of Frobenius $\phi_G$ and let  $\alpha \in L_{\d}(G)$.
  Then  the equality
  $$l\d(u)=\alpha$$ 
  will be referred to as a {\it $\d_G$-linear equation} with unknown  $u\in G$.
  \end{definition}

 Finally, let us discuss compatibility of lifts of Frobenius with 
 
 1. left and right translation by subgroups; 
 
 2. ``outer" automorphisms (in particular involutions) and their symmetric spaces;
 
 3. ``inner" automorphisms and conjugacy classes.
  
  We start with the first type of compatibility. Recall from the Introduction the following:
 
  \begin{definition}
  Let $H\subset G$ be a smooth closed subgroup scheme.
 We say that $\phi_G$ (or, equivalenty $\d_G$) is 
 {\it left compatible} with $H$ if 
  $H$ is $\d_G$-horizontal and  the diagram \ref{vigil1} is commutative.
We say that $\phi_G$ (or, equivalenty $\d_G$) is 
 {\it right compatible} with $H$ if 
  $H$ is $\d_G$-horizontal and  the diagram \ref{vigil2} is commutative.
  We say that $\phi_G$ is {\it bicompatible} with $H$ if it is both left and right compatible with $H$.
Note that the vertical maps in the diagrams \ref{vigil1} and \ref{vigil2} are morphisms of formal schemes over $\bZ_p$ and not over $R$!
 \end{definition}

 \begin{remark}
 By a split torus over $R$ we will understand a group scheme of the form $T=Spec\ R[t_1,t_1^{-1},...,t_n,t^{-1}_n]$ with $t_i$ group like elements, i.e. elements defining multiplicative characters. There is exactly one lift of Frobenius $\phi_T$ on $\widehat{T}$ which is bicompatible with $T$ itself.
 This $\phi_T$  acts as $\phi_T(t_i)=t_i^p$.
 \end{remark}
 
 \begin{remark}
 By a {\it finite constant group scheme} over $R$ we will understand a group scheme which is a disjoint union $W$ of schemes $Spec\ R$ indexed by some finite abstract group, with the obvious structure of group scheme. For any finite constant group scheme $W$ there is exactly one lift of Frobenius  on $\cO(W)\h= R\times...\times R$, hence on $\widehat{W}$. Clearly this lift of Frobenius is bicompatible with $W$ itself.
 \end{remark}
 
 \begin{remark}
 More generally one can consider semidirect products $N$ of a split torus $T$ and a finite constant group scheme $W$ acting on $T$. Then there is exactly one lift of Frobenius on $\widehat{N}$ that is bicompatible with $N$. 
 \end{remark}

\begin{lemma}
\label{deal}
Assume one of the following holds:

1)  $\phi_G$ is left compatible with $H$, $a\in H$, $b\in G$;

2)  $\phi_G$ is right compatible with $H$, $a\in G$, $b\in H$;

Then 
$\{a,b\}_G=1$ and, in particular, $$l\d(ab)=(a\star_{\d} l\d(b)) \cdot l\d(a).$$
\end{lemma}

{\it Proof}. Clear from the definitions. 
\qed

\bigskip

Next we discuss compatibility with automorphisms. Recall from the Introduction the following:

\begin{definition}
\label{sleep}
Let $G$ be an affine smooth group scheme over $R$ equipped with  a lift of Frobenius $\phi_{G,0}$ on $\widehat{G}$.
For an automorphism $\tau$ of $G$ one can attach the closed subgroup $G^+\subset G$, the closed subscheme $G^-\subset G$,  the {\it quadratic map} $\cH$, and the {\it bilinear map} ${\mathcal B}$ as in \ref{Gtaujos}, \ref{Gtausus}, \ref{hash},
 \ref{scorpion4}. Also for $g\in G$ consider the maps $\cH_g$ and ${\mathcal B}_g$ in \ref{gilmore} and \ref{scorpion42} respectively.
 The identity component $S:=(G^+)^{\circ}$ is called the {\it subgroup  defined by $\tau$}. On the other hand $G^-$ is not a subgroup of $G$; recall, however, that $G^-$ is stable under the maps $G\ra G$, $a\mapsto a^{\nu}$ for all $\nu\in \bZ$.
 Consider   lifts of Frobenius $\phi_G$ and $\phi_{G,0}$ on $\widehat{G}$.
We say that $\phi_{G}$ is {\it ${\mathcal H}_g$-horizontal} with respect to $\phi_{G,0}$ if the diagram  \ref{sleepless} is commutative; we say that $\phi_G$ is  {\it ${\mathcal B}_g$-symmetric} with respect to $\phi_{G,0}$ if the diagram  \ref{seatle} is commutative.
 Note that if $\phi_G$ is ${\mathcal H}$-horizontal with respect to $\phi_{G,0}$ and $\phi_{G,0}(1)=1$ then $S:=(G^+)^{\circ}$ is $\phi_{G}$-horizontal and $\phi_G(1)\in S$; in particular there is an induced lift of Frobenius $\phi_S$ on $\widehat{S}$.  \end{definition}

\begin{remark}
Assume in this remark that $\tau$ is an involution.
The geometry of $G,G^-,G^+$ above  is, morally, analogous to that appearing in the theory of symmetric spaces \cite{helgason}, p.209. Let us examine this geometry more closely in what follows.
Indeed, in the notation above, for any $R$-point $a\in G$ we have an equality of sets of $R$-points $\cH^{-1}(\cH(a))=G^+\cdot a$. Moreover  the morphism $\cH:G\ra G$ factors through the embedding $G^-\subset G$ hence induces a morphism $\cH:G\ra G^-$. In particular $\cH$ induces an injection of sets of $R$-points
$$G^+\backslash G\ra G^-.$$ 
Classically the pair $(G,G^+)$ is referred to as a {\it symmetric pair} and the quotient $G^+\backslash G$ is referred to as a {\it symmetric space}.
Furthermore $(G^-)^{\tau}\subset G^-$, and we have a commutative diagram
$$
\begin{array}{rcl}
G & \stackrel{\tau}{\longrightarrow} & G\\
\cH  \downarrow & \  & \downarrow \cH\\
G^- & \stackrel{\tau}{\longrightarrow} & G^-
\end{array}
$$
Denote by an upper bar the reduction mod $p$ functor. Then the set theoretic fibers $\overline{\cH}^{-1}(\overline{\cH}(\overline{a}))$ of the morphism $\overline{\cH}:\overline{G}\ra \overline{G}$ are still the set theoretic right cosets $\overline{G^+}\cdot \overline{a}$. If we denote by $L(\ )$ the tangent space at the identity $1$ of a closed subscheme of $\overline{G}$ that passes through $1$ then $L(\overline{G^+})$ identifies with the $+1$ eigenspace $L(\overline{G})^+$ of $L(\overline{\tau})$ in $L(\overline{G})$ while
$L(\overline{G^-})$ identifies with the $-1$ eigenspace $L(\overline{G})^-$ of $L(\overline{\tau})$ in $L(\overline{G})$. Furthermore (since the characteristic is not $2$) we have a direct sum decomposition 
$$L(\overline{G})=L(\overline{G})^+\oplus L(\overline{G})^-.$$
On the above space  $L(\overline{\tau})$ acts as $(b^+,b^-)\mapsto (b^+,-b^-)$;
moreover the map $L(\overline{\cH}):L(\overline{G})\ra L(\overline{G^-})$ sends $(b^+,b^-)\mapsto 2b^-$ so it is surjective with kernel $L(\overline{G^+})$. Assume in what follows that $\overline{G}$ is connected and $\overline{G^+}$ is smooth and  denote by $(\overline{G^-})^{\circ}$ the connected component of $\overline{G^-}$ containing $1$. Then we claim that 
$\overline{\cH}:\overline{G}\ra (\overline{G^-})^{\circ}$ is dominant and its image
is contained in the smooth locus 
$(\overline{G^-})^{\circ}_{s}$ of $(\overline{G^-})^{\circ}$.
 Indeed the set theoretic image of
$\overline{\cH}:\overline{G}\ra (\overline{G^-})^{\circ}$ can be identified with  the quotient $\overline{G^+}\backslash\overline{G}$; on the other hand 
we have
$$\begin{array}{rcl}
\dim(\overline{G^+}\backslash\overline{G}) & = & \dim(\overline{G})-\dim(\overline{G^+})\\
\  & = & \dim(L(\overline{G}))-\dim(L(\overline{G})^+),\ \ \text{by smoothness of $\overline{G^+}$}\\
\  & = & \dim(L(\overline{G})^-)\\
\  & \geq & \dim_1(\overline{G^-})\\
\  & \geq &  \dim(\overline{G^+}\backslash\overline{G})
\end{array}
$$
which forces the inequalities above to be equalities. This implies that the local ring at $1$ of the scheme $\overline{G^-}$ is regular so $1\in (\overline{G^-})^{\circ}_{s}$. So  $(\overline{G^-})^{\circ}$  is irreducible and of the same dimension as $\overline{G^+}\backslash\overline{G}$, hence the latter is dense in $(\overline{G^-})^{\circ}$.
Since $\overline{G}$ acts on the right on $(\overline{G^-})^{\circ}$ via $(a,b)\mapsto b^{-\tau}ab$ we get that the whole of the set
$\overline{G^+}\backslash\overline{G}$ is contained in 
$(\overline{G^-})^{\circ}_{s}$. So we have an induced dominant
morphism 
$\overline{\cH}:\overline{G}\ra (\overline{G^-})^{\circ}_s$.
 Since  $(\overline{G^-})^{\circ}_s$ is an integral scheme the map $\cO((\overline{G^-})^{\circ}_s)\ra \cO(\overline{G})$ induced by $\overline{\cH}$ is injective. Hence  the map
$\cO((\widehat{G^-})^{\circ}_s)\ra \cO(\widehat{G})$ induced by $\widehat{\cH}:\widehat{G}\ra (\widehat{G^-})^{\circ}_s$ is injective. 
In particular, if there exists $\phi_G$ that is ${\mathcal H}$-horizontal with respect to   $\phi_{G,0}$ then $\phi_{G,0}$ sends $(\widehat{G^-})^{\circ}_s$ into itself.
\end{remark}

\begin{remark}
Here is an alternative approach to the commutativity of diagrams \ref{sleepless} and  \ref{seatle}. Consider the map 
\begin{equation}
\label{Q}
{\mathcal Q}:G\times G \ra G,\ \ \ {\mathcal Q}(x,y)=x^{-1}y.
\end{equation}
 For any two $\phi$-linear morphisms $f_1,f_2:\widehat{G}\ra \widehat{G}$
of $p$-formal schemes  one can consider the $\phi$-linear morphism
$f_1\times f_2:\widehat{G}\ra \widehat{G}\times \widehat{G}$ (where $\times=\times_R$) and hence the composition
$$f_1^{-1}\cdot f_2:=\widehat{\mathcal Q}\circ (f_1\times f_2):\widehat{G}\ra \widehat{G}\times \widehat{G}\ra \widehat{G}.$$
Similarly if ${\mathcal Q}^{\prime}:G\times G \ra G,\ \ \ {\mathcal Q}^{\prime}(x,y)=xy^{-1}$
then we set
$$f_1\cdot f_2^{-1}:=\widehat{\mathcal Q}^{\prime}\circ (f_1\times f_2):\widehat{G}\ra \widehat{G}\times \widehat{G}\ra \widehat{G}.$$

\end{remark}

We have the following:

\begin{lemma}
\label{policeman}
Let $G$ be an affine smooth group scheme over $R$ and $g\in G$.
Let $\tau$ be an involution and let $\phi_{G,0}, \phi_G$ be  lifts of Frobenius on $\widehat{G}$. 

1) $\phi_G$ is ${\mathcal B}_g$-symmetric with respect to $\phi_{G,0}$ if and only if \begin{equation}
\label{sqqq}
(\phi_G \cdot \phi_{G,0}^{-1})^{-\tau}=\phi_G\cdot \phi_{G,0}^{-1}.\end{equation}

2) Assume  $\phi_G$ is ${\mathcal H}_g$-horizontal with respect to $\phi_{G,0}$.
Then $\phi_G$ is ${\mathcal B}_g$-symmetric with respect to $\phi_{G,0}$ if and only if we have the following equality of maps $\widehat{G}\ra \widehat{G}$:
\begin{equation}
\label{square}
(\phi_{G,0}^{-1}\cdot \phi_G)^2=({\mathcal H}_g\circ \phi_{G,0})^{-1}\cdot (\phi_{G,0}\circ \cH_g).\end{equation}
\end{lemma}

In the above statement 
$(\phi_G\cdot \phi_{G,0}^{-1})^{-\tau}$ is the composition of $\phi_G \cdot \phi_{G,0}^{-1}:\widehat{G}\ra \widehat{G}$ with the  map $\widehat{G}\ra \widehat{G}$, $x\mapsto x^{-\tau}$.
Similarly
$(\phi_{G,0}^{-1}\cdot \phi_G)^2=\phi_{G,0}^{-1}\cdot \phi_G\cdot \phi_{G,0}^{-1}\cdot \phi_G$ is the composition of $\phi_{G,0}^{-1}\cdot \phi_G:\widehat{G}\ra \widehat{G}$ with the squaring map $\widehat{G}\ra \widehat{G}$, $x\mapsto x^2$.

\medskip

{\it Proof}. To check 1) note that, by diagram \ref{seatle}, we have 
$$\phi_G(a)^{-\tau}\cdot  \phi_{G,0}(a)=\phi_{G,0}(a)^{-\tau}\cdot \phi_G(a),\ \ \ a\in G.$$
This immediately implies \ref{sqqq} evaluated at $a$.
Similarly, to check 2)
note that we have ${\mathcal B}_g(x,y)={\mathcal H}_g(x)\cdot {\mathcal Q}(x,y)$. Hence 
$${\mathcal B}_g\circ (\phi_G\times \phi_{G,0})=({\mathcal H}_g\circ \phi_G) \cdot (\phi_{G}^{-1}\cdot \phi_{G,0})=(\phi_{G,0}\circ {\mathcal H}_g) \cdot (\phi_{G}^{-1}\cdot \phi_{G,0}),$$
$${\mathcal B}_g\circ (\phi_{G,0}\times \phi_{G})=({\mathcal H}_g\circ \phi_{G,0}) \cdot (\phi_{G,0}^{-1}\cdot \phi_{G}).$$
So the diagram \ref{seatle} is commutative if and only if equality \ref{square} holds.
\qed

\medskip

We will use Lemma \ref{policeman} in combination with  the following:

\begin{lemma}
\label{below}
Let ${\mathcal A}$ be a $p$-adically separated flat $R$-algebra, let $\nu\in \bZ$ be an integer not divisible by $p$,  and let  $M_1,M_2\in GL_n({\mathcal A})$ two matrices such that $M_1^{\nu}=M_2^{\nu}$ and $M_1\equiv M_2\equiv 1$ mod $p$. Then $M_1=M_2$.
\end{lemma}

{\it Proof}.
We prove by induction on $n$ that $M_1\equiv M_2$ mod $p^n$. Indeed assume the latter and write $M_1=M_2+p^n M$. Then
$$M_2^{\nu}=M_1^{\nu}\equiv M_2^{\nu}+p^n \left(\sum_{i+j=\nu-1} M_2^iMM_2^j\right)\ \ \ mod\ \ \ p^{n+1}.$$
Since $M_2\equiv 1$ mod $p$ we get $p^n\nu M\equiv 0$ mod $p^{n+1}$.
Since ${\mathcal A}$ is $R$-flat we get $\nu M\equiv 0$ mod $p$ and hence, since $\nu\in R^{\times}$, we get $M\equiv 0$ mod $p$ so $M_1\equiv M_2$ mod $p^{n+1}$.
\qed

\medskip

In the next corollary a group scheme is called linear if it is isomorphic to a closed subgroup scheme of $GL_n$.

\begin{corollary}
\label{unique}
Let $G$ be an smooth linear group scheme over $R$. Let ${\tau}:G\ra G$ be an involution and let $\phi_{G,0},\phi_{G,1},\phi_{G,2}$ be  lifts of Frobenius on $\widehat{G}$ such that $\phi_{G,1}$ and $\phi_{G,2}$ are ${\mathcal H}_g$-horizontal  and ${\mathcal B}_g$-symmetric with respect to $\phi_{G,0}$.  Then $\phi_{G,1}=\phi_{G,2}$.
\end{corollary}

{\it Proof}.
By Lemma \ref{policeman} we have 
$$(\phi_{G,0}^{-1}\cdot \phi_{G,1})^2=(\phi_{G,0}^{-1}\cdot \phi_{G,2})^2.$$
We also have 
$$\phi_{G,0}^{-1}\cdot \phi_{G,1}\equiv \phi_{G,0}^{-1}\cdot \phi_{G,2}\equiv 1 \ \ \ mod\ \ \ p.$$
View $G$ as embedded into some $GL_n$. Interpreting the latter congruence as a congruence between matrices with entries in the ring $\cO(G)\h$ it follows from
 Lemma \ref{below}   that
$$\phi_{G,0}^{-1}\cdot \phi_{G,1}= \phi_{G,0}^{-1}\cdot \phi_{G,2}$$
and we are done.
\qed

\medskip

Here is another consequence of Lemma \ref{policeman}.  Assume $\tau$ is an involution and  consider the naturally induced involution $L_{\d}(\tau)$ on $J^1(G)$; it preserves $L_{\d}(G)$ so it induces an involution $L_{\d}(\tau)$ on $L_{\d}(G)$. When there is no danger of confusion we continue to denote by $\tau$ the involution $L_{\d}(\tau)$. Set, as usual, 
$$L_{\d}(G)^+=\{b\in L_{\d}(G);b^{\tau}=b\},$$
$$L_{\d}(G)^-=\{b\in L_{\d}(G);b^{\tau}=b^{-1}\}.$$
Again $L_{\d}(G)^+$ is a subgroup of $L_{\d}(G)$ but $L_{\d}(G)^-$ is not a subgroup of $L_{\d}(G)$ in general; however  $L_{\d}(G)^-$ is, again, stable under the  maps $L_{\d}(G)\ra L_{\d}(G)$, $b\mapsto b^{\nu}$ for all $\nu\in \bZ$. Note also that $L_{\d}(G^+)$ is a subgroup of $L_{\d}(G)^+$.
 Let $\sigma_G, \sigma_{G,0}:G\ra J^1(G)$ be the sections
defined by the lifts of Frobenius $\phi_G,\phi_{G,0}$ and let $l\d,l\d_0:G\ra L_{\d}(G)$ be the corresponding arithmetic logarithmic derivatives. 

\begin{corollary}
\label{deputy}
The following are equivalent:

1) $\phi_G$ is ${\mathcal B}$-symmetric with respect to $\phi_{G,0}$.

2) The image of the map $\phi_G\cdot \phi_{G,0}^{-1}:G\ra G$ is contained in $G^-$.

3) The image of the map $\sigma_G\cdot \sigma_{G,0}^{-1}:G\ra L_{\d}(G)$
is contained in $L_{\d}(G)^-$.

4)  The image of the $\d$-map $(l\d_0)^{-1}\cdot (l\d):G\ra L_{\d}(G)$ is contained in $L_{\d}(G)^-$.
\end{corollary}

{\it Proof}.
The equivalence of 1) and 2) is a rephrasing of assertion 1 in Lemma \ref{policeman}. The rest of the Corollary
is a trivial exercise.
\qed

\medskip

On the other hand here is a consequence of Lemma \ref{below}:

\begin{corollary}
\label{wolf}
Let $G$ be a smooth linear group scheme. The following hold:

1) For all $\nu\in \bZ$ not divisible by $p$ the following maps are injective: $$L_{\d}(G)\ra L_{\d}(G),\ \ \ a\mapsto a^{\nu}.$$

2) If $\tau$ is an involution of $G$ then the multiplication map below is injective:
$$L_{\d}(G)^+\times L_{\d}^-(G)\ra L_{\d}(G),\ \ \ (b,c)\mapsto bc.$$ 
\end{corollary}

{\it Proof}.
To check 1) embed $G$ into some $GL_n$. Then there is a natural injective homomorphism between groups of points, $\epsilon:L_{\d}(G)\ra GL_n$; it is given by the map
$$\cO(GL_n)\h \ra \cO(G)\h\stackrel{\phi}{\longrightarrow} \cO(J^1(G))\ra \cO(L_{\d}(G)),$$
where $\phi(a)=a^p+p\d a$, $\d$ the universal $p$-derivation.  Any matrix in the image of $\epsilon$ is $\equiv 1$ mod $p$. So if $a^{\nu}=b^{\nu}$ with $a,b\in L_{\d}(G)$ then $\epsilon(a)^{\nu}=\epsilon(b)^{\nu}$. By Lemma \ref{below} we get $\epsilon(a)=\epsilon(b)$ and hence $a=b$.

To check 2) assume $b_1c_1=b_2c_2$ with $b_i^{\tau}=b_i$ and $c_i^{\tau}=c_i^{-1}$. Set $b=b_2^{-1}b_1$. Then $b^{\tau}=b$ and $b=c_2c_1^{-1}$. Hence
$$c_1c_2^{-1}=b^{-1}=b^{-\tau}=(c_2c_1^{-1})^{-\tau}=c_1^{\tau}c_2^{-\tau}=c_1^{-1}c_2.$$
Hence $c_1^2=c_2^2$. By part 1) we get $c_1=c_2$. Hence $b=1$ and hence $b_1=b_2$.
\qed

\begin{remark}
In the notation of Corollary \ref{wolf}, for any $a\in L_{\d}(G)$ in the image of the map $L_{\d}(G)^+\times L_{\d}^-(G)\ra L_{\d}(G)$ we have a unique representation $a=a^+a^-$ with $a^{\pm}\in L_{\d}(G)^{\pm}$ which we shall refer to as the {\it Cartan decomposition} of $a$ with respect to $\tau$. We will later show (cf. Proposition \ref{cartan}) that for $G=GL_n$ and $\tau$ certain remarkable involutions the map $L_{\d}(G)^+\times L_{\d}^-(G)\ra L_{\d}(G)$ is a bijection; so in that case  we will have that any element of $L_{\d}(G)$ has a Cartan decomposition.\end{remark}

\begin{remark}
The arithmetic logarithmic derivatives $l\d:G\ra L_{\d}(G)$ attached to various lifts of Frobenius $\phi_G$ all derive from a certain $\d$-map $l^1 \d:J^1(G)\ra L_{\d}(G)$ which we describe in what follows. This map was introduced in \cite{je} and will not play any role later in the present paper; but it is conceptually tightly linked to our discussion so we present below a quick overview of its main features. 

 Let $G$ be any smooth group scheme over $R$ and denote by $\pi:J^1(G)\ra G$ the canonical projection. One can consider the $\d$-map, of order $1$,
 $l^1\d:J^1(G)\ra J^1(G)$ defined on points by $l^1\d a:=\nabla^1(\pi(a))\cdot a^{-1}$, for $a\in J^1(G)$. Clearly this map takes values in $L_{\d}(G)$ so it induces a $\d$-map of order $1$
 $$l^1\d:J^1(G)\ra L_{\d}(G).$$
 The following are trivial to check:
 
 1) The map $l^1\d:J^1(G)\ra L_{\d}(G)$ is a $\d$-cocycle for the $\star$-action in the sense that
 $$l^1\d(ab)=(l^1\d a)\cdot (a\star (l^1\d b)),\ \ \ a,b\in J^1(G).$$
 
 2) For $a,b \in J^1(G)$ one has $l^1\d a=l^1\d b$ if and only if there exists $c\in G$ such that $ab^{-1}=\nabla^1 c$. In particular the composition of  $l^1\d:J^1(G)\ra L_{\d}(G)$ with the map $\nabla^1:G\ra J^1(G)$ is the constant map with value $1\in L_{\d}(G)$.
 
 3) The composition of $l^1\d:J^1(G)\ra L_{\d}(G)$ with any section $s:G\ra J^1(G)$ of the projection $\pi:J^1(G)\ra G$ equals the arithmetic logarithmic derivative $l\d:G\ra L_{\d}(G)$ attached to the lift of Frobenius defined by $s$. 
 
 4) The composition of $l^1\d:J^1(G)\ra L_{\d}(G)$ with the inclusion $L_{\d}(G)\subset J^1(G)$ is the antipode
 $L_{\d}(G)\ra L_{\d}(G)$, $a\mapsto a^{-1}$.
 
 5) The map of formal schemes 
 $$l^{1}\d^1:J^1(J^1(G))\ra L_{\d}(G)$$
 defining the $\d$-map $l^1\d:J^1(G)\ra L_{\d}(G)$ (i.e. with the property that $l^1\d=(l^{1}\d^1)_*$) can be described as follows.  
  Let $\pi^1_G:J^1(G)\ra J^0(G)$ denote the natural projection. Then one may consider
 the homomorphism
 \begin{equation}
\label{vint}
\pi^{11}_G:=\pi^1_{J^1(G)}\times J^1(\pi^1_G):J^1(J^1(G))\ra J^1(G)\times_{J^0(G)} J^1(G).
\end{equation}
On the other hand 
the quotient map
$$J^1(G)\times J^1(G)\ra J^1(G),\ \ (a,b)\mapsto ba^{-1}$$
induces a morphism 
\begin{equation}
\label{quot}J^1(G)\times_{J^0(G)} J^1(G)\ra L_{\d}(G).\end{equation}
Then the map $l^{1}\d^1$ is obtained by composing the map \ref{quot}  above with the morphism
$
\pi^{11}_G$
defined in \ref{vint}.   It was checked in
 \cite{je} that $(l^{1}\d^1)^{-1}(1)=J^2(G)$. Moreover
  the restriction of $l^{1}\d^1$ to $L_{\d}(J^1(G))$ is obtained by composing the map
$L_{\d}(\pi^1_G):L_{\d}(J^1(G))\ra L_{\d}(G)$ (induced by $J^1(\pi^1_G)$) with the antipode map $L_{\d}(G)\ra L_{\d}(G)$.

6) For any section $s:\widehat{G}\ra J^1(G)$ of the projection $\pi^1_G:J^1(G)\ra \widehat{G}$ the corresponding arithmetic logarithmic derivative $l\d:G \ra L_{\d}(G)$ is induced by the map
$$l^1\d^1\circ J^1(s):J^1(G)\ra J^1(J^1(G))\ra L_{\d}(G);$$
i.e. $(l^1\d^1\circ J^1(s))_*=l\d:G\ra L_{\d}(G)$.
\end{remark}

Next we discuss compatibility of lifts of Frobenius with conjugation.
Recall the following definition from the Introduction.

\begin{definition}
\label{late}
  Let $G$ be a smooth affine group scheme over $R$
  and let ${\mathcal C}:G \times G\ra G$,  ${\mathcal C}(h,g)=g^{-1}hg$, be the conjugation map. Let
   $H$ be a closed smooth subscheme, and  $G^*\subset G$ an open set which is invariant under the action of $G$ on $G$ by conjugation. Let $H^*=H\cap G^*$ and let 
  ${\mathcal C}:H^*\times G\ra G^*$ be the induced map. Let $\phi_{G,0}$ be a lift of Frobenius on $\widehat{G}$ and $\phi_{G^*}$ a lift of Frobenius on $\widehat{G^*}$. We say that $\phi_{G,0}$ is {\it ${\mathcal C}$-horizontal} with respect to $\phi_{G^*}$
  if $H$ is $\phi_{G,0}$-horizontal and, upon denoting by $\phi_{H^*,0}$ the  lift of Frobenius on $\widehat{H^*}$ induced by $\phi_{G,0}$, the  diagram \ref{love} in the Introduction is commutative.
 \end{definition}
 
 There is an easy relation between ${\mathcal C}$-horizontality and ${\mathcal H}$-horizontality as follows. 
 
 \begin{lemma}
 \label{soap}
 Assume in definition \ref{late} that   $\phi_{G,0}$ is ${\mathcal C}$-horizontal with respect to $\phi_{G^*}$. Let $h\in H^*$ be $\phi_{G,0}$-horizontal and let $V_h$ be the $p$-adic completion  of the Zariski closure   of the image of ${\mathcal C}_h:G\ra G^*$, $g\mapsto {\mathcal C}_h(g)=g^{-1}hg$. The following hold:
 
 1) $V_h$ is $\phi_{G^*}$-horizontal.
 
 2) Assume $G=GL_n$ and let  $\tau$ is the automorphism $x^{\tau}=h^{-1}xh$ of $G$. Assume furthermore that $V_h$ is closed in $\widehat{G}$. Then there exists a lift of Frobenius $\phi_G$ on $\widehat{G}$ such that $\phi_{G,0}$ is ${\mathcal H}$-horizontal with respect to $\phi_G$. 
 \end{lemma}
 
 {\it Proof}. To prove assertion 1 note that by combining the commutative diagram \ref{love} with the the commutative diagram
 \begin{equation}
 \label{colloquium}
 \begin{array}{rcl}
 \widehat{G} & \stackrel{\phi_{G,0}}{\longrightarrow} & \widehat{G}\\
 h\times id\downarrow & \  & \downarrow h\times id\\
 \widehat{H^*}\times \widehat{G} & \stackrel{\phi_{H^*,0}\times \phi_{G,0}}{\longrightarrow} & \widehat{H^*}\times \widehat{G}
 \end{array}\end{equation}
 one gets a commutative diagram
  \begin{equation}
 \label{colloquium1}
 \begin{array}{rcl}
 \widehat{G} & \stackrel{\phi_{G,0}}{\longrightarrow} & \widehat{G}\\
 {\mathcal C}_h\downarrow & \  & \downarrow {\mathcal C}_h\\
 \widehat{G^*} & \stackrel{\phi_{G^*}}{\longrightarrow} & \widehat{G^*}
 \end{array}\end{equation}
which immediately implies assertion 1.
 Let us check assertion 2. Putting 
 together the commutative diagram \ref{colloquium1} and the diagram
 $$\begin{array}{rcl}
 \widehat{G^*} & \stackrel{\phi_{G^*}}{\longrightarrow} & \widehat{G^*}\\
 L_h^{-1}\downarrow &\  & \downarrow L_h^{-1}\\
 h^{-1}\widehat{G^*} & \stackrel{L_h^{-1}\phi_{G^*}L_h}{\longrightarrow} & h^{-1}\widehat{G^*}
 \end{array}$$
 and noting that $L_h^{-1}\phi_{G^*}L_h$ is a lift of Frobenius we get a commutative diagram
 \begin{equation}
 \label{colloquium2}
 \begin{array}{rcl}
 \widehat{G} & \stackrel{\phi_{G,0}}{\longrightarrow} & \widehat{G}\\
 \downarrow & \  & \downarrow \\
h^{-1}V_h & \stackrel{L_h^{-1}\phi_{G^*}L_h}{\longrightarrow} & h^{-1}V_h
 \end{array}\end{equation}
Now since $V_h$ is closed in $\widehat{G}$  so is $h^{-1}V_h$ hence  by Lemma \ref{liftingg} the bottom arrow of \ref{colloquium2} extends to a lift of Frobenius $\phi_G$ on $\widehat{G}$ and we are done. 
 \qed
 
 \bigskip
 
 We end our discussion of abstract Lie theory by discussing brackets.
 Let $G$ be a smooth group scheme or group $p$-formal scheme. We may consider the
 group $p$-formal scheme, which we call the {\it order $r$ $\d$-Lie algebra} of $G$,
 $$L_{\d}^r(G):=Ker(\pi:J^r(G)\ra J^{r-1}(G))$$
 where $\pi$ is as in \ref{defofpi}.
 In particular $L_{\d}^1(G)=L_{\d}(G)$. Also we consider the
  composition 
 \begin{equation}
 \label{defofE}
 ex^r:L_{\d}^r(G)\subset J^r(G)\stackrel{\pi_{\phi}}{\longrightarrow}J^{r-1}(G)
 \stackrel{\pi_{\phi}}{\longrightarrow}...\stackrel{\pi_{\phi}}{\longrightarrow}J^1(G)
 \stackrel{\pi_{\phi}}{\longrightarrow}
  \widehat{G},
 \end{equation}
 where $\pi_{\phi}$ is as in \ref{defofpiphi}; $ex^r$ can be thought of as an ``exponential" and is a $\phi^r$-linear morphism of $p$-formal schemes. It is clearly functorial
 in $G$ with respect to homomorphisms of group ($p$-formal) schemes. As usual we still denote by $ex^r:L_{\d}^r(G)\ra G$ the map induced on sets of $R$-points. We also set $ex=ex^1$.
 In the next section we will prove the following:
 
 \begin{proposition}
 \label{bracketprop}
 Let $G$ be a smooth linear  group scheme with commutator map $[\ ,\ ]:G\times G\ra G$, $[a,b]=aba^{-1}b^{-1}$ and let  $r,s\geq 1$. 
 
 1) There exists a unique $\d$-map 
 $$[\ ,\ ]_{\d}:L_{\d}^{r}(G)\times L_{\d}^{s}(G)\ra L_{\d}^{r+s}(G)$$
 of order $\max\{r,s\}$ making the following diagram of sets of $R$-points commutative
 \begin{equation}
 \label{telephone}
 \begin{array}{rcl}
 L_{\d}^{r}(G)\times L_{\d}^{s}(G) & \stackrel{[\ ,\ ]_{\d}}{\longrightarrow} &
  L_{\d}^{r+s}(G)\\
  ex^{r}\times ex^{s} \downarrow & \  & \downarrow ex^{r+s}\\
  G \times G & \stackrel{[\ ,\ ]}{\longrightarrow} & G
 \end{array}
 \end{equation}
 
 2) The $\d$-map $[\ ,\ ]_{\d}$ is functorial in $G$ with respect to homomorphisms $f:G\ra G'$ of group schemes, i.e., for the induced maps $L^i_{\d}(f):L_{\d}^i(G)\ra L_{\d}^i(G)$, $i=r,s$, we have 
 $$[L_{\d}^r(\alpha),L_{\d}^s(\beta)]_{\d}=L_{\d}^{r+s}([\alpha,\beta]_{\d}).$$
 
 3) If $\alpha_1,\alpha_2\in L_{\d}^r(G)$ and $\beta_1,\beta_2\in L_{\d}^s(G)$ are such that $\alpha_1\equiv \alpha_2$ mod $p^m$ and $\beta_1\equiv \beta_2$ mod $p^m$ then
 $$[\alpha_1,\beta_1]_{\d}\equiv [\alpha_2,\beta_2]_{\d}\ \ \ \text{mod}\ \ p^m.$$
 
 4) (Linearity) If $\alpha_1,\alpha_2\in L^r_{\d}(G)$, $\beta\in L^s_{\d}(G)$ then
 $$[\alpha_1+_{\d,*}\alpha_2,\beta]_{\d}\equiv [\alpha_1,\beta]_{\d}+_{\d,*}[\alpha_2,\beta]_{\d}\ \ \ \text{mod}\ \ \ p.$$

 5) (Antisymmetry) If $\alpha\in L^r_{\d}(G)$, $\beta\in L_{\d}^s(G)$  then
 $$[\alpha,\beta]_{\d}+_{\d,*}[\beta,\alpha]_{\d}\equiv 0 \ \ \text{mod}\ \ p.$$

6) (Jacobi identity) If $\alpha\in L^r_{\d}(G)$, $\beta\in L_{\d}^s(G)$,  $\gamma\in L_{\d}^t(G)$  then
$$[[\alpha,\beta]_{\d},\gamma]_{\d}+_{\d,*}[[\beta,\gamma]_{\d},\alpha]_{\d}+_{\d,*}
[[\beta,\gamma]_{\d},\alpha]_{\d}\equiv 0 \ \ \text{mod}\ \ p.$$

 7) Let $G=GL_n$ and let $[\ ,\ ]:{\mathfrak g}{\mathfrak l}_n\times {\mathfrak g}{\mathfrak l}_n\ra {\mathfrak g}{\mathfrak l}_n$ be the commutator map $[a,b]=ab-ba$. Then for any $\alpha\in L_{\d}^{r}(GL_n)$ and $\beta\in L_{\d}^{s}(GL_n)$ we have a congruence in  ${\mathfrak g}{\mathfrak l}_n$:
 $$[\alpha,\beta]_{\d} \equiv [\alpha^{(p^{s})}, \beta^{(p^{r})}]\ \ \text{mod}\ \ p.$$
 \end{proposition}

As a matter of notation used in the statement above if $a$ and $b$ are two $R$-points of an $R$-scheme we write $a\equiv b$ mod $p^m$ if the $R/p^mR$-points  induced by $a$ and $b$ coincide.
Note also  that the map of sets $ex^{r}\times ex^{s}$ is not generally induced by a morphism of $p$-formal schemes if $r\neq s$; but it is induced by a ($\phi^{r}$-linear) morphism of $p$-formal schemes if $r=s$. On the other hand the map $ex^{r+s}$ is induced by a $\phi^{r+s}$-linear morphism of $p$-formal schemes. 

\begin{remark}
The higher order $\d$-Lie algebras $L_{\d}^r(G)$ have an alternative description as follows. Let us define by induction $L_{\d}^{\circ 1}(G):=L_{\d}(G)$ and
$$L_{\d}^{\circ r}(G):=L_{\d}(L_{\d}^{\circ (r-1)}(G)).$$
By the above discussion we have natural morphisms
$$ex^{\circ r}:L_{\d}^{\circ r}(G)\stackrel{ex}{\longrightarrow} L_{\d}^{\circ (r-1)}\stackrel{ex}{\longrightarrow} ...\stackrel{ex}{\longrightarrow}\widehat{G}.$$ Then we claim there exist unique isomorphisms  $L_{\d}^{\circ r}(G)  \simeq  L_{\d}^r(G)$ 
of group $p$-formal schemes fitting into 
commutative diagrams of abstract groups
$$
\begin{array}{rcl}
L_{\d}^{\circ r}(G) & \simeq & L_{\d}^r(G)\\
ex^{\circ r}\downarrow & \  & \downarrow ex^{r}\\
\widehat{G} & = & \widehat{G}
\end{array}
$$
This is an easy consequence of \cite{char}, Proposition 2.2.
\end{remark}

\section{General linear group}

In this section we specialize the ``abstract" arithmetic Lie theory to the case of  $GL_n$. So let
 $$G=GL_n=Spec\ R[x,\det(x)^{-1}]$$ be viewed as a group scheme over $R$, where $x=(x_{ij})$ is a matrix of indeterminates.  Consider an arbitrary  lift of Frobenius 
 $$\phi_G:R[x,\det(x)^{-1}]\h\ra R[x,\det(x)^{-1}]\h.$$
 Set $$\phi_G(x_{ij})=\Phi_{ij}(x)\in R[x,\det(x)^{-1}]\h,$$ 
 $$\d_G(x_{ij})=\Delta_{ij}(x)\in R[x,\det(x)^{-1}]\h.$$ 
 So
 $\Phi_{ij}(x)=x^p_{ij}+p\Delta_{ij}(x)$.
 Write $\Phi(x)$ and $\Delta(x)$ for the matrices $(\Phi_{ij}(x))$ and  $(\Delta_{ij}(x))$ respectively. We refer to $\Delta(x)$ as the {\it Christoffel symbol}. (Conversely any choice of a matrix $\Delta$ defines a unique lift of Frobenius $\phi_{G}$.) As usual, view the Lie algebra 
$$L(G)={\mathfrak g}={\mathfrak g}{\mathfrak l}_n=Spec\ R[x']$$ as a vector group over $R$ where $x'$ is a matrix of indeterminates. We identify
$$J^1(G)=Spf\ R[x,x',\det(x)^{-1}]\h,\ \ \d x=x'.$$
If $u,v$ are two matrices of indeterminates which are coordinates on $G\times G$ then the multiplication  $G\times G\ra G$ is given, of course, by the map $x\mapsto uv$. This map 
induces, by universality,  a  multiplication map on $J^1(G)$ induced by the map
$$\mu:R[x,x',\det(x)^{-1}]\h\ra R[u,v,u',v',\det(u)^{-1},\det(v)^{-1}]\h$$
with $\mu(x)=uv$ and $\mu(x')=\d(uv)$. 
If $u=(u_{ij})$ then we set $u^{(p)}=(u_{ij}^p)$ the matrix whose entries are the $p$-th powers of the entries of $u$, we set $\phi(u)=(\phi(u_{ij}))$,  and we set $\d u=u'=(\d u_{ij})$ the matrix whose entries are obtained from the entries of $u$ by applying $\d$; hence we have $$\phi(u)=u^{(p)}+p\d u;$$ and similarly for $v$. Also, using the above notation we may write
\begin{equation}
\label{phig}
\Phi(x)=x^{(p)}+p\Delta(x).\end{equation}
In particular we have $\phi_G(x)=\Phi(x)$. From the identity of matrices $\phi(uv)=\phi(u)\phi(v)$ we get: 
\begin{equation}
\label{you}
\d(uv)=u^{(p)}v'+u' v^{(p)} +p u'v'+p^{-1}(u^{(p)}v^{(p)}-(uv)^{(p)}),\end{equation}
\begin{equation}
\label{mee} 
\d(u^{-1})=-\phi(u)^{-1}(u'(u^{-1})^{(p)}+p^{-1}(u^{(p)}(u^{-1})^{(p)}-1)).
\end{equation}
We have
$$L_{\d}(G):=Ker(J^1(G)\ra \widehat{G})=Spf\ R[x']\h$$ 
and consequently we have (canonical !) identifications as $p$-formal schemes (but not as groups)
$$J^1(G)\simeq \widehat{G}\times \widehat{{\mathfrak g}},\ \ \ L_{\d}(G) \simeq \widehat{{\mathfrak g}}.$$ 
(For groups other than $GL_n$ there is still an identification as above but it is not canonical in general !)
The set  of points of $J^1(G)$ will be identified with pairs $(a_0,a_1)\in G\times {\mathfrak g}$; the set of $R$-points of $L_{\d}(G)$ will be identified with the set of pairs $(1,a)$ with $a\in {\mathfrak g}$ and hence with ${\mathfrak g}$. 
More generally, and similarly, we have a natural identification of $p$-formal schemes
$$J^r(G)\simeq Spf\ R[x,x',...,x^{(r)},\det(x)^{-1}]\h\simeq \widehat{G}\times \widehat{\mathfrak g}^r,$$
$$L_{\d}^r(G)\simeq \widehat{\mathfrak g}.$$
Under these identifications:

$\bullet$ the group structure on $J^1(G)$ induces a group structure on the set  $G\times {\mathfrak g}$; we denote the multiplication and inverse on $G \times {\mathfrak g}$ by $\circ$ and $\iota$.

$\bullet$  the  action  $\star_{\d}:G\times L_{\d}(G)\ra L_{\d}(G)$ induces the action (still denoted by)
 $$\star_{\d}:G\times {\mathfrak g}\ra {\mathfrak g}$$

 $\bullet$ the group structure 
on $L_{\d}^r(G)$ induces a group structure on ${\mathfrak g}$, 
$$+_{\d,r}:{\mathfrak g}\times {\mathfrak g}\ra {\mathfrak g}$$

$\bullet$   $\nabla^1:\widehat{G}\ra J^1(G)$ induces a homomorphism (still denoted by) 
$$\nabla^1:G\ra G\times {\mathfrak g}$$

$\bullet$  the map $\sigma:\widehat{G}\ra J^1(G)$ induces a map (still denoted by) 
$$\sigma:G\ra G \times {\mathfrak g}$$

$\bullet$ for $\sigma:\widehat{G}\ra J^1(G)$ as above, and $\alpha\in L_{\d}(G)={\mathfrak g}$,  the map $\sigma^{\alpha}:=\alpha \cdot \sigma:\widehat{G}\ra J^1(G)$ induces a map (still denoted by) $$\sigma^{\alpha}:G\ra G\times {\mathfrak g}.$$

$\bullet$ for $\sigma:\widehat{G}\ra J^1(G)$ as above, and $u\in G$,  the map $\sigma^{\star u}:\widehat{G}\ra J^1(G)$ defined by $\sigma^{\star u}(v)=\nabla(u)^{-1}\cdot \sigma(uv)$ induces a map (still denoted by) $$\sigma^{\star u}:G\ra G\times {\mathfrak g}.$$

$\bullet$ the arithmetic logarithmic derivative
$l\d:G\ra L_{\d}(G)$
 induces the map (still denoted by)
$$l\d:G\ra {\mathfrak g}$$

$\bullet$  the map $\{\ ,\ \}:G\times G\ra L_{\d}(G)$ induces a map (still denoted by)
$$\{\ ,\ \}:G\times G\ra {\mathfrak g}.$$

$\bullet$ The map $ex^r:L_{\d}^r(G)\ra G$ induces a group homomorphism
$$ex^r:{\mathfrak g}\ra G.$$

The following proposition provides a list of formulas computing the maps above.
In the proposition below $+$ and $\cdot$ (or simply juxtaposition) in ${\mathfrak g}$ denote addition and multiplication in the associative ring ${\mathfrak g}$ of $n\times n$ matrices.

\begin{proposition} \label{hair}
Let $G=GL_n$.

 1) The multiplication $\circ$ on $G\times {\mathfrak g}$ is given
by:
$$
(a_0,a_1) \circ (b_0,b_1)= (a_0b_0,a_0^{(p)}b_1+a_1 b_0^{(p)} +p a_1 b_1+p^{-1}(a_0^{(p)}b_0^{(p)}-(a_0b_0)^{(p)})).$$

2) The identity  on $G\times {\mathfrak g}$ is the pair $$(1,0).$$

3) The inverse $\iota$ map on $G\times {\mathfrak g}$ is given by
$$\iota(a_0,a_1)=(a_0^{-1},
-(a_0^{(p)}+pa_1)^{-1}(a_1(a_0^{-1})^{(p)}+p^{-1}(a_0^{(p)}(a_0^{-1})^{(p)}-1))).
$$

4) The homomorphism $\nabla^1:G\ra G\times {\mathfrak g}$  is given   by
$$\nabla^1(a)=(a,\d a).$$

5) The multiplication $+_{\d,r}$ on ${\mathfrak g}$
 is given by
$$a+_{\d,r}b:=a+b+p^rab.$$

6) The action
 $\star_{\d}:G\times {\mathfrak g}\ra {\mathfrak g}$ is given by
 $$a\star_{\d} b=\phi(a)\cdot b\cdot \phi(a)^{-1}.$$

7)
 The arithmetic logarithmic derivative  $l\d:G\ra {\mathfrak g}$ satisfies
$$
\label{cris}l\d(ab)=(\phi(a)\cdot l\d(b) \cdot \phi(a)^{-1}) +_{\d} l\d(a) +_{\d} \{a,b\}.$$

 8) The map
 $\sigma:G\ra G \times {\mathfrak g}$    is given by
$$\sigma(a)=(a,\Delta(a)).$$ 

9) The map
$l\d:G\ra {\mathfrak g}$ is given by
$$l\d a:=\frac{1}{p}\left(\phi(a)\Phi(a)^{-1} -1 \right)=(\d a-\Delta(a))(a^{(p)}+p\Delta(a))^{-1}.$$

10) The map
$\{\ ,\ \}:G\times G\ra {\mathfrak g}$
is given by
$$\{a,b\}=p^{-1}(\Phi(a)\Phi(b)\Phi(ab)^{-1}-1).$$

11) If $\epsilon=1+p\alpha$ then the $\d_G$-linear equation 
$l\d (u)=\alpha$ is equivalent to the equation
\begin{equation}\label{clackk} 
\d u = \Delta^{\alpha}(u)\end{equation}
and also to the equation
\begin{equation}
\phi(u) = \Phi^{\alpha}(u),
\end{equation}
where $\Delta^{\alpha}(x)  :=   \alpha \cdot \Phi(x) +\Delta(x)$, $\Phi^{\alpha}(x) = \epsilon \cdot \Phi(x)$,
 $\epsilon=1+p\alpha$. (We shall refer to  Equation \ref{clackk} as a $\Delta$-linear equation.)
 
 12) For $\alpha\in {\mathfrak g}$ the map $\sigma^{\alpha}:G\ra G\times {\mathfrak g}$ is given by
 $$\sigma^{\alpha}(a)=(a,\Delta^{\alpha}(a)).$$
 (So we can refer to $\{\Delta^{\alpha}(x);\alpha\in {\mathfrak g}\}\subset {\mathfrak g}(R[x,\det(x)^{-1}]\h)$ as the linear equivalence class of $\Delta(x)$.)
 
 13) For $u\in G$ the map $\sigma^{\star u}:G\ra G\times {\mathfrak g}$ is given by
 $$\sigma^{\star u}(a)=(a,\Delta^{\star u}(a)),$$
 where $\Delta^{\star u}(x)=p^{-1}(\Phi^{\star u}(x)-x^{(p)})$, $\Phi^{\star u}(x)=\phi(u)^{-1}\Phi(ux)$.
 (So in particular if $\d u=\Delta^{\alpha}(u)$ then $(\Phi^{\alpha})^{\star u}(x)=\Phi(u)^{-1}\Phi(ux)$.)
 
 14) The group homomorphism $ex^r:{\mathfrak g}\ra G$ is given by
 $$ex^r(a)=1+p^r\phi^{-r}(a).$$
\end{proposition} 

{\it Proof}. Assertions 1, 2, 3  follow directly from \ref{you} and \ref{mee}. Assertion 4 follows directly from definitions.  Assertions 12, 13 follow from assertion 1.
Assertion 5 for $r=1$ follows from assertion 1; for arbitrary $r$ assertion 5 follows from \cite{char}, Proposition 2.2, which gives the expression for  $+_{\d,r}$ in terms of the formal group law.
 Assertion 7 follows from assertion 6.  Assertion 8 is clear.  Assertion 12 is an easy exercise. Assertion 11 follows from assertion 9. Assertion 14 follows from the formula
 $$\phi^r(x)_{|x=1,x'=...=x^{(r-1)}=0}=1+p^rx^{(r)}.$$
To prove the rest of the assertions one uses the following adaptation of a standard trick in Witt vector theory. Define the ``ghost map" 
$$w:G \times {\mathfrak g}\ra G \times {\mathfrak g},\ \ \ w(a_0,a_1)=(a_0,a_0^{(p)}+pa_1).$$
Then $w$ is a  homomorphism of monoids with identity where the source $G\times {\mathfrak g}$ is the monoid (actually group) equipped with multiplication $\circ$  and the target $G\times {\mathfrak g}$ is equipped with the structure of direct product of multiplicative monoids (where ${\mathfrak g}$ is viewed as a monoid with respect to multiplication in the associative algebra ${\mathfrak g}$). The homomorphism $w$ injective. So  in order to prove the identities in assertions 6, 9, 10 it is enough to prove them after composition with $w$; this however follows trivially using the equalities: $$w(\nabla^1(a))=(a,\phi(a)),\ \ w(\sigma(a))=(a,\Phi(a)).$$
\qed

\begin{remark}
Here are some facts about involutions on $GL_n$  that will be become relevant in what follows. Let $q\in GL_n$, $n\geq 2$. Then the following hold:

1) $x\mapsto q^{-1}(x^t)^{-1}q$ is an involution  if and only if $q^t=\pm q$.

2) $x\mapsto q^{-1}xq$ is an involution  if and only if $q=cu$ where $c$ is in the center of $GL_n$  and $u^2=1$.

3) $x\mapsto \det(x)q^{-1}(x^t)^{-1}q$ is an involution  if and only if $n=2$ and  $q^t=\pm q$.

4) $x\mapsto \det(x)^{-1}q^{-1}xq$ is an involution  if and only if $n=2$ and  $q=cu$, where $c$ is in the center of $GL_2$ and $u^2=1$.

5) Any involution of $GL_n$ over an algebraically closed field is of one of the forms 1), 2), 3), 4) above. This follows easily from the well known fact  \cite{dieudonne} that the  automorphism group of $SL_n$ over a algebraically closed field is generated by the inner automorphisms together with $x\mapsto (x^t)^{-1}$. 

Involutions of the form 1) will be referred to as {\it outer} involutions. Involutions of the form 2) will be referred to as {\it inner} involutions; more generally, of course, automorphisms of the form $x\mapsto q^{-1}xq$ are referred to as {\it inner} automorphisms. Outer involutions are  used to define the classical symmetric spaces. Inner automorphisms are, on the other hand, closely connected to conjugacy classes. Each of these two types of homogeneous spaces will be considered from the viewpoint of our theory in what follows. The involutions of the form 3) and 4) will not be considered in this paper but can be treated along the same lines as the other types 1) and 2).
\end{remark}

\begin{definition}
 Denote by $T\subset GL_n$  the split torus of diagonal matrices, by $N$ the normalizer of $T$ in $GL_n$,  and by $W\subset GL_n$  the group of permutation matrices. Then  $N=WT=TW$ and $W$ is isomorphic to the {\it Weyl group} $N/T$. In what follows we will view $T,N,W$ either as  abstract subgroups of the abstract group $GL_n$ or as subgroup schemes of the group scheme $GL_n$; all three group schemes are smooth over $R$.
 Finally one can consider the involution $x^{\tau}=x$ on $GL_n$;  the group defined by $\tau$ is, of course, $GL_n$ itself and we refer to $\tau$ as the canonical involution defining $GL_n$. Note that the associated maps $\cH$ and ${\mathcal B}$ are given by $\cH(x)=1$ and 
 ${\mathcal B}(x,y)=x^{-1}y$.
 Writing $G=GL_n$ we have
$$N=\{b\in G; (ab)^{(p)}=a^{(p)}b^{(p)}\ \ \text{for all}\ \ a\in G\},$$
$$N=\{b\in G;
(ba)^{(p)}=b^{(p)}a^{(p)}\ \ \text{for all}\ \ a\in G\}.$$
Also note that for any $w\in W$ we have $w^t=w^{-1}$, $w^{(p)}=w$, $\phi(w)=w$.
Also, for any $d\in T$, $d^t=d$ and $d^{(p)}=d^p$. A lift of Frobenius on $GL_n$ is (left, right, bi-) compatible with $N$ if and only if it is (left, right, bi-) compatible with both $T$ and $W$.
\end{definition}

In the Proposition below $\widehat{GL_n}$ is viewed as an open subset of $\widehat{{\mathfrak g}{\mathfrak l}_n}$.

\begin{proposition}
\label{characterization}
The lift of Frobenius $\phi_{GL_n,0}(x)=x^{(p)}$ on $\widehat{GL_n}$ is the unique lift of Frobenius  that is bicompatible with $N$  and extends to a lift of Frobenius
 on $\widehat{{\mathfrak g}{\mathfrak l}_n}$. \end{proposition}

{\it Proof}. Set $G=GL_n=Spec\ R[x,\det(x)^{-1}]$.
Let $T=Spec\ R[s_1,s_1^{-1},...,s_n,s_n^{-1}]$ and consider another copy of $T$ given by $T=Spec\ R[t_1,t_1^{-1},...,t_n,t_n^{-1}]$. The multiplication map $\mu:T\times G \times T\ra G$ is given by the map 
$$\mu:R[x,\det(x)^{-1}]\ra R[x,\det(x)^{-1},s_1,s_1^{-1},...,s_n,s_n^{-1},t_1,t_1^{-1},...,t_n,t_n^{-1}]$$
$$\mu x= sxt$$
where $s=\text{diag}(s_1,...,s_n)$, $t=\text{diag}(t_1,...,t_n)$.
Let $\phi_{GL_n}$ be a  lift of Frobenius on $\widehat{GL_n}$ that is bicompatible with $T$ and $W$ and extends to a lift of Frobenius
 on $\widehat{{\mathfrak g}{\mathfrak l}_n}$. Then $\phi_{GL_n}(x)=\Phi(x)$, $\Phi(x)=(\Phi_{ij}(x))$, $\Phi_{ij}(x)\in R[x]\h$. Bicompatibility with $T$ implies 
 $$\Phi(sxt)=s^p\Phi(x)t^p.$$
 For any $n\times n$ matrix $M=(m_{ij})$ with entries $m_{ij}\in \bZ_+$ (non-negative integers) we write $x^M=\prod_{ij} x_{ij}^{m_{ij}}$. We say that $x^M$ has weight $(u,v)
 \in \bZ_+^n\times \bZ_+^n$, $u=(u_1,...,u_n)$, $v=(v_1,...,v_n)$,  if
 $u_i=\sum_j m_{ij}$ for all $i$ and  $v_j=\sum_i m_{ij}$
 for all  $j$. An element in $R[x]\h$ is said to have weight $(u,v)$ if it is a (possibly infinite) sum of monomials $\lambda x^M$ with $\lambda\in R$ and $x^M$ of weight $(u,v)$.
 Write $\Phi_{ij}(x)=\sum_{u,v\in \bZ_+^n} \Phi_{ijuv}$ with $\Phi_{ijuv}$ of weight $(u,v)$.
Write $s^u=s_1^{u_1}...s_n^{u_n}$ and similarly for $t^v$. Then we have
  $$\mu \Phi_{ij}=\sum_{u,v} s^u t^v \Phi_{ijuv}.$$
  On the other hand
  $$\mu \Phi_{ij}=s_i^pt_j^p\Phi_{ij}=\sum_{u,v} s_i^pt_j^p \Phi_{ijuv}.$$
  This implies that $\Phi_{ij}$ has weight $(pe_i^t,pe_j^t)$ where $e_1,...,e_n$ are the columns of the identity matrix $1_n=[e_1,...,e_n]$. There is, however, only one monomial $x^M$ with weight $(pe_i^t,pe_j^t)$, namely $x^M=x_{ij}^p$. So $\Phi_{ij}=\lambda_{ij}x_{ij}^p$ for some $\lambda_{ij}\in R$. Now we use the bicompatibility with $W$ to show that all $\lambda_{ij}$ are equal; this will end the proof because the common value of the $\lambda_{ij}$s must be $1$ (which follows by looking at the restriction of $\phi_{GL_n}$ to $T$). To use bicompatibility with $W$ denote by $w_{\sigma}\in W$ the matrix $w_{\sigma}=[e_{\sigma(1)},...,e_{\sigma(n)}]$ obtained by permuting the columns of the identity matrix according to the permutation $\sigma$. For any matrix $a=(a_{ij})$ we have the formulae $(aw_{\sigma})_{ij}=a_{i\sigma(j)}$, $(w_{\sigma}a)_{\sigma(i)j}=a_{ij}$. Bicompatibility with $W$ yields 
  \begin{equation}
  \label{doggy}
  \Phi(xw_{\sigma})=\Phi(x)w_{\sigma}\ \ \ \text{and}\ \ \ \Phi(w_{\sigma}x)=w_{\sigma}\Phi(x).\end{equation}
  The first of these equalities yields $(\Phi(xw_{\sigma}))_{ij}=(\Phi(x)w_{\sigma})_{ij}$.
  Now
  $$(\Phi(xw_{\sigma}))_{ij}=\Phi_{ij}(xw_{\sigma})=\lambda_{ij}(xw_{\sigma})_{ij}=\lambda_{ij}x_{i\sigma(j)}^p;$$
  on the other hand
  $$(\Phi(x)w_{\sigma})_{ij}=\Phi_{i\sigma(j)}(x)=\lambda_{i\sigma(j)}x_{i\sigma(j)}^p.$$
  We conclude that $\lambda_{ij}=\lambda_{i\sigma(j)}$. Similarly, using the second of the equalities in \ref{doggy} we get $\lambda_{ij}=\lambda_{\sigma(i)j}$. This ends the proof.
\qed

\begin{remark}
Proposition \ref{characterization} fails if one removes the condition that the lift of Frobenius extend to one on $\widehat{{\mathfrak g}{\mathfrak l}_n}$. Indeed assume $n=2$ and consider $u(x)\in R[x,\det(x)^{-1}]\h$  defined by
$$u(x)=\frac{x_{11}x_{12}x_{21}x_{22}}{\det(x)^2}.$$
Note that $u(sxt)=u(x)$ for all $s,t\in T$. 
Define the lift of Frobenius $\phi_{GL_2}$ on $\widehat{GL_2}$ by 
$\phi_{GL_2}(x)=(\Phi_{ij}(x))$ where $\Phi_{ij}(x)=u(x) \cdot x^p_{ij}$. Then clearly $\phi_{GL_2}$ is bicompatible with $T$ and $W$. 
\end{remark}

\medskip

For the rest of the paper we let $\phi_{GL_n,0}$ be the lift of Frobenius on $\widehat{GL_n}$  given by $\phi_{GL_n,0}(x)=x^{(p)}$. 

\begin{proposition}
\label{hepatitis}  Let $x^{\tau}=x$ on $GL_n$. The following hold.

i) $\phi_{GL_n,0}$ is the unique lift of Frobenius that is ${\mathcal H}$-horizontal  and ${\mathcal B}$-symmetric  with respect to $\phi_{GL_n,0}$.

ii) The arithmetic logarithmic derivative 
$l\d:GL_n\ra {\mathfrak g}{\mathfrak l}_n$ attached to $\phi_{GL_n,0}$ is given by
$$l\d a:=\d a \cdot (a^{(p)})^{-1}.$$

iii) The map
$\{\ ,\ \}:GL_n\times GL_n\ra {\mathfrak g}{\mathfrak l}_n$ attached to $\phi_{GL_n,0}$
is given by
$$\{a,b\}=p^{-1}(a^{(p)}\cdot b^{(p)} \cdot ((ab)^{(p)})^{-1}-1).$$

iv) If $\epsilon=1+p\alpha$ then the  equation 
$l\d (u)=\alpha$ is equivalent to either of the following equations:
$$\begin{array}{rcl}
\d u & = & \alpha\cdot u^{(p)},\\
\phi(u) & = & \epsilon \cdot u^{(p)}.\end{array}$$

v) For all $a,b \in GL_n$ with either $a$ or $b$ in $N$ we have
  $$
\label{cris}l\d(ab)=(\phi(a)\cdot l\d(b) \cdot \phi(a)^{-1}) +_{\d} l\d(a).$$
\end{proposition}

{\it Proof}.  Uniqueness in i)  follows from Corollary \ref{unique}.
ii), iii), iv) are a specialization of Proposition \ref{hair}.
v) follows from Lemma \ref{deal}\qed

\begin{remark}
\label{iridor}
For a $p$-adically complete ring ${\mathcal A}$
we  will need to use, in what follows, various fractional powers of matrices with entries in the ring ${\mathcal A}$. So we introduce the following convention. Let $U\in GL_M({\mathcal A})$ with $U\equiv 1=1_M$ mod $p$,
and let $\nu\in \bZ$ be an integer not divisible by $p$. Then let $V=\frac{U-1}{p}$ and set
\begin{equation}
\label{rirdor}
U^{1/\nu}:=(1+pV)^{1/\nu}:=\sum_{m=0}^{\infty}\left(\begin{array}{c} 1/\nu\\m\end{array}\right) p^m V^m\in GL_M({\mathcal A}).\end{equation} 
We refer to the above $U^{1/\nu}$ as the {\it principal} $1/\nu$-root of $U$. As a consequence we get:
\end{remark}

\begin{corollary}
Let $G=GL_n$. Then, for any integer $\nu\in \bZ$ coprime to $p$,
the map $L_{\d}(G)\ra L_{\d}(G)$, $b\mapsto b^{\nu}$ is bijective.
\end{corollary}

{\it Proof}.
Injectivity follows from Corollary \ref{wolf}. Surjectivity follows from Remark \ref{iridor}.
\qed

\begin{remark}\label{kill}
Let $G\subset GL_n$ be a smooth closed subgroup scheme and let $\phi_{GL_n}$ be a lift of Frobenius on $\cO(GL_n)\h=R[x,\det(x)^{-1}]\h$ such that $G$ is $\d_{GL_n}$-horizontal. (Such a $\phi_{GL_n}$ always exists due to the smoothness of $G$ and Lemma \ref{liftingg}.) Let $\phi_G$ be the induced lift of Frobenius on $\cO(G)\h$.
Then the embedding $G\subset GL_n$ induces a closed embedding  $J^1(G)\subset J^1(GL_n)$ and hence a closed embedding $L_{\d}(G)\subset L_{\d}(GL_n)$. 
Composing the latter with the natural identification $L_{\d}(GL_n)= {\mathfrak g}{\mathfrak l}_n$
we get a natural embedding of groups $L_{\d}(G)\subset {\mathfrak g}{\mathfrak l}_n$ where ${\mathfrak g}{\mathfrak l}_n$ is viewed as a group with respect to $+_{\d}$. Note that the usual Lie algebra $L(G)$ is also embedded as  a subgroup  into $L(GL_n)={\mathfrak g}{\mathfrak l}_n$ where this time ${\mathfrak g}{\mathfrak l}_n$ is equipped with the usual addition $+$ of matrices; but note that, as subsets of ${\mathfrak g}{\mathfrak l}_n$ we have $L_{\d}(G)\neq L(G)$ in general! Indeed it is trivial to check that if 
$$G=Spec\ R[x,\det(x)^{-1}]\h/(f_1(x),...,f_m(x))$$
and if
$$\d_1f_i(x'):=p^{-1}f_i^{(\phi)}(1+px')$$ 
where $f^{(\phi)}_i$ is obtained from $f_i$ by applying $\phi$ to the coefficients then
$$L_{\d}(G)=Spf\ R[x']\h/(\d_1f_1(x'),...,\d_1f_m(x'))$$
hence the $R$-points of the latter are the zeroes in ${\mathfrak g}{\mathfrak l}_n$ of the polynomials $\d_1 f_i(x')$. On the other hand  the Lie algebra $L(G)$
is given by 
$$L(G)=Spec\ R[x']/(d_1f_1(x'),...,d_1f_m(x'))$$
where
$$d_1f_i(x')=``\epsilon^{-1}"f_i(1+\epsilon x').$$
So the set of $R$-points of $L(G)$ is the zero locus in ${\mathfrak g}{\mathfrak l}_n$ of the polynomials $d_1f_i$.
 Cf. our discussion of examples in the next section.
Going back to our discussion  it follows from  horizontality that
$$l\d:G\ra L_{\d}(G)$$
is the restriction of the map
$$l\d:GL_n\ra L_{\d}(GL_n)={\mathfrak g}{\mathfrak l}_n.$$
Similarly the map
$$\{\ ,\ \}_{G}:G\times G\ra L_{\d}(G)$$
 is the restriction of the map
 $$\{\ ,\ \}_{GL_n}:GL_n\times GL_n\ra L_{\d}(GL_n)={\mathfrak g}{\mathfrak l}_n.$$
 In particular  the formulas involving $l\d$ and $\{\ ,\ \}$ in Proposition \ref{hair} continue to be valid if one replaces $G=GL_n$  by $G\subset GL_n$ an arbitrary smooth closed subgroup of $GL_n$, provided $G$ is $\d_{GL_n}$-horizontal and provided we make the identifications explained above.
 
 Finally note that under the identifications above we have $$L_{\d}^r(GL_n)=Spf\ R[x^{(r)}]\h,$$ 
 $$L_{\d}^r(G)=Spf\ R[x^{(r)}]\h/(p^{-r}f_1^{(\phi^r)}(1+p^rx^{(r)}),...,p^{-r}f_m^{(\phi^r)}(1+p^rx^{(r)})).$$
 In particular the map $ex^r:L_{\d}^r(GL_n)\ra GL_n$, $ex^r(a)=1+p^r\phi^{-r}(a)$, satisfies:
 \begin{equation}
 \label{director}
 (ex^r)^{-1}(G)=L_{\d}^r(G).\end{equation}
 \end{remark}
 
 We are ready to give the
 
 \medskip
 
 {\it Proof of Proposition \ref{bracketprop}.}
 For any $\alpha,\beta\in {\mathfrak g}{\mathfrak l}_n$ the right hand side of \ref{snuggles} is trivially seen to be a restricted power series in $\phi^s(\alpha)$ and $\phi^r(\beta)$. Hence assertion 3 in our proposition follows for $G=GL_n$. Also we have an induced $\d$-map 
 \begin{equation}
 \label{om}
 [\ ,\ ]_{\d}:L_{\d}^r(GL_n)\times L_{\d}^s(GL_n)\ra L_{\d}^{r+s}(GL_n)
 \end{equation}
 such that the following diagram of sets is commutative:
  \begin{equation}
 \label{telephones}
 \begin{array}{rcl}
 L_{\d}^{r}(GL_n)\times L_{\d}^{s}(GL_n) & \stackrel{[\ ,\ ]_{\d}}{\longrightarrow} &
  L_{\d}^{r+s}(GL_n)\\
  ex^{r}\times ex^{s} \downarrow & \  & \downarrow ex^{r+s}\\
  GL_n \times GL_n & \stackrel{[\ ,\ ]}{\longrightarrow} & GL_n
 \end{array}
 \end{equation}
 By \ref{director} the upper horizontal arrow in \ref{telephones} sends $L_{\d}^r(G)\times L_{\d}^s(G)$ into $L_{\d}^{r+s}(G)$ so the $\d$-map \ref{om} induces a $\d$-map as in assertion 1 of the proposition. (Uniqueness follows from the injectivity of the $ex$ maps.)
 Assertion 2 of the proposition is clear from the functoriality of $L_{\d}^r(G)$ and $ex^r$.
 Assertion 7 is a direct computation. Now to prove assertions 3, 4, 5, 6 it is enough to prove them in the case $G=GL_n$.
 We already know assertion 3 holds in this case. Assertion 5 for $G=GL_n$ follows directly from assertion 7.
 Assertion 6 follows from assertions 3 and 7 because of the usual Jacobi identity plus the following computation:
 $$
 \begin{array}{rcll}
 [[\alpha,\beta]_{\d},\gamma]_{\d} & \equiv & [[\alpha^{(p^s)},\beta^{(p^r)}],\gamma]_{\d}& \text{mod}\ \ p\\
 \  & \  & \  \\
 \  & \equiv &  [[\alpha^{(p^s)},\beta^{(p^r)}]^{(p^t)},\gamma^{(p^{r+s})}]& \text{mod}\ \ p\\
 \  & \  & \  \\
 \  & \equiv &  [[\alpha^{(p^{s+t})},\beta^{(p^{r+t})}],\gamma^{(p^{r+s})}]& \text{mod}\ \ p.
 \end{array}
 $$
 Assertion 4 can be checked in the same way.
 \qed
 
 \section{Outer involutions}
 
 In this section we consider  outer involutions on $GL_n$ and prove our main results about them. We start with the case of $SL_n$; we will then continue by investigating the case of $SO_n$ and $Sp_n$ (or more generally the groups $SO(q)$ attached to symmetric/antisymmetric matrices $q$; see below).

\begin{example}
\label{SL} ($SL_n$).
Consider the group scheme
$$SL_n=Spec\ R[x]/(\det(x)-1).$$
Embed $SL_n$ as a closed subgroup scheme of $GL_n$ in the natural way and equip $\widehat{GL_n}$ with the lift of Frobenius $\phi_{GL_n,0}(x)=x^{(p)}$. 
Recall that $T,N,W$ are the torus of diagonal matrices in $GL_n$, its normalizer in $GL_n$,  and the subgroup of the permutation matrices in $GL_n$.  Let $T_{SL_n}:=T\cap SL_n$ and $N_{SL_n}:=N\cap SL_n$. 
 Then $T_{SL_n}$ is a maximal torus in $SL_n$ and $N_{SL_n}$ is the normalizer of $T_{SL_n}$ in $SL_n$. The groups $T_{SL_n}$ and $N_{SL_n}$ are smooth. Moreover the Weyl group $W_{SL_n}:=N_{SL_n}/T_{SL_n}$ is isomorphic to $W$ (but the map $W\cap SL_n\ra W_{SL_n}$ is not an isomorphism!)  Also recall that the group of characters $Hom(T,{\mathbb G}_m)$ of $T$ has a $\bZ$-basis $\chi_1,...,\chi_n$ where 
 $$\chi_i(\text{diag}(d_1,...,d_n))=d_i.$$
 The restrictions of $\chi_1,...,\chi_n$ to $T_{SL_n}$, which will still be denoted by $\chi_1,...,\chi_n$, generate the group of characters of $T_{SL_n}$ and satisfy one relation $\chi_1+...+\chi_n=0$ (where we use, as usual, the additive notation for characters). The {\it roots} of $SL_n$ are
 $$\chi_i-\chi_j,\ \ i\neq j.$$
 The {\it root subgroup}  $U_{\chi_i-\chi_j}$  comes from an embedding of ${\mathbb G}_a$ over $R$ into $SL_n$ and, as an abstract group
 $$U_{\chi_i-\chi_j}=\{\mu E_{ij}\}$$
 where $\mu\in R$ and $E_{ij}$ is the matrix with $1$s on the diagonal, $1$ on position $(i,j)$ and $0$ everywhere else. The embedding ${\mathbb G}_a\ra SL_n$ is given by $\mu\mapsto \mu E_{ij}$. Note that $U_{\chi_i-\chi_j}$ is $\phi_{GL_n,0}$-horizontal in $GL_n$ and the lift of Frobenius induced by $\phi_{GL_n,0}$ on $\widehat{{\mathbb G}_a}=Spf\ R[z]\h$  is given by $z\mapsto z^p$.
 
 Finally recall from the Introduction that,
 in discussing  $SL_n$, we need to require $p \not| n$ and we need to 
 consider the \'{e}tale cover $\pi:G=GL_n'\ra GL_n$ defined in \ref{scorpion1}. On $G$ we have the involution $\tau$ in \ref{scorpion2} and hence the quadratic map ${\mathcal H}$ and the bilinear map ${\mathcal B}$ attached to $\tau$.   Then $\pi$ induces an isomorphism $(G^+)^{\circ}\simeq SL_n$; we call $\tau$ the {\it canonical involution} defining $SL_n$.   Of course, $SL_n$ is smooth over $R$.
 Let us equip $\widehat{G}$ with the unique lift of Frobenius $\phi_{G,0}$ extending the lift of Frobenius $\phi_{GL_n,0}$ on $\widehat{GL_n}$.
 Consider a lift of Frobenius $\phi_{GL_n}$ on $\widehat{GL_n}$. We say that
 $\phi_{GL_n}$ is {\it ${\mathcal H}$-horizontal} (respectively {\it ${\mathcal B}$-symmetric}) with respect to  $\phi_{GL_n,0}$ 
 if the unique lift of Frobenius $\phi_G$ on $\widehat{G}$ extending $\phi_{GL_n}$ is ${\mathcal H}$-horizontal  (respectively ${\mathcal B}$-symmetric) with respect to    $\phi_{G,0}$.
 For any such $\phi_{GL_n}$ the group $SL_n$ is easily seen to be $\phi_{GL_n}$-horizontal so there is an induced lift of Frobenius $\phi_{SL_n}$ on $\widehat{SL_n}$.
\end{example}

In the next statement fractional powers are principal roots as in \ref{rirdor}.

\begin{proposition}\label{night}
\ 

i) There is a unique lift of Frobenius
 $\phi_{GL_n}$ on $\widehat{GL_n}$ 
 that is ${\mathcal H}$-horizontal  and ${\mathcal B}$-symmetric with respect to $\phi_{GL_n,0}$.
  It is given by
  $\phi_{GL_n}(x)=  \lambda(x)\cdot x^{(p)}$ where $\lambda(x)\in R[x,\det(x)^{-1}]\h$,
   $$
\lambda(x) = \left( \frac{\det(x^{(p)})}{(\det(x))^p}\right)^{-1/n}.$$

ii)    $\phi_{GL_n}$ is 
bicompatible with 
 $N$.
 
  iii)  $\phi_{GL_n}$ and  $\phi_{GL_n,0}$ coincide on the set $SL_n\cap\phi_{GL_n,0}^{-1}(SL_n)$; in particular the root subgroups $U_{\chi_i-\chi_j}$ are $\phi_{GL_n}$-horizontal and the  lifts of Frobenius induced by $\phi_{GL_n}$ on $\widehat{{\mathbb G}_a}=Spf\ R[z]\h$  are given by $z\mapsto z^p$.

 iv) The arithmetic logarithmic derivative
$l\d:GL_n\ra   {\mathfrak g}{\mathfrak l}_n$  attached to $\phi_{GL_n}$
 is given by
$$l\d a=\d a \cdot (\lambda(a)\cdot a^{(p)})^{-1}+\frac{\lambda(a)^{-1}-1}{p}\cdot 1_n$$ 
and it satisfies 
$$
\label{cris}l\d(ab)=(\phi(a)\cdot l\d(b) \cdot \phi(a)^{-1}) +_{\d} l\d(a)$$
 for all $a, b\in GL_n$ with either $a$ or $b$ in $N$.

v) The map
$\{\ ,\ \}:GL_n\times GL_n\ra  {\mathfrak g}{\mathfrak l}_n$
is given by
$$\{a,b\}=p^{-1}\left(\frac{\lambda(a)\lambda(b)}{\lambda(ab)}a^{(p)}b^{(p)}((ab)^{(p)})^{-1}-1_n\right).$$

vi) If  $\epsilon=1+p\alpha$ then the  equation 
$l\d (u)=\alpha$ is equivalent to either of the following equations:
$$\begin{array}{rcl}
\d u & = & \left(\lambda(u)\cdot \alpha+\frac{\lambda(u)-1}{p}\cdot 1_n\right) \cdot u^{(p)},\\
\  & \  & \  \\
\phi(u) & = & \epsilon \cdot \lambda(u) \cdot  u^{(p)}.\end{array}$$

 vii)  Let $\alpha\in L_{\d}(SL_n)$, let  $\phi^{\alpha}_{GL_n}$ be the lift of Frobenius
on  $\widehat{GL_n}$ defined by $\phi^{\alpha}_{GL_n}(x)=\epsilon\cdot \phi_{GL_n}(x)$, $\epsilon=1+p\alpha$, and let $\d^{\alpha}_{GL_n}$ be the $p$-derivation on $\cO(GL_n)\h$ associated to $\phi^{\alpha}_{GL_n}$. Let $\cH^*\in \cO(GL_n)\h$ be defined as $\cH^*(x)=\det(x)$. Then
$$\d_{GL_n}^{\alpha}(\cH^*)=0.$$
\end{proposition}

 {\it Proof}. 
 We have
$$\cO(G)=R[x,\det(x)^{-1},y]/(y^n-\det(x)^2)=R[x,\det(x)^{-1},t],$$
where $t$ is  the class of $y$. Denoting, as usual, by $\tau,\cH,{\mathcal B}$ the ring maps induced by the corresponding scheme maps, we have
$x^{\tau}=t^{-1}x$, $t^{\tau}=t^{-1}$, ${\mathcal B}(x)=t_1x^{-1}_1x_2$, $\cH(x)=t\cdot 1_n$. (Here $t_1,x_1,x_2$ are $t\otimes 1$, $x\otimes 1$ and $1\otimes x$ respectively.)
Now from $t^n=\det(x)^2$ we get
\begin{equation}
\label{formulae2}
\begin{array}{rcl}
\phi_{G,0}(t)^n & = & \phi_{G,0}(\det(x))^2\\
\  & = & \det(x^{(p)})^2\\
\  & = & \lambda(x)^{-2n} \det(x)^{2p}\\
\  & = & (\lambda(x)^{-2} t^p)^n\end{array}
\end{equation}
where $\lambda(x)$ is defined as in assertion i).
Since $\phi_{G,0}(t)\equiv \lambda(x)^{-2} t^p$ mod $p$ it follows that
\begin{equation}
\label{formulae1}
\phi_{G,0}(t)=\lambda(x)^{-2} t^p.\end{equation}
 Let now $\phi_{GL_n}(x)=\lambda(x)\cdot  x^{(p)}$.
Then 
\begin{equation}
\begin{array}{rcl}
\phi_G(t)^n & = & \phi_G(\det(x))^2\\
\  & = & \det(\lambda(x) \cdot x^{(p)})^2\\
\  & = & \lambda(x)^{2n}\det(x^{(p)})^2 \\
\  & = & (\det(x))^{2p}\\
\  & = &  t^{np}.\end{array}
\end{equation}
Since $\phi_G(t)\equiv  t^p$ mod $p$ it follows that
$\phi_G(t)=   t^p$.
Consequently
$$\cH(\phi_{G,0}(x))=\cH(x^{(p)})= t^p\cdot 1.$$
On the other hand
$$
\phi_G(\cH(x))=\phi_G(t\cdot 1)= t^p\cdot 1.$$
This proves  the commutativity of \ref{sleepless}. 
To check the commutativity of \ref{seatle} note that
$$(\phi_{G,0}\times \phi_G)({\mathcal B}(x))=(\phi_{G,0}\times \phi_G)(t_1\cdot x_1^{-1}\cdot x_2)=
(\lambda(x)^{-2}t^p)\cdot (x^{(p)})^{-1}\cdot (\lambda(x)x^{(p)}),$$
$$(\phi_{G}\times \phi_{G,0})({\mathcal B}(x))=(\phi_{G}\times \phi_{G,0})(t_1\cdot x_1^{-1}\cdot x_2)=
t^p \cdot (\lambda(x) x^{(p)})^{-1}\cdot x^{(p)};$$
the right hand sides of the two equations above are, of course, equal.
 So the existence part of assertion i) is proved. The uniqueness follows 
from Corollary \ref{unique}. 

To check assertion ii)  note first that if $a\in N$   then
$$\begin{array}{rcl}
\lambda(xa) & = &\left( \frac{\det((xa)^{(p)})}{(\det(xa))^p}\right)^{-1/n}\\
\  & = &
\left( \frac{\det(x^{(p)}a^{(p)})}{\det(x)^p\cdot \det(a)^p}\right)^{-1/n}\\
\  & = & 
\left( \frac{\det(x^{(p)})\cdot \det(a)^p}{\det(x)^p\cdot \det(a)^p}\right)^{-1/n}\\
\  & = & \lambda(x).\end{array}
$$
Similarly $\lambda(ax)=\lambda(x)$.  Hence 
$$\phi_{GL_n}(xa)=\lambda(xa)\cdot (xa)^{(p)}=\lambda(x)\cdot x^{(p)}\cdot a^{(p)} =\phi_{GL_n}(x)\cdot a^{(p)}.$$
Similarly 
$$\phi_{GL_n}(ax)=\lambda(ax)\cdot (ax)^{(p)}=\lambda(x)\cdot a^{(p)}\cdot x^{(p)}=a^{(p)}\cdot \phi_{GL_n}(x).$$
In particular $\phi_{GL_n}(a)=\lambda(a)\cdot a^{(p)}=a^{(p)}$ which is in $N$ if $a$ is in $N$. So
$N$ is $\d_{GL_n}$-horizontal and $\phi_{GL_n}$ is bicompatible with $N$.

To check assertion iii) let $g\in SL_n\cap \phi_{GL_n,0}^{-1}(SL_n)$.
So $\det(g)=1$ and 
$$1=\det(\phi_{GL_n,0}(g))=\det(\phi^{-1}(g^{(p)}))=\phi^{-1}(\det(g^{(p)}))$$ hence $\det(g^{(p)})=1$.
So $\lambda(g)=1$ and so
$$\phi_{GL_n}(g)=\phi^{-1}(\lambda(g))\phi^{-1}(x^{(p)})=\phi^{-1}(\lambda(g))\phi_{GL_n,0}(g)=\phi_{GL_n,0}(g).$$
which ends the proof of the first assertion in iii). The second assertion in iii) follows from Lemma \ref{pat}.

Assertions iv), v), vi) follow from Proposition \ref{hair} and Lemma \ref{deal}.

To check assertion vii) note that the equality $\d_{GL_n}^{\alpha}(\det(x))=0$, which is equivalent to $\phi_{GL_n}^{\alpha}(\det(x))=\det(x)^p$, follows from the following computation (where we use $\det(\epsilon)=1$):
$$
\begin{array}{rcl}
\phi_{GL_n}^{\alpha}(\det(x)) & = & \det(\phi_{GL_n}^{\alpha}(x))\\
\  & = & \det(\lambda(x) \cdot \epsilon \cdot x^{(p)})\\
\  & = & \lambda(x)^n\cdot \det(\epsilon) \cdot \det(x^{(p)})\\
\  & = & \det(x)^p.
\end{array}
$$
\qed

\medskip
\begin{example}
\label{O} ($SO(q)$).
Let again $G=GL_n$ and consider the lift of Frobenius $\phi_{GL_n,0}(x)=x^{(p)}$ on $\widehat{G}$.
 Let $x^t$ be the transpose of $x$ and let $q\in GL_n$ be either symmetric or antisymmetric, i.e. $q^t=\pm q$. Then the formula $x^{\tau}=q^{-1}(x^t)^{-1}q$ defines an involution on $GL_n$. The associated quadratic map ${\mathcal H}$, normalized quadratic map $\cH_q$,  bilinear map ${\mathcal B}$, and normalized bilinear map ${\mathcal B}_q$ are 
 $$\cH(x)=q^{-1}x^tqx,\ \ \ \cH_q(x)=x^tqx,\ \ \ {\mathcal B}(x,y)= q^{-1}x^tqy,\ \ {\mathcal B}_q(x,y)=x^tqy.$$
  Let $SO(q)=(G^+)^{\circ}$ be the group defined by $\tau$; we refer to $\tau$ as the {\it canonical involution} defining $SO(q)$. The group $SO(q)$ is trivially seen to be  smooth  over $R$.
 Also $G^-=q^{-1}G_{\pm}$ where $G_{\pm}=\{x;x^t=\pm x\}$;  so $G^-$ is smooth (an open set of an affine space). 
 If $T\subset GL_n$ is the torus of diagonal matrices in $GL_n$ and $N\subset GL_n$ is the normalizer of $T$ in $GL_n$ we set $T_{SO(q)}:=T\cap SO(q)$ and $N_{SO(q)}:=N\cap SO(q)$. 
 It is easy to see that $T_{SO(q)}$ and $N_{SO(q)}$ are smooth over $R$. 
  Note that the $\d$-Lie algebra $L_{\d}(SO(q))$ identifies with the set of matrices
$a\in {\mathfrak g}{\mathfrak l}_n$ satisfying 
$$a^t\phi(q)+\phi(q)a+pa^t\phi(q)a=0.$$
On the other hand the Lie algebra $L(SO(q))$ identifies with the set of matrices
$a\in {\mathfrak g}{\mathfrak l}_n$ satisfying 
$$a^tq+qa=0.$$

We say $q$ (or $SO(q)$) is {\it split} if $q$ is one of the matrices \ref{scorpion3}. If $q$ is split then   $qq^t=1$, $q^t=\pm q$, $q\in N$,  $q^{(p)}=q$, $\phi(q)=q$, $\d q=0$, and, for all $b\in G$,  $(qb)^{(p)}=qb^{(p)}$ and $(bq)^{(p)}=b^{(p)}q$. We claim that $G^-$ is $\phi_{GL_n,0}$-horizontal; indeed, since $G^-$ is smooth it is enough to check that if $a\in G^-$ then $a^{(p)}\in G^-$ which is a trivial exercise.

Let us assume, until further notice, that $q$ is split. In this case
 $SO(q)$ is  denoted by
  $Sp_{2r},\ \ SO_{2r},\ \ SO_{2r+1}$ respectively.  The groups $T_{SO(q)}$ are then the {\it  split maximal tori}:
 $$\begin{array}{lcl}
 T_{Sp_{2r}}& = & \{d\in T; d=\text{diag}(d_1,...,d_r,d_1^{-1},...,d_r^{-1})\}\\
 T_{SO_{2r}} & = & \{d\in T; d=\text{diag}(d_1,...,d_r,d_1^{-1},...,d_r^{-1})\}\\
 T_{SO_{2r+1}} & = & \{d\in T; d=\text{diag}(1,d_1,...,d_r,d_1^{-1},...,d_r^{-1})\}.
 \end{array}
 $$
 One also has that $N_{SO(q)}$ is the normalizer of $T_{SO(q)}$ in $SO(q)$; this follows from the fact that if $a\in SO(q)$ normalizes $T_{SO(q)}$ and $d\in T_{SO(q)}$ is a regular matrix (i.e. has distinct mod $p$ elements on the diagonal) then $a^{-1}da=wdw^{-1}$ for some $w\in W$ hence $aw$ is in the centralizer of $d$; by regularity of $d$ we get $aw\in T$, hence $a\in N$, hence $a\in N_{SO(q)}$.
 
 Note now that
  $N_{SO(q)}$ contains all the matrices of the form
 \begin{equation}
 \label{gorilla}\left( 
\begin{array}{ll} w_r & 0\\0 & w_r\end{array}\right),\ \ \left( 
\begin{array}{ll} w_r & 0\\0 & w_r\end{array}\right),\ \ 
\left( \begin{array}{lll} 1 & 0 & 0\\
0 & w_r & 0\\
0 & 0 & w_r\end{array}\right),
\end{equation}
respectively, where $w_r\in W_r$ is an arbitrary permutation matrix in $GL_r$.
The matrices \ref{gorilla} modulo $T_{SO(q)}$ do not exhaust the {\it Weyl group} $W_{SO(q)}:=N_{SO(q)}/T_{SO(q)}$. A system of representatives of the  Weyl group can be obtained by adding to the matrices in \ref{gorilla} the matrices in the set $W'$ which we describe below. Let $e_1,...,e_n$ be the columns of the identity matrix $1_n$. Then, for $Sp_{2r}$, $W'$ consists of all products
$w'_{i_1}...w'_{i_s}$ where $s\geq 1$ and 
$$w'_i:=[e_1,...,e_{i-1},e_{r+i},e_{i+1},...,e_{r+i-1},-e_i,e_{r+i+1},...].$$
For $SO_{2r}$, $W'$ consists of all products
$w'_{i_1}...w'_{i_s}$ where $s\geq 2$ is even and 
$$w'_i:=[e_1,...,e_{i-1},e_{r+i},e_{i+1},...,e_{r+i-1},e_i,e_{r+i+1},...].$$
For $SO_{2r+1}$, $W'$ consists of all products
$w'_{i_1}...w'_{i_s}$ where $s\geq 1$  and 
$$w'_i:=[-e_1,e_2,...,e_{i},e_{r+i+1},e_{i+2},...,e_{r+i},e_{i+1},e_{r+i+2},...].$$
(Cf. \cite{GW}, section 2.5.1; the forms of the classical groups used in loc. cit. are different from ours; however, in order to conclude our claims,  the only information needed from loc. cit. is the cardinality of the Weyl groups.) The bottom-line of the above description of the Weyl groups is the following:

\begin{lemma}
\label{monkey}
Assume $q$ is split. Then for any $v\in N_{SO(q)}$ we have $v^{(p)}\in N_{SO(q)}$.
\end{lemma}

{\it Proof}.
Indeed any element $w'_i$ in the description above can be written as $w'_i=rw$ with $r\in T$, $r^2=1$, $w\in W$. So 
$$(w'_i)^{(p)}=(rw)^{(p)}=r^{(p)}w^{(p)}=rw=w'_i.$$
Now if $v\in N_{SO(q)}$ then either $v=tw$ or $v=tw'_{i_1}...w'_{i_s}$ with $t\in T_{SO(q)}$, $w\in W$, hence either 
$$v^{(p)}=t^{(p)}w^{(p)}=t^pw=t^{p-1}v\in N_{SO(q)}$$
or
$$v^{(p)}=t^{(p)}(w'_{i_1})^{(p)}...(w'_{i_s})^{(p)}=t^pw'_{i_1}...w'_{i_s}=t^{p-1}v\in N_{SO(q)}.$$
\qed

In what follows we review characters and roots of the split classical groups.
The groups of characters of $T_{Sp_{2r}}$ and $T_{SO_{2r}}$ have bases $\chi_1,...,\chi_r$ where
$$\chi_i(\text{diag}(d_1,...,d_r,d_1^{-1},...,d_r^{-1}))=d_i.$$
The group of characters of $T_{SO_{2r+1}}$ has basis $\chi_1,...,\chi_r$ where
$$\chi_i(\text{diag}(1,d_1,...,d_r,d_1^{-1},...,d_r^{-1}))=d_i.$$
The {\it roots} of $Sp_{2r}$ are $\pm 2\chi_i$, $\pm \chi_i\pm \chi_j$, $i\neq j$. The {\it roots} of $SO_{2r}$ are  $\pm \chi_i\pm \chi_j$, $i\neq j$. The {\it roots} of $SO_{2r+1}$ are $\pm \chi_i$, $\pm \chi_i\pm \chi_j$, $i\neq j$.
Recall that if $G$ is any of the groups $Sp_{2r},SO_{2r},SO_{2r+1}$ and if $\chi$ is a root of $G$ then there is an attached {\it root subgroup}  $U_{\chi}$ of $G$ over the algebraic closure $K^a$ of the field of fractions $K$ of $R$; it is the unique  subgroup isomorphic to ${\mathbb G}_a$ over $K^a$ normalized by the corresponding maximal torus on which this maximal torus acts via $\chi$. In all these cases $U_{\chi}$ comes from an embedding of ${\mathbb G}_a$ into $G$ over $R$ and has the following explicit description which we need to review. 
We consider $r\times r$ matrices as follows. For $i=1,...,r$ let $F_{ii}$ be the matrix with $1$ on position $(i,i)$ and $0$ everywhere else. Let now $i\neq j$ between $1$ and $r$. Let $E_{ij}$ be the matrix with $1$s on the diagonal, $1$ on position $(i,j)$ and $0$ everywhere else. Let $F_{ij}=E_{ij}+E_{ji}$ and $G_{ij}=E_{ij}-E_{ji}$. Finally let $e_i$ be the columns of the identity matrix $1_r$. Then the groups $U_{\chi}$ are given as follows (where $\mu \in R$). 

For $Sp_{2r}$ we have:
$$\begin{array}{rcl}
U_{2\chi_i} & = & \{\left( \begin{array}{cc}
1 & \mu F_{ii}\\
0 & 1\end{array}\right)\},\ \ U_{-2\chi_i}=U_{2\chi_i}^t\\
\  & \  &\ \\
U_{\chi_i-\chi_j} & = &  \{\left( \begin{array}{cc}
1+\mu E_{ij} & 0\\
0 & 1-\mu E_{ji}\end{array}\right)\}\\
\  & \  & \  \\
U_{\chi_i+\chi_j} & = & \{\left( \begin{array}{cc}
1 & \mu F_{ij}\\
0 & 1\end{array}\right)\},\ \ U_{-\chi_i-\chi_j}=U_{\chi_i+\chi_j}^t
\end{array}$$

For $SO_{2r}$ we have:
$$\begin{array}{rcl}
U_{\chi_i+\chi_j} & = & \{\left( \begin{array}{cc}
1 & \mu G_{ij}\\
0 & 1\end{array}\right)\},\ \ U_{-\chi_i-\chi_j}=U_{\chi_i+\chi_j}^t\\
\  & \  &\ \\
U_{\chi_i-\chi_j} & = &  \{\left( \begin{array}{cc}
1+\mu E_{ij} & 0\\
0 & 1-\mu E_{ji}\end{array}\right)\}
\end{array}$$

For $SO_{2r+1}$ we have:
$$\begin{array}{rcl}
U_{\pm\chi_i\pm \chi_j} & =  & \{\left( \begin{array}{cc}
1 & 0\\
0 & \star \end{array}\right)\},\ \ \text{where $\star$ is as for $SO_{2r}$}\\
\  & \  & \  \\
U_{\chi_i} & = &  \{\left( \begin{array}{ccc}
1 & \mu e_i^t & 0\\
0 & 1& 0\\
-\mu e_i & -\frac{\mu^2}{2}e_i e_i^t & 1
\end{array}\right)\},\ \ U_{-\chi_i}=U_{\chi_i}^t.
\end{array}$$
Remark that if $G$ is $Sp_{2r}$ and $\chi$ is any root or if $G$ is $SO_{2r}$ and $\chi$ is any root or if $G$ is $SO_{2r+1}$ and $\chi\neq \pm \chi_i$ we have that $U_{\chi}$ is $\phi_{GL_n,0}$-horizontal. On the other hand if $G=SO_{2r+1}$ and $\chi=\pm\chi_i$ then $U_{\chi}$ is not $\phi_{GL_n,0}$-horizontal (because it is not stable under the map $u\mapsto u^{(p)}$). For obvious reasons we shall refer to $\pm\chi_i$ as the {\it short roots}
of $SO_{2r+1}$. A root $\chi$ of a split $SO(q)$ will be called {\it abnormal} if $SO(q)=SO_{2r+1}$ and $\chi$ is a short root of the latter. 
\end{example}

\begin{proposition}
\label{cartan}
Assume $q$ is split and consider the involution
$$x^{\tau}=q^{-1}(x^t)^{-1}q$$ on $G=GL_n$ and the induced involution $b\mapsto b^{L_{\d}(\tau)}$ on $L_{\d}(G)={\mathfrak g}{\mathfrak l}_n$. Then the following hold:

1) For all $b\in L_{\d}(G)$ we have $b^{L_{\d}(\tau)}=-_{\d}(q^{-1}b^tq)$.

2)  The multiplication map $L_{\d}(G)^+\times L_{\d}(G)^-\ra L_{\d}(G)$ is a bijection. (So any element  $a\in L_{\d}(G)$ has a Cartan decomposition $a=a^++_{\d}a^-$ with respect to $\tau$.)
\end{proposition}

{\it Proof}.
 To prove 1) note that
 if $-\tau$ is the map $x\mapsto x^{-\tau}$ then
  the action of $L_{\d}(-\tau)$ on $L_{\d}(G)=Spf\ R[x']\h$ is given by the following computation, where we use the fact that $x^{t(p)}=x^{(p)t}$, $q\in N$, and $q^{(p)}=q$:
$$\begin{array}{rcl}
(x')^{L_{\d}(-\tau)} & = &\d(x^{-\tau})_{|x=1}\\
\  & = & \d(q^{-1}x^tq)_{|x=1}\\
\  & = & p^{-1}(\phi(q^{-1}x^tq)-(q^{-1}x^tq)^{(p)})_{|x=1}\\
\  & = & p^{-1}(q^{-1}(x^{(p)}+px')^tq-q^{-1}x^{t(p)}q)\\
\  & = & q^{-1}(x')^tq.\end{array}$$
We conclude that
$$(x')^{L_{\d}(\tau)}=-_{\d}((x')^{L_{\d}(-\tau)})=-_{\d}(q^{-1}(x')^tq).$$
Then 1) follows because $x'$ are the coordinates used in the identification $L_{\d}(G)={\mathfrak g}{\mathfrak l}_n$.

Let us prove 2). 
We already know this map is injective; cf. Corollary \ref{wolf}.
To prove surjectivity let $a\in L_{\d}(G)={\mathfrak g}{\mathfrak l}_n$ 
and set $A=1+pa$. We claim first that $A=UV$ with $U,V\in GL_n$, $U^{\tau}=U$, $V^{\tau}=V^{-1}$. This is standard and  analogous to the argument for the existence of the ``polar decomposition" \cite{chevalley}, p. 14. Indeed if we set
$V=(A^{-\tau}A)^{1/2}$ (cf. Equation \ref{iridor}) and $U=AV^{-1}$ then the series expression of $V$ shows that $V^{\tau}=V^{-1}$; on the other hand
$U^{-\tau}U=(V^{-\tau})^{-1}A^{-\tau}AV^{-1}$ and the latter equals $1$
if and only if $A^{-\tau}A=V^{-\tau}V$. But $V^{-\tau}V=V^2=A^{-\tau}A$ so we conclude that $U^{-\tau}U=1$ and our claim is proved.
Going back to the proof of the proposition set 
$u=\frac{U-1}{p}$ and $v=\frac{V-1}{p}$. We get that
$1+pa=(1+pu)(1+pv)$ hence $a=u+_{\d}v$.
On the other hand, by part 1), 
$$u^{L_{\d}(-\tau)}=q^{-1}u^tq=q^{-1} \frac{U^t-1}{p} q=\frac{U^{-\tau}-1}{p}=\frac{U^{-1}-1}{p}=-_{\d}u.$$
Hence 
$u^{L_{\d}(\tau)}=u$.
Similarly
$v^{L_{\d}(\tau)}=-_{\d}v$ which ends our proof. 
\qed

\medskip

In the next statement (and its proof) fractional powers are, again, principal roots as in \ref{rirdor}. Also recall that we denote by $GL_r^c$ the centralizer of $q_0:=\left(\begin{array}{rl} 0 & 1 \\ - 1 & 0\end{array}\right)$ in $GL_{2r}$.

\begin{proposition}
\label{existences}
Assume $q^t=\pm q$ and  let $x^{\tau}=q^{-1}(x^t)^{-1}q$ be the canonical  involution on $GL_n$ defining $SO(q)$. Then the following holds:
 
 i) There exists a unique lift of Frobenius $\phi_{GL_n}$  on $\widehat{GL_n}$   that is ${\mathcal H}_q$-horizontal
 and  ${\mathcal B}_q$-symmetric with respect to  $\phi_{GL_n,0}$. It is given by $\phi_{GL_n}(x)=x^{(p)}\cdot \Lambda(x)$, where
 $$\Lambda(x) =(((x^{(p)})^t\phi(q) x^{(p)})^{-1}(x^tqx)^{(p)})^{1/2}.$$
 If, in addition, $n=2r$ and $q^t=q\in GL_r^c$; then $GL_r^c$ is $\phi_{GL_n}$-horizontal. 

 Assume in what follows that $q$ is split (so $\phi_{GL_n}$ is also the unique lift of Frobenius that is ${\mathcal H}$-horizontal and ${\mathcal B}$-symmetric with respect to $\phi_{GL_{n,0}}$). Then the following hold:

ii)  $\phi_{GL_n}$ is right
compatible with $N$ and also left compatible (and hence bicompatible)  with
 $N_{SO(q)}$. 
 
 iii)   $\phi_{GL_n}$ and $\phi_{GL_n,0}$ coincide on the set $SO(q)\cap \phi_{GL_n,0}^{-1}(SO(q))$; in particular if  $\chi$ is not abnormal then the root subgroup $U_{\chi}$ is $\phi_{GL_n}$-horizontal and the lift of Frobenius on $\widehat{{\mathbb G}_a}=Spf\ R[z]\h$ induced by $\phi_{GL_n}$ is given by $z\mapsto z^p$.

iv)  If $l\d:GL_n\ra {\mathfrak g}{\mathfrak l}_n$ is the arithmetic logarithmic derivative attached to $\phi_{GL_n}$ then  for all $a\in  N_{SO(q)}$  and $b\in GL_n$ (alternatively for all $a\in GL_n$ and $b\in N$) 
  we have
  $$
\label{cris}l\d(ab)=(\phi(a)\cdot l\d(b) \cdot \phi(a)^{-1}) +_{\d} l\d(a).$$

v) 
If $l\d_0:GL_n\ra {\mathfrak g}{\mathfrak l}_n$ is the arithmetic logarithmic derivative attached to $\phi_{GL_n,0}$ then for any $a\in GL_n$ we have
$$l\d(a)\in l\d_0(a)+_{\d}{\mathfrak g}{\mathfrak l}_n^-.$$
In particular if $a\in SO(q)$ and 
$l\d_0(a)=(l\d_0(a))^++_{\d}(l\d_0(a))^-$
is the Cartan decomposition of $l\d_0(a)$ then we have
$$l\d(a)=(l\d_0(a))^+.$$

vi)  Let $\alpha\in L_{\d}(SO(q))$, let  $\phi^{\alpha}_{GL_n}$ be the lift of Frobenius
on  $\cO(GL_n)\h$ defined by $\phi^{\alpha}_{GL_n}(x)=\epsilon\cdot \phi_{GL_n}(x)$, $\epsilon=1+p\alpha$, and let $\d^{\alpha}_{GL_n}$ be the $p$-derivation on $\cO(GL_n)\h$ associated to $\phi^{\alpha}_{GL_n}$. Then $\phi_{GL_n}^{\alpha}$ is ${\mathcal H}$-horizontal with respect to $\phi_{GL_n,0}$; equivalently,
$$\d_{GL_n}^{\alpha}(\cH)=0.$$

\end{proposition}

\medskip

{\it Proof}.
Consider the matrices
\begin{equation}
\label{AB}
A  :=  (x^{(p)})^t \cdot \phi(q)\cdot  x^{(p)},\ \ \ \ B  :=  (x^t  q  x)^{(p)}.\end{equation}
 Clearly $A^t=\pm A$, $B^t=\pm B$, according as $q^t=\pm q$. Note that $A\equiv B$ mod $p$; set 
 $$C=p^{-1}(B-A).$$ 
 Define the matrix
\begin{equation}
\label{teta}
\Lambda=\Lambda(x):=(A^{-1}B)^{1/2}:=(1+pA^{-1}C)^{1/2}:=\sum_{i=0}^{\infty} \left( \begin{array}{c}1/2\\ i\end{array}\right) p^i(A^{-1}C)^i.
\end{equation}
Clearly $\Lambda\equiv 1$ mod $p$. Note that
\begin{equation}
\label{bob}
(A \Lambda)^t = \pm A\Lambda;
\end{equation}
this follows because 
$$(A(A^{-1}C)^i)^t=((CA^{-1})^{i-1}C)^t=\pm C(A^{-1}C)^{i-1}=\pm A(A^{-1}C)^i.$$
Also note that 
\begin{equation}
\label{anna}
  \Lambda^t A \Lambda= B;\end{equation}
this follows from the following computation:
$$
\begin{array}{rcl}
\Lambda^t A\Lambda & = & \left( \sum_{i=0}^{\infty}\left(\begin{array}{c} 1/2\\ i\end{array}\right)
p^i (CA^{-1})^i\right) A 
\left( \sum_{j= 0}^{\infty}\left(\begin{array}{c} 1/2\\ j\end{array}\right)
p^j (A^{-1}C)^i\right)\\
\  & \  & \  \\
\  & = & A+\sum_{n=1}^{\infty} p^n (CA^{-1})^{n-1}C \sum_{i+j=n} 
\left(\begin{array}{c} 1/2\\ i\end{array}\right)
\left(\begin{array}{c} 1/2\\ j\end{array}\right)\\
\  & \  & \ \\
\  & = & A+pC\\
\  & \  & \  \\
\  & = & B.
\end{array}
$$
Let us consider next the right action of $SO(q)$ on $R[x,\det(x)^{-1}]\h$
$$F(x) \bullet c:=F(xc),\ \ F\in R[x,\det(x)^{-1}]\h,\ \ c\in SO(q),$$
 and the left action of $SO(q)$ on $R[x,\det(x)^{-1}]\h$
$$c \bullet F(x):=F(cx),\ \ F\in R[x,\det(x)^{-1}]\h,\ \ c\in SO(q).$$
Let $u$ be in $N$ and, in case $q$ is split,  let $v$ be in  $N_{SO(q)}$.
Using Lemma \ref{monkey} one checks that:
\begin{equation}
\label{firstclue}
A\bullet u=(u^{(p)})^t A u^{(p)},\ \ v \bullet A=A,\ \ B\bullet u=(u^{(p)})^t B u^{(p)},\ \ v \bullet B=B,
\end{equation}
and, if in addition $q$ is split, 
\begin{equation}
\label{sunny}
A(u)=B(u).
\end{equation}
Hence
\begin{equation}
\label{rainy}
C\bullet u=(u^{(p)})^t C u^{(p)},\ \ v \bullet C=C,
\end{equation}
\begin{equation}
\label{firstclues}
\Lambda \bullet u=(u^{(p)})^{-1}\Lambda u^{(p)}, \ v \bullet \Lambda=\Lambda,
\end{equation}
and, if in addition $q$ is split,
\begin{equation}
\label{satisfaction}
C(u)=0,\ \ \ \Lambda(u)=1.
\end{equation}
 
 Let us prove i). 
To prove the existence part in i)  define the lift of Frobenius $\phi_{GL_n}$ by  $\phi_{GL_n}(x)=\Phi(x)=x^{(p)}\cdot \Lambda(x)$. 
The commutativity of \ref{sleepless} for $g=q$ is equivalent to the condition
$$\Phi(x)^t\cdot \phi(q)\cdot \Phi(x)=(x^tqx)^{(p)},$$
which is, of course equivalent to equation \ref{anna}.
On the other hand the commutativity of \ref{seatle} is equivalent to
$$
(x^{(p)})^t\cdot  \phi(q)\cdot  \Phi(x)=
\Phi(x)^t \cdot \phi(q) \cdot x^{(p)},$$
which is equivalent to  equality \ref{bob}. This ends the proof of the existence part.

The uniqueness part in i) follows from Corollary \ref{unique}. 

The statement about $GL_r^c$ follows from the fact that, if $q_0qq_0^{-1}=q$ then $\phi_{GL_{2r}}$ sends any $R$-point of $GL_r^c$ into a point of $GL^c_r$, as can be seen directly by using the formula in i) and the fact that conjugation by $q_0$ commutes with the maps $GL_{2r}\ra GL_{2r}$ given by  $u\mapsto u^{(p)}$,
$u\mapsto u^t$, $u\mapsto u^{-1}$, and $u\mapsto \phi(u)$.

To check ii)  we check first that $N$ is $\phi_{GL_n}$-horizontal.
For this it is enough to show that $\phi_{GL_n}(u)\in N$  if $u\in N$; but this follows from the fact that $\Lambda(u)=1$ in
\ref{firstclues}.

In order to check right compatibility of $\phi_{GL_n}$  with $N$ and left 
compatibility of $\phi_{GL_n}$ with $N_{SO(q)}$ it is sufficient to notice that, by \ref{firstclues}, we get 
$$\begin{array}{rcl}
\phi_{GL_n}(xu) & = & \phi_{GL_n}(x) \bullet u\\
\  & = & (x^{(p)}\Lambda)\bullet u\\
\  & = & (xu)^{(p)}(\Lambda \bullet u)\\
\  & = & (x^{(p)} u^{(p)})((u^{(p)})^{-1} \Lambda u^{(p)})\\
\  & = & x^{(p)}  \Lambda u^{(p)}\\
\  & = & \phi_{GL_n}(x) \cdot \phi_{N}(u)
\end{array}$$
and
$$\begin{array}{rcl}
\phi_{GL_n}(vx) & = & v \bullet \phi_{GL_n}(x)\\
\  & = & v \bullet (x^{(p)}\Lambda)\\
\  & = & (vx)^{(p)} (v\bullet \Lambda)\\
\  & = & v^{(p)} x^{(p)} \Lambda\\
\  & = & \phi_{N_{SO(q)}}(v) \cdot \phi_{GL_n}(x).  
\end{array}$$

To check assertion iii) let 
$g\in SO(q)\cap \phi_{GL_n,0}^{-1}(SO(q))$. Then $g^tqg=q$ hence $B(g)=q^{(p)}=q$; and also $\phi_{GL_n,0}(g)^t q \phi_{GL_n,0}(g)=q$, hence 
$$q=\phi(q)=\phi( \phi_{GL_n,0}(g)^t q \phi_{GL_n,0}(g))=(g^{(p)})^t\phi(q)g^{(p)}=
(g^{(p)})^t q g^{(p)},$$ so $A(g)=q$. Hence $C(g)=0$, hence $\Lambda(g)=1$,  
and hence
$$\phi_{GL_n}(g)=\phi^{-1}(x^{(p)})\phi^{-1}(\Lambda(g))=\phi_{GL_n,0}(g).$$
This ends the proof of the first assertion in iii). The second assertion in iii) follows from Lemma \ref{pat}.

Assertion  iv) follows from Lemma \ref{deal}.

Assertion v) follows from Lemma  \ref{deputy} plus the fact that if $a\in SO(q)$ then 
$l\d(a)\in L_{\d}(SO(q))\subset {\mathfrak g}{\mathfrak l}_n^+$.

To check assertion vi) note that, by the proof of assertion i), and since $q$ is split, we have
\begin{equation}
\label{noise}
q^{-1}\Phi(x)^t q \Phi(x)=(q^{-1}x^tqx)^{(p)}.
\end{equation}
We conclude by the  following computation (where we use $\epsilon\in SO(q)$):
$$
\begin{array}{rcl}
\phi_{GL_n}^{\alpha}(q^{-1}x^t qx) & = & q^{-1}(\phi_{GL_n}^{\alpha}(x))^t q \phi_{GL_n}^{\alpha}(x)\\
\   & \  & \  \\
\  & = & q^{-1}(\epsilon  \Phi(x))^t q \epsilon \Phi(x)\\
\  & \  & \  \\
\  & = & q^{-1}\Phi(x)^t \epsilon^t q \epsilon \Phi(x)\\
\  & \  & \  \\
\  & = &q^{-1} \Phi(x)^t q \Phi(x)\\
\  & \  &\  \\
\  & = & q^{-1}(x^t qx)^{(p)}\\
\  & \  & \  \\
\  & = & (q^{-1}x^tqx)^{(p)}.
\end{array}
$$
\qed

\begin{remark}
In contrast with assertion 3 in  Proposition \ref{existences} above if $\chi$ is an abnormal root (i.e. $SO(q)=SO_{2r+1}$ and $\chi=\pm \chi_i$) then one can check (by direct computation involving the explicit formula for $\phi_{GL_n}$ in the above proof) that $U_{\chi}$ is not $\phi_{GL_n}$-horizontal. \end{remark}

\begin{remark}
Here are some computations for the case $n=2$.

1) The closed subgroup schemes $SL_2$ and $Sp_2$ of $GL_2$ coincide. This by itself does not directly imply that the lifts of Frobenius $\phi_{GL_2}$ on $\widehat{GL_2}$ attached to $SL_2$ and $Sp_2$ in Propositions \ref{night} and \ref{existences} coincide. Nevertheless a trivial computation shows
that these lifts of Frobenius do indeed coincide. In other words, for
$$ q=\left(\begin{array}{rr} 0 & 1\\ -1 & 0\end{array}\right),$$  we have
$$\begin{array}{rcl}
 x^{(p)}(((x^{(p)})^tq x^{(p)})^{-1}(x^tqx)^{(p)})^{1/2}
  & = & \left( \frac{\det(x^{(p)})}{(\det(x))^p}\right)^{-1/2}x^{(p)}.\end{array}$$
 
 2) Let
$$x=\left(\begin{array}{rr} a & b\\ c & d\end{array}\right),\ \ \ q=\left(\begin{array}{rr} 0 & 1\\ 1 & 0\end{array}\right).$$ 
Then the lift of Frobenius $\phi_{GL_2}$ on $\widehat{GL_2}$ attached to $SO_2$ in Proposition \ref{existences} is given by
 $$\phi_{GL_2}(x)=\left(\begin{array}{rr} a^p & b^p\\ c^p & d^p\end{array}\right)
 \left(\begin{array}{rcl} u & v\\w & u\end{array}\right)^{1/2}$$
 where
 $$\begin{array}{rcl}
 u & = & \frac{(a^pd^p+b^pc^p)(ad+bc)^p-2^{p+1}a^pb^pc^pd^p}{(a^pd^p-b^pc^p)^2},\\
 \  & \  & \  \\
 v & = & b^pd^p\cdot \frac{2^p(a^pd^p+b^pc^p)-2(ad+bc)^p}{(a^pd^p-b^pc^p)^2},\\
 \  & \  & \  \\
 w & = & a^pc^p\cdot \frac{2^p(a^pd^p+b^pc^p)-2(ad+bc)^p}{(a^pd^p-b^pc^p)^2}.\end{array}
 $$
 
 3) Consider an arbitrary $q\in GL_2(R)$ with $q^t=q$, write
 $$q=\left(\begin{array}{rr} \alpha & \beta\\ \beta & \gamma \end{array}\right)$$
 where $\alpha,\beta,\gamma\in R$, and consider the lift of Frobenius  
 $\phi_{GL_n}(x)=x^{(p)}\Lambda(x)$ in assertion i) of Proposition \ref{existences}.
 We would like to ``compute" $\Lambda(1)$.
 First note that we have
 $$\Lambda(1)=\left( 1+\frac{p}{\alpha^p\gamma^p-\beta^{2p}}V\right)^{-1/2},\ \ \ V:=\left(
 \begin{array}{cc}
 \gamma^p\d \alpha-\beta^p\d \beta & \gamma^p\d \beta-\beta^p\d \gamma\\
 \ & \    \\
 \alpha^p\d \beta-\beta^p \d \alpha & \alpha^p\d \gamma -\beta^p \d \beta\end{array}\right).$$
 Set 
 $$\begin{array}{rcl}
 \{\alpha,\beta\}_{\d} & := & \alpha^p\d \beta-\beta^p\d \alpha,\\
 \  & \  & \  \\
 \{\alpha,\beta,\gamma\}_{\d} & := & \frac{1}{2}(\alpha^p\d \gamma-2\beta^p\d \beta+\gamma^p\d \alpha),\\
 \  & \ & \  \\
 D & :=& \{\alpha,\gamma\}^2_{\d}-4\{\alpha,\beta\}_{\d}\{\beta,\gamma\}_{\d}.\end{array}$$
 Assume $D\not\equiv 0$ mod $p$. Then the 
 eigenvalues $\lambda_1,\lambda_2$ of $V$ are in $R$ and  are non-congruent modulo $p$; they are given by
 $$\lambda_{1,2}=\{\alpha,\beta,\gamma\}_{\d}\pm\frac{\sqrt{D}}{2}.$$
  Let $u_i$
 be an eigenvector of $V$ with coefficients in $R$ for the eigenvalue $\lambda_i$;
  dividing $u_i$ by appropriate powers of $p$ we may assume $u_i\not\equiv 0$ mod $p$. Since $\lambda_1,\lambda_2$  are non-congruent modulo $p$ it follows that the reductions mod $p$ of $u_1,u_2$ are linearly independent; in other words
 the $2\times 2$ matrix $U=(u_1 u_2)$ is in $GL_2(R)$ and in particular,
 $V=U\left(\begin{array}{cc} \lambda_1 & 0\\0&\lambda_2\end{array}\right)U^{-1}$.
 We conclude that
 \begin{equation}
 \label{lambda1}
 \Lambda(1)=U
\left(\begin{array}{cc} \varphi_1^{1/2} & 0\\ 0 & \varphi_2^{1/2} \end{array}\right)
U^{-1},\ \ \ \varphi_i:=(1+\frac{p\lambda_i}{\alpha^p\gamma^p-\beta^{2p}})^{-1},\end{equation}
where the fractional powers are principal roots.

 Note that if we assume   $\alpha,\beta,\gamma$ are algebraic over ${\mathbb Q}$ then   the coefficients of $V$ are also algebraic over ${\mathbb Q}$ and hence so are $\lambda_1,\lambda_2$. Then one can arrange $U$ to have coefficients algebraic over ${\mathbb Q}$ and, hence, $\Lambda(1)$ has coefficients algebraic over ${\mathbb Q}$.
 
 Let us   assume, on the other hand,   that $\alpha,\beta,\gamma\in \bZ_p$, denote by $\left(\frac{\ \ }{p}\right):\bZ_p^{\times}\ra \{\pm 1\}$ the Legendre symbol, and 
 assume that $\left(\frac{D}{p}\right)=1$. Then $\lambda_1,\lambda_2\in \bZ_p$ and,  in particular, if $\epsilon_i\in \bZ_p$ are {\it any} elements
 satisfying $\varphi_i=\epsilon_i^{p-1}$, then 
 $$\varphi_i^{1/2}=\left(\frac{\epsilon_i}{p}\right) \cdot \epsilon_i^{(p-1)/2}.$$
 
 As a special case of the above discussion let us assume that $\alpha,\beta\in \bZ_{(p)}$, $\gamma=\alpha$,  $\alpha\not\equiv \pm \beta$ mod $p$, and $\{\alpha,\beta\}_{\d}\not\equiv 0$ mod $p$.
 Then one can take, in \ref{lambda1},
 $$\varphi_1=\frac{\alpha^p-\beta^p}{\alpha-\beta},\ \ 
 \varphi_2=\frac{\alpha^p+\beta^p}{\alpha+\beta},\ \ 
 U=\left(\begin{array}{rr}
 1 & 1\\-1 & 1\end{array}\right).
 $$
 
 4) Assume $q=1\in GL_2(R)$ is the identity matrix and let $\phi_{GL_1^c}$ be the lift of Frobenius on 
 $$\cO(\widehat{GL_1^c})=R[a,b,\frac{1}{a^2+b^2}]\h$$ induced by the lift of Frobenius $\phi_{GL_2}$ on $\widehat{GL_2}$ attached to $q$; here, as usual, we set $x=\left(\begin{array}{cc}a & b\\c & d\end{array}\right)$ and we view $\cO(GL_1^c)=R[a,b,c,d,\det(x)^{-1}]/(a-d,b+c)$. A trivial computation gives 
 \begin{equation}
 \label{crocodile}
 \phi_{GL_1^c}\left(\begin{array}{rr}a & b\\-b & a\end{array}\right)=
 \left(\begin{array}{rr}U(a,b) & V(a,b)\\-V(a,b) & U(a,b)\end{array}\right),
 \end{equation}
 \begin{equation}
 \label{alligator}
 U(a,b) = a^p K(a,b),\ \ 
  V(a,b)= b^p K(a,b),\end{equation}
  \begin{equation}
  \label{cococo}
   K(a,b) := \left(\frac{(a^2+b^2)^p}{a^{2p}+b^{2p}}\right)^{1/2}.
 \end{equation}
\end{remark}

The latter computation can be used to prove the following proposition
that settles assertion 3 in Threorem \ref{pretzel}. In this proposition we let
 $$c:GL_r^c=Spec\ R[a,b,\det(x)^{-1}]\ra GL_r=Spec\ R[v,\det(v)^{-1}]$$ be the $R$-homomorphism
 given by $v\ra  x^c:=a+\sqrt{-1}\cdot b$, where $x=\left(\begin{array}{rr}a & b\\ -b & a\end{array}\right)$. 
 
 \begin{proposition}
\label{miramar}
Let $\phi_{GL_r^c}$ be the lift of Frobenius on $\widehat{GL_r^c}$ attached to $q=1_{2r}\in GL_{2r}(R)$.
 Then there is no lift of Frobenius $\phi_{GL_r}$ on $\widehat{GL_r}$ such that $\phi_{GL_r^c}$ is $c$-horizontal with respect to $\phi_{GL_r}$. 
\end{proposition}

{\it Proof}.
We first prove our proposition in case $r=1$. 

In this case $\cO(GL_1^c)=R[a,b,\frac{1}{a^2+b^2}]$ and $\cO(GL_1)=R[v,v^{-1}]$, where $a,b,v$ are $3$ variables.
Assume there is a lift of Frobenius $\phi_{GL_1}$  such that $\phi_{GL_1^c}$ is $c$-horizontal with respect to $\phi_{GL_1}$ and seek a contradiction.
Set $\phi_{GL_1}(v)=F(v)\in R[v,v^{-1}]\h$. Then the horizontality condition gives (with notation as in \ref{crocodile}),
$$F(a+\sqrt{- 1}\cdot b)=U(a,b)+\epsilon \cdot \sqrt{- 1}\cdot V(a,b)$$
where $\phi(\sqrt{-1})=\epsilon\cdot \sqrt{-1}$, $\epsilon=\pm 1$.
Taking $\frac{\partial}{\partial a}$ and $\frac{\partial}{\partial b}$ in the previous equality we get
$$\frac{\partial F}{\partial v}(a+\sqrt{-1}\cdot b)=\frac{\partial U}{\partial a}+\epsilon \cdot \sqrt{-1}\cdot \frac{\partial V}{\partial a},$$
$$\frac{\partial F}{\partial v}(a+\sqrt{-1}\cdot b)\times \sqrt{-1}=\frac{\partial U}{\partial b}+\epsilon \cdot \sqrt{-1}\cdot \frac{\partial V}{\partial b},$$
hence
\begin{equation}
\label{permanent}
\frac{\partial U}{\partial a}=\epsilon\cdot \frac{\partial V}{\partial b},\ \ \frac{\partial V}{\partial a}=- \epsilon\cdot \frac{\partial U}{\partial b}.\end{equation}
Now by  \ref{alligator}, \ref{cococo}, one trivially gets:
\begin{equation}
\label{immanent}
\frac{\partial U}{\partial a}=pa^{p-1}K-pa^{p+1}b^2K^{-1}L,\ \ \  \frac{\partial V}{\partial b}=pb^{p-1}K+pb^{p+1}a^2K^{-1}L,\end{equation}
$$L:=\frac{(a^2+b^2)^{p-1}(a^{2p-2}-b^{2p-2})}{(a^{2p}+b^{2p})^2}.$$
Taking the image of \ref{permanent} via the map 
$$R[a,b,\frac{1}{a^2+b^2}]\h\ra R[b,b^{-1}]\h,\ \ a\mapsto 0,$$
one gets, by \ref{immanent}, that $pb^{p-1}=0$ in $R[b,b^{-1}]\h$, a contradiction. This ends the proof of the case $r=1$ of our proposition. 

Assume now that $r\geq 2$ and assume there is a lift of Frobenius $\phi_{GL_r}$  such that $\phi_{GL_r^c}$ is $c$-horizontal with respect to $\phi_{GL_r}$; again, we seek a contradiction. Consider the commutative diagram
$$\begin{array}{rcl}
GL_1^c & \stackrel{i}{\longrightarrow} & GL_r^c\\
c_1\downarrow & \  & \downarrow c_r\\
GL_1 & \stackrel{j}{\longrightarrow} & GL_r
\end{array}
$$ 
where $c_1,c_r$ are the natural projections and 
$$j(z)=z\cdot 1_r,\ \ 
i\left(\begin{array}{rr} a & b\\ -b & a \end{array}\right)=
\left(\begin{array}{rr} a\cdot 1_r & b\cdot 1_r\\  -b\cdot 1_r & a\cdot 1_r\end{array}\right).$$
Let $\phi_{GL_1^c}$ be the lift of Frobenius on $\widehat{GL_1^c}$
induced by the Chern connection attached to $1\in GL_{2r}$. Note that $i$ is horizontal with respect to $\phi_{GL_1^c}$ and $\phi_{GL_r^c}$; this follows directly from the formula 
in assertion i) of Proposition \ref{existences},
giving $\phi_{GL_r^c}$, and from the fact that the group $i(GL_1^c)$ is invariant under the operations $x\mapsto x^{(p)}$, $x\mapsto x^t$, and $(1+py)\mapsto (1+py)^{1/2}$. Since $\phi_{GL_r^c}$ sends $i(GL_1^c)$ into itself and since $c_1$ is surjective it follows that
$\phi_{GL_r}$ sends $j(GL_1)$ into itself so it induces a lift of Frobenius $\phi_{GL_1}$ on $\widehat{GL_1}$. Clearly $\phi_{GL_1^c}$ is then $c_1$-horizontal with respect to $\phi_{GL_1}$, contradicting the case $r=1$ of the proposition.
This ends the proof of the case $r\geq 2$ of the proposition.
\qed

\begin{remark}
The lift of Frobenius $\phi_{GL_n}$ in assertion i) of Proposition \ref{existences} depends, of course, on $q$; we would like to discuss, in this remark, the nature of this dependence.
Let $G_{\pm}$ be the locus of all $q\in G=GL_n$ such that $q^t=\pm q$ and for each $q\in G_{\pm}$ let $\sigma_q:\widehat{G}\ra J^1(G)$ be the section of $J^1(G)\ra \widehat{G}$
corresponding to the lift of Frobenius $\phi_{GL_n}$ attached to $q$ in assertion i) of Proposition \ref{existences}. Then by the formula in assertion i) of Proposition \ref{existences} the map of sets
$$G_{\pm}\times G\ra J^1(G),\ \ \ \ (q,a)\mapsto \sigma_q(a)$$
is a $\d$-map of order $1$ and actually comes from a morphism of $p$-formal schemes
\begin{equation}
\label{clockwork}
J^1(G_{\pm})\times \widehat{G}\ra J^1(G).\end{equation}
Explicitly, if $y$ is a matrix of indeterminates, and we view
$$G_{\pm}=Spec\ R[y]/(y^t\mp y),$$
then the map \ref{clockwork} 
 is given by the map of $R[x,\det(x)^{-1}]\h$-algebras
$$R[x,x',\det(x)^{-1}]\h\ra R[x,y,y',\det(x)^{-1}]\h/(y^t\mp y, (y')^t\mp y')$$
sending 
$$x'\mapsto \frac{1}{p}x^{(p)}\{\{(x^{(p)t}(y^{(p)}+py')x^{(p)})^{-1}(x^tyx)^{(p)}\}^{1/2}-1\}.$$
\end{remark}

\section{Inner automorphisms}

In this section we consider inner automorphisms of $GL_n$ and prove our main results about them. Throughout this section we continue to consider the lift of Frobenius $\phi_{GL_n,0}(x)=x^{(p)}$ on $\widehat{GL_n}$. 

Let us start by considering an inner  involution $x^{\tau}=q^{-1}xq$, $q\in T\subset GL_n$ with $q^2=1$, and let ${\mathcal H}(x)=q^{-1}x^{-1}qx$ and ${\mathcal B}(x,y)=q^{-1}x^{-1}qy$ be the attached quadratic and bilinear maps. (Recall that $T$ is the maximal torus of diagonal matrices.)  For such an involution  the group $G^+$ can be identified with the centralizer of $q$ and the corresponding factor $G^+\backslash G$ can be identified with the conjugacy class of $q$. One can ask if, in analogy with the case of outer involutions (Theorem \ref{laugh}), there exists a lift of Frobenius $\phi_{GL_n}$ on $\widehat{GL_n}$ that is ${\mathcal H}$-horizontal and ${\mathcal B}$-symmetric with respect to $\phi_{GL_n,0}$. In the trivial case $q=\pm 1$ the answer is yes (and then the necessarily unique $\phi_{GL_n}$ equals $\phi_{GL_n,0}$ itself).  But, except for this case, the answer is {\it no}; cf. the next proposition. 
 
 \begin{proposition}
 \label{bade}
Assume $q\in T\subset GL_n$, $q^2=1$, $q\neq \pm 1$.
Consider the involution $x^{\tau}=q^{-1}xq$ on $GL_n$, ${\mathcal H}(x)=q^{-1}x^{-1}qx$, and the lift of Frobenius $\phi_{GL_n,0}(x)=x^{(p)}$ on $\widehat{GL_n}$. 
Then
there is no lift of Frobenius $\phi_{GL_n}$ on $\widehat{GL_n}$ that is ${\mathcal H}$-horizontal  with respect to $\phi_{GL_n,0}$.
\end{proposition}

{\it Proof}.
Let
 $A$ and $B$ be defined by 
\begin{equation}
A=(x^{(p)})^{-1}qx^{(p)},\ \ \ B=(x^{-1}qx)^{(p)}.
\end{equation}
Assume $\phi_{GL_n}(x)=x^{(p)}\cdot \Lambda(x)$ is a lift of Frobenius on $\widehat{GL_n}$ that is ${\mathcal H}$-horizontal  with respect to $\phi_{GL_n,0}$ and seek a contradiction. 
By ${\mathcal H}$-horizontality we get
$$q^{-1}\Phi(x)^{-1}q\Phi(x)=(q^{-1}x^{-1}qx)^{(p)},$$
hence  
$\Lambda^{-1}A\Lambda=B$.
Taking determinants we get
\begin{equation}
\label{culprit}
\det(B)=\det(q).\end{equation}
In order to get a contradiction 
it is enough to show that \ref{culprit} fails for some value of  $x$. In order to check the latter we may
assume 
$$q=\left(
\begin{array}{ll}
q' & 0 \\
0 & q''
\end{array}
\right),\ \ q'\in GL_{n-2},\ \ q''=\left(
\begin{array}{rr}
  1 & 0\\
 0 & -1
\end{array}
\right).$$ 
Let us assume \ref{culprit} holds for all matrices $x$ of the form 
$$x=\left(\begin{array}{ccc} 1_{n-2} & 0 & 0\\
 0 & a & b\\
 0 & c & d\end{array}\right)\in GL_n(R)$$
 and seek a contradiction. Then \ref{culprit} yields the equality
 $$
 (ad+bc)^{2p}-2^{2p}a^pb^pc^pd^p=(ad-bc)^{2p}
 $$
 valid for all $a,b,c,d\in R$ with $ad-bc\in R^{\times}$. Since $R/pR$ is algebraically closed the above equality is valid as an equality between polynomials in $4$ variables $a,b,c,d$. Picking out the coefficients of $a^{2p-1}d^{2p-1}bc$ we get a contradiction. \qed

\begin{remark}
The above proposition shows that inner automorphisms $x^{\tau}=q^{-1}xq$ to not fit into the paradigm of Theorem \ref{laugh}. The next proposition shows that such inner automorphisms are, on the other hand, naturally compatible with ${\mathcal C}$-horizontality; cf. Definition \ref{late}.\end{remark} 

So let, in what follows, $D^*(x)\in R[x]$ be the discriminant of the characteristic polynomial $\det(s\cdot 1-x)\in R[x][s]$ and let $GL_n^*\subset GL_n$ be the open subset where $D^*(x)$ is invertible, i.e. $$GL_n^*=Spec\ R[x,\det(x)^{-1},D^*(x)^{-1}].$$
So $GL_n^*$ is the locus of {\it regular matrices}. Then $T^*=T\cap GL_n^*$ is the locus of {\it regular diagonal matrices}. 
If $T=Spec\ R[t,\det(t)^{-1}]$, with $t=\text{diag}(t_1,...,t_n)$ then 
$\det(t)=t_1...t_n$
and 
$$T^*=Spec\ R[t,\det(t),D(t)^{-1}],\ \ \ D(t)=\prod_{i<j}(t_i-t_j)^2.$$
Consider the morphism 
\begin{equation}
\label{defofc}
{\mathcal C}:T^*\times GL_n\ra GL_n^*,\ \ {\mathcal C}(d,a)=a^{-1}da.\end{equation}
It is defined by the map
$$R[x,\det(x)^{-1},D^*(x)^{-1}]\ra R[t,\det(t)^{-1},D^*(t)^{-1},x,\det(x)^{-1}],\ \ \ x\mapsto x^{-1}tx.$$
Note that, under this map, $\det(x)\mapsto \det(t)$, $D^*(x)\mapsto D^*(t)$.
Also consider, as usual, the lift of Frobenius $\phi_{GL_n,0}(x)=x^{(p)}$ on $\widehat{GL_n}$. Recall that $\phi_{G,0}=\phi_{GL_n,0}$ is called ${\mathcal C}$-compatible with a lift of Frobenius $\phi_{G^*}=\phi_{GL_n^*}$ on $\widehat{G^*}=\widehat{GL_n^*}$ if the diagram \ref{love} is commutative.
Then we will prove:

\begin{proposition}
\label{douche}
There exists a unique lift of Frobenius $\phi_{GL_n^*}$ on $\widehat{GL_n^*}$ such that $\phi_{GL_n,0}$ is ${\mathcal C}$-horizontal with respect to $\phi_{GL_n^*}$.
\end{proposition}

By Proposition \ref{douche} and Lemma \ref{soap} one trivially gets:

\begin{proposition}
\label{sick}
Let  $q\in T$ be regular (i.e. $q\in T^*$) with $q^{(p)}=q$, let $\tau$ be the automorphism of $GL_n$ defined by $x^{\tau}=q^{-1}xq$, and let ${\mathcal H}(x)=q^{-1}x^{-1}qx$.
Consider the lift of Frobenius $\phi_{GL_n,0}(x)=x^{(p)}$ on $\widehat{GL_n}$. Then there exists a lift of Frobenius $\phi_{GL_n}$ on $\widehat{GL_n}$ such that $\phi_{GL_n,0}$ is ${\mathcal H}$-horizontal with respect to $\phi_{GL_n}$.
\end{proposition}

To prove Proposition \ref{douche} we need some preliminaries. 
In the discussion below we will 
identify $k$-varieties ($k=R/pR$) with their sets of $k$-points. Furthermore we let $GL_{n,k}$, $T_k$, $N_k$ be the corresponding algebraic groups over $k$. Note that $N$ acts on the left on $T^*\times GL_{n}$ by the rule
\begin{equation}
\label{equat1}
n \cdot (d,a)= (ndn^{-1},na),\ \ \ n\in N,\ d\in T^*,\ a\in GL_n.\end{equation}
In particular $N_k$ acts on 
$$\cO(T_k^*\times GL_{n,k})=k[t,\det(t)^{-1},D^*(t)^{-1},x,\det(x)^{-1}].$$

\begin{lemma}\label{wind}
The ring of $N_k$-invariant elements in $\cO(T_k^*\times GL_{n,k})$ is given by
$$k[t,\det(t)^{-1},D^*(t)^{-1},x,\det(x)^{-1}]^{N_k}=k[x^{-1}tx,\det(t)^{-1},D^*(t)^{-1}].$$
\end{lemma}

{\it Proof}.
The variety $GL_{n,k}$ is a Zariski locally trivial principal homogeneous space for $T_k$ over the affine variety $U_k$ obtained from $({\mathbb P}_k^{n-1})^n$ by removing the zero locus of the determinant. (Here the elements of ${\mathbb P}_k^{n-1}$ correspond to the rows of our matrices.)
Then  
the variety $X_k:=T^*_k\times GL_{n,k}$ is a Zariski locally trivial principal homogeneous space over $Y_k:=T^*_k\times U_k$.
We have that $T_k$ acts on $X_k$ via its left action on $GL_{n,k}$, $Y_k$ is set theoretically the quotient $Y_k=T_k\backslash X_k$, and $\cO(X)^{T_k}=\cO(Y_k)$.
Consider further the quotient $Z_k=W_k \backslash Y_k$ of $Y_k$ by the finite permutation group $W_k$ acting via the action induced by that of $N_k$; so $\cO(Z_k)=\cO(Y_k)^{W_k}$, cf. \cite{serre}, p. 48.
Hence $Z_k=N_k\backslash X_k$ as sets and $\cO(Z_k)=\cO(X_k)^{N_k}$.
Set $V_k=GL_{n,k}^*$.
The morphism ${\mathcal C}_k:X_k\ra V_k$ induced by ${\mathcal C}$ 
is surjective and its set-theoretic  fibers coincide with the orbits of the $N_k$-action on $X_k$.
 Since the entries of the matrix $x^{-1}tx$ are $N_k$-invariant it follows that the map ${\mathcal C}_k$ factors as
$${\mathcal C}_k:X_k\stackrel{\alpha}{\longrightarrow} Y_k \stackrel{\beta}{\longrightarrow} Z_k \stackrel{\gamma}{\longrightarrow} V_k$$
and $\gamma$ is a bijection.

 {\it Claim 1. The scheme theoretic fiber of ${\mathcal C}_k$ through any point of the form $(d,1)\in X_k$ is smooth at $(d,1)$. }
 
 Indeed the  tangent space at $1$ of the above scheme theoretic fiber identifies with the linear space of all pairs $(a_1,d_1)$
of $n\times n$ matrices with entries in $k$ such that
$$(1+\epsilon a_1)^{-1}(d+\epsilon d_1)(1+\epsilon a_1)=d,$$
where $\epsilon^2=0$; the above condition reads
$da_1-a_1d+d_1=0$ which implies that $d_1=0$ and $a_1$ is diagonal; so the dimension of this linear space is $n$, equal to the dimension of the fiber
of ${\mathcal C}_k$ through $(d,1)$. This ends the proof of 
Claim 1.

{\it Claim 2. $\gamma$ is birational.}

 Indeed assume  $\gamma$ is not birational. Then all the scheme theoretic fibers of $\gamma$ are non-reduced. Hence all scheme theoretic fibers of ${\mathcal C}_k$ above points in a Zariski open set of $V_k$ must be non-reduced. But this contradicts Claim 1 and ends the proof of Claim 2.
 
 By Claim 2 and   Zariski's Main Theorem it follows that $\gamma$ is an isomorphism. So $\cO(V_k)=\cO(X_k)^{N_k}$
and we are done.
\qed

\medskip

{\it Proof of Proposition \ref{douche}}. 
Uniqueness is clear because the map $$\cO(GL_n^*)\h\ra \cO(T^*\times GL_n)\h$$ induced by ${\mathcal C}$  is injective. Let's prove the existence part.
Let $\phi_{GL_n^*}$  be any lift of Frobenius on $\widehat{GL_n^*}$, $\phi_{GL_n^*}(x)=x^{(p)}\Lambda(x)$ 
where $\Lambda$ is an $n\times n$ matrix with entries in $R[x,\det(x)^{-1},D^*(x)^{-1}]\h$, $\Lambda\equiv 1$ mod $p$.
Set $$A(t,x)=(x^{(p)})^{-1}t^{(p)}x^{(p)}, B(t,x)=(x^{-1}tx)^{(p)}\in \cO(T^*\times GL_n).$$
The condition that $\phi_{GL_n^*}$ is ${\mathcal C}$-horizontal with respect to $\phi_{GL_n,0}$ translates into the equality
\begin{equation}
\label{functionaleq}
\Lambda(x^{-1}tx)=B(t,x)^{-1}A(t,x)\end{equation}
in the ring $R[t,\det(t)^{-1}, D^*(t)^{-1},x,\det(x)^{-1}]\h$.
Note that $A(t,x)$ and $B(t,x)$ are $N$-invariant for the action \ref{equat1}.
We shall construct a sequence of matrices $\Lambda_{\nu}(x)$ with entries in $R[x,\det(x)^{-1},D^*(x)^{-1}]\h$, with $\Lambda_0=1$,  satisfying the following properties:

1) $\Lambda_{\nu+1}(x)\equiv \Lambda_{\nu}(x)$ mod $p^{\nu+1}$ for $\nu\geq 0$;

2) $\Lambda_{\nu}(x^{-1}tx)\equiv B(t,x)^{-1}A(t,x)$ mod $p^{\nu+1}$ for $n\geq 0$.

This will end the proof by taking 
\begin{equation}
\label{Lambda}
\Lambda(x)=\lim \Lambda_{\nu}(x).\end{equation}
 To construct $\Lambda_{\nu}$ we proceed by induction. Assume $\Lambda_{\nu}$ has been constructed for some ${\nu}\geq 0$. Write
$$\Lambda_{\nu}(x^{-1}tx)=B(t,x)^{-1}A(t,x)+p^{\nu+1}C(t,x).$$
Then $C(t,x)$ is $N$-invariant. By Lemma \ref{wind} there exists $F\in R[x,\det(x)^{-1},D^*(x)^{-1}]$ such that 
$$C(x)\equiv F(x^{-1}qx)\ \ \ \text{mod}\ \ \ p.$$
Setting 
$\Lambda_{\nu+1}=\Lambda_{\nu}(x)-p^{\nu+1}F(x)$
 one immediately checks that
$$\Lambda_{\nu+1}(x^{-1}qx)\equiv B(x)^{-1}A(x)\ \ \ \text{mod}\ \ \ p^{\nu+2}$$
which ends our construction and hence the proof of the proposition.
\qed

\medskip

 In contrast with the situation for regular matrices described in  Proposition \ref{sick}, we have the following results for non-regular matrices: 
 
\begin{proposition}
\label{sickness}
Let  $q\in T$ with $q^{(p)}=q$, let $\tau$ be the automorphism of $GL_n$ defined by $x^{\tau}=q^{-1}xq$, and let ${\mathcal H}(x)=q^{-1}x^{-1}qx$.
Consider the lift of Frobenius $\phi_{GL_n,0}(x)=x^{(p)}$ on $\widehat{GL_n}$. Assume $q$ is not regular (i.e. $q\not\in T^*$)  and not scalar. Then there exists no lift of Frobenius $\phi_{GL_n}$ on $\widehat{GL_n}$ such that $\phi_{GL_n,0}$ is ${\mathcal H}$-horizontal with respect to $\phi_{GL_n}$.
\end{proposition}

{\it Proof}.
Let $\phi_{GL_n}$ be any lift of Frobenius on $\widehat{GL_n}$, $\phi_{GL_n}(x)=x^{(p)}\Lambda(x)$ and set $$A(x)=(x^{(p)})^{-1}qx^{(p)},\ \ \ B(x)=(x^{-1}qx)^{(p)}.$$
The condition that $\phi_{GL_n,0}$ is ${\mathcal H}$-horizontal with respect to $\phi_{GL_n}$ means by definition that the following diagram is commutative:
$$\begin{array}{rcl}
\widehat{GL_n} & \stackrel{\phi_{GL_n,0}}{\longrightarrow} & \widehat{GL_n}\\
\cH \downarrow & \  & \downarrow \cH\\
\widehat{GL_n} & \stackrel{\phi_{GL_n}}{\longrightarrow} & \widehat{GL_n}
\end{array}$$
 Hence this condition translates into the following equality:
\begin{equation}
\label{impo}
\Lambda(q^{-1}x^{-1}qx)=B(x)^{-1}A(x).\end{equation}
Let $q\in T$ be not regular and not scalar, with $q^{(p)}=q$. Let us assume  \ref{impo} and seek a contradiction. We may assume $$q=\text{diag}(q_1\cdot 1_{r_1},...,q_s\cdot 1_{r_s})$$ with $q_1,...,q_s$ distinct and $s\geq 2$, $r_1\geq 2$.
It is enough then to obtain a contradiction in the case $n=3$, $r_1=2$, $r_2=1$.
Let $H\subset GL_3$ be the group of matrices of the form $\left(\begin{array}{cc}
a & b \\ c & d\end{array}\right)$ with $a\in GL_2$; hence $q\in T\cap H$ and the centralizer of $q$ in $GL_3$ is $H$. Note that $B(hx)=B(x)$ and $\Lambda(q^{-1}(hx)^{-1}q(hx))=\Lambda(q^{-1}x^{-1}qx)$ for all $h\in H$. From \ref{impo} it follows that
$A(hx)=A(x)$ for $h\in H$, hence
\begin{equation}
((hx)^{(p)})^{-1}q(hx)^{(p)}=(x^{(p)})^{-1}qx^{(p)},\ \ h\in H.\end{equation}
This implies that $(hx)^{(p)}(x^{(p)})^{-1}$ is in the centralizer of $q$ and hence 
\begin{equation}
\label{aaron}
(hx)^{(p)}(x^{(p)})^{-1}=\left(\begin{array}{cc}
w_{11} & 0 \\ 0 & w_{22}\end{array}\right)
\end{equation}
with $w_{ij}$ functions of $x$, $w_{11}$ a $2\times 2$ matrix.
Write 
$$h=\left(\begin{array}{cc}
h_{11} & 0 \\ 0 & h_{22}\end{array}\right), \ \ x=\left(\begin{array}{cc}
x_{11} & x_{12} \\ x_{21} & x_{22}\end{array}\right),$$
with $h_{11},x_{11}$ $2 \times 2$ matrices. Then \ref{aaron} yields:
$$((h_{11}x_{11})^{(p)})^{-1}(h_{11}x_{12})^{(p)}=(x_{11}^{(p)})^{-1}x_{12}^{(p)}$$
for any $h_{11}\in GL_2$. Let $x^1_{12},x^2_{12}$ be column vectors of indeterminates and let $x_{12}^*$ be the $2\times 2$ matrix with columns $x^1_{12},x^2_{12}$.
Also let $y_{11}$ be a $2\times 2$ matrix of indeterminates. Then we get
$$((y_{11}x_{11})^{(p)})^{-1}(y_{11}x^*_{12})^{(p)}=(x_{11}^{(p)})^{-1}(x^*_{12})^{(p)}$$
Setting $x^*_{12}=y_{11}^{-1}$ and $x_{11}=1$ we get
$$((x_{12}^*)^{-1})^{(p)})^{-1}=(x^*_{12})^{(p)},$$
which is easily seen to be a contradiction.
\qed

\bigskip

By Proposition \ref{sickness} and (the proof of) Proposition \ref{douche}
we will get the following:

\begin{proposition}
\label{leak}
Assume $n\geq 3$. Then the lift of Frobenius $\phi_{GL_n^*}$ on $\widehat{GL_n^*}$ in Proposition \ref{douche} does not extend to a lift of Frobenius on $\widehat{GL_n}$.
\end{proposition}

{\it Proof}.
Assume the 
lift of Frobenius $\phi_{GL_n^*}$ on $\widehat{GL_n^*}$  extends to a lift of Frobenius $\phi_{GL_n}$ on $\widehat{GL_n}$ and seek a contradiction.
We get that the matrix $\Lambda(x)$ in \ref{Lambda} has coefficients in $R[x,\det(x)^{-1}]\h$. Hence the equation \ref{functionaleq} holds in the ring $R[t,\det(t)^{-1},x,\det(x)^{-1}]\h$ because the map $$R[t,\det(t)^{-1},x,\det(x)^{-1}]\h\ra R[t,\det(t)^{-1},D^*(t)^{-1},x,\det(x)^{-1}]\h$$
is injective. Now pick any $q\in T$ with $q^{(p)}=q$  which is not regular and not scalar; this is possible because $n\geq 3$. Let $\tau$ be the automorphism of $GL_n$ defined by $x^{\tau}=q^{-1}xq$. Setting $t=q$ in equality \ref{functionaleq} we get an equality 
$$\Lambda(x^{-1}qx)=B(q,x)^{-1}A(q,x)$$
in the ring $R[x,\det(x)^{-1}]\h$. But this is immediately seen to imply that
the lift of Frobenius $\phi_{GL_n}$ on $\widehat{GL_n}$ defined by 
$\phi_{GL_n}(x)=x^{(p)}\Lambda(qx)$ has the property that $\phi_{GL_n,0}$ is ${\mathcal H}$-horizontal with respect to $\phi_{GL_n}$. This, however, contradicts Proposition \ref{sickness}, which ends our proof.
\qed

\bigskip

We examine now the interaction between lifts of Frobenius and the map
$${\mathcal P}:GL_n\ra {\mathbb A}^n=Spec\ R[z_1,...,z_n]$$ 
with coordinates ${\mathcal P}_1,...,{\mathcal P}_n\in R[x]$ given by the formula \ref{characteristic}; so 
$${\mathcal P}_1(x)=\text{tr}(x),...,{\mathcal P}(x)=\det(x).$$
Recall from the Introduction that
we are considering the lift of Frobenius $\phi_{{\mathbb A}^n,0}$ on $\widehat{{\mathbb A}^n}$ defined by $\phi_{{\mathbb A}^n,0}(z_i)=z_i^{p}$ and  we say that a lift of Frobenius $\phi_{G^{**}}$ on an open set $\widehat{G^{**}}$ of $\widehat{G}=\widehat{GL_n}$ is ${\mathcal P}$-horizontal with respect to $\phi_{{\mathbb A}^n,0}$ if the diagram \ref{loving} in the Introduction is commutative. On the other hand we will say that $\phi_{G^{**}}$ is a {\it $T$-deformation} of $\phi_{G,0}(x)=x^{(p)}$ if there is a commutative diagram \ref{hating}; cf. the Introduction. If  $G^{**}$ is invariant under conjugation by $N$ we say $\phi_{G^{**}}$ is ${\mathcal C}$-compatible with $N$ if the diagram \ref{drinking} in the Introduction is commutative.

\begin{proposition}
\label{crazy}
 \ 
 
 1) There exists an open set $G^{**}$ of $G=GL_n$  and  a  lift of Frobenius $\phi_{G^{**}}$ on $\widehat{G^{**}}$ such that $G^{**}$ is invariant under conjugation by $N$,  $G^{**}\cap T=T^*$,
  $\phi_{G^{**}}$
   is ${\mathcal P}$-horizontal with respect to $\phi_{{\mathbb A}^n,0}$, and   $\phi_{G^{**}}$ is a $T$-deformation of $\phi_{G,0}$.
  
  2) For any $G^{**}$ as above, $\phi_{G^{**}}$ is unique with the above properties and  is ${\mathcal C}$-compatible with $N$. 
  
  3) For any matrix $\alpha(x)$ with coefficients in $\cO(G^{**})\h$  the lift of Frobenius $\phi_{G^{**}}^{(\alpha)}$ on $\widehat{G^{**}}$ defined by 
   $$\phi_{G^{**}}^{(\alpha)}(x)=\epsilon(x)\cdot \phi_{G^{**}}(x) \cdot \epsilon(x)^{-1}, \ \ \epsilon(x):=1+p\alpha(x)$$
   is ${\mathcal P}$-horizontal with respect to $\phi_{{\mathbb A}^n,0}$. So if $\d^{(\alpha)}_{G^{**}}$ is the $p$-derivation on $\cO(G^{**})\h$ attached to $\phi_{G^{**}}^{(\alpha)}$ then 
   $$\d^{(\alpha)}_{G^{**}} {\mathcal P}_i=0,\ \ i=1,...,n.$$
\end{proposition}

{\it Proof}. Consider the action of $N$ on $R[x\det(x)^{-1}]\h$ by conjugation,
$n\star x=nxn^{-1}$, $n\in N$. 
Let $e_1,...,e_n$ be the columns of the $n\times n$ identity matrix, i.e.,  $1_n=[e_1,...,e_n]$, and let
$$1_{n,j}=[e_1,...,e_{j-1},0,e_{j+1},...,e_n],$$
$$\det(s\cdot 1_{n,j}-x)=\sum_{i=0}^{n-1}(-1)^i{\mathcal P}_{ij}(x)s^{n-i}\in R[x][s],$$
$$D^{**}(x):=\det({\mathcal P}_{ij}(x))\in R[x].$$
Clearly 
\begin{equation}
\label{**P}
d\star {\mathcal P}_{ij}(x)={\mathcal P}_{ij}(x),\ \ \ d\star D^{**}(x)=D^{**}(x)\end{equation}
 for $d\in T$. One also easily checks that
 \begin{equation}
 \label{kelly}
 w_{\sigma} \star {\mathcal P}_{ij}(x)={\mathcal P}_{i \sigma^{-1}(j)}(x),\ \ 
 w_{\sigma}\star D^{**}(x)=\pm D^{**}(x),
 \end{equation}
 for all $w_{\sigma}=[e_{\sigma(1)},...,e_{\sigma(n)}]\in W$.
Define $G^{**}$ as the locus in $G$ where $D^{**}$ is invertible; so $G^{**}$ is invariant under conjugation by $N$.
Note that if $\mu=\text{diag}(\mu_1,...,\mu_n)\in T$ then
$D^{**}(\mu)\in R^{\times}$ if and only if the polynomials 
$$\prod_{i\neq 1}(s-\mu_i),...,\prod_{i\neq n}(s-\mu_i)$$
are linearly independent mod $p$ which holds if and only if $\mu_1,...,\mu_n$ are distinct mod $p$; this shows that $G^{**}\cap T=T^*$. Let $\phi_{G^{**}}$ be a lift of Frobenius on $\widehat{G^{**}}$,
$\phi_{G^{**}}(x)=x^{(p)}\Lambda(x)$, where $\Lambda(x)$ is a matrix $\equiv 1$ mod $p$ with coefficients in $R[x,\det(x)^{-1},D^{**}(x)^{-1}]\h$. The commutativity of \ref{loving} is equivalent to the equalities 
\begin{equation}
\label{white}
{\mathcal P}_i(x^{(p)}\Lambda(x))={\mathcal P}_i(x)^p\end{equation}
 for all $i=1,...,n$, hence to the equality
$$V(x,s):=\sum_{i=0}^n(-1)^i{\mathcal P}_i(x)^ps^{n-i}=\sum_{i=0}^n(-1)^i{\mathcal P}_i(x^{(p)}\Lambda(x))s^{n-i},$$
hence to the equality
$$V(x,s^p)=\sum_{i=0}^n(-1)^i{\mathcal P}_i(x^{(p)}\Lambda(x))s^{p(n-i)}.$$
Note that, by \ref{white}, if $\phi_{G^{**}}$ is
${\mathcal P}$-horizontal with respect to $\phi_{{\mathbb A}^n,0}$, then
for any matrix $\alpha(x)$ with coefficients in $\cO(G^{**})\h$ the lift of Frobenius $\phi_{G^{**}}^{(\alpha)}$ in the statement of the proposition
   is ${\mathcal P}$-horizontal with respect to $\phi_{{\mathbb A}^n,0}$.
Note also that  $\phi_{G^{**}}$ is a $T$-deformation of $\phi_{G,0}$ if and only if $\Lambda(x)$ is a diagonal matrix. Finally one easily checks that $\phi_{G^{**}}$   is ${\mathcal C}$-compatible with $N$ if and only if $d\star\Lambda(x)=\Lambda(x)$ for $d\in T$ and $w\star \Lambda(x)=w\Lambda(x)w^{-1}$ for $w\in W$. In what follows we prove the existence of a diagonal matrix $\Lambda(x)$ as above.
The uniqueness part of the proposition  is proved via a similar argument. To prove the existence of $\Lambda(x)$ it is enough to construct a sequence $\Lambda_{\nu}(x)$ of diagonal matrices satisfying the following properties:

1) $\det(s^p\cdot 1_n-x^{(p)}\Lambda(x))\equiv V(x,s^p)$ mod $p^{\nu+1}$, $\nu\geq 0$.

2) $\Lambda_{\nu+1}(x)\equiv \Lambda_{\nu}(x)$ mod $p^{\nu+1}$, $\nu\geq 0$.

3) $d\star \Lambda_{\nu}(x)=\Lambda_{\nu}(x)$ for $d\in T$.

4) $w\star \Lambda_{\nu}(x)=w\Lambda_{\nu}(x)w^{-1}$ for $w\in W$.

We take $\Lambda_0(x)=1$. Assume now $\Lambda_{\nu}(x)$ was constructed and seek $\Lambda_{\nu+1}(x)$ of the form $\Lambda_{\nu}(x)-p^{\nu+1}Z(x)$ where $Z(x)$ is a diagonal matrix. Write
$$\begin{array}{rcl}
A & = & s^p\cdot 1_n-x^{(p)}\Lambda_{\nu}(x)=(a_{ij}),\\
\  & \  & \  \\
\det(A) & = & V(x,s^p)+p^{\nu+1}U(x,s^p),\\
\  & \  & \  \\
U(x,s) & = & \sum_{i=0}^{n-1}U_i(x)s^i,\\
\  & \  & \  \\
Z & = & \text{diag}(z_1,...,z_n),\\
\  & \  & \  \\
 B & = & x^{(p)}Z(x)=(b_{ij})=(z_jx_{ij}^p).\end{array}$$
 Since $n\star V(x,s)=V(x,s)$ and $n\star A=A$ for $n\in N$ it follows that
 \begin{equation}
 \label{**U}
 n\star U_i(x)=U_i(x)\end{equation}
  for all $i$ and $n\in N$.
Then $$\det(s^p\cdot 1_n-x^{(p)}\Lambda_{\nu+1}(x))= \det\left(
\begin{array}{ccc}
a_{11}+p^{\nu+1}b_{11} & a_{12}+p^{\nu+1}b_{12} & *\\
a_{21}+p^{\nu+1}b_{21} & a_{22}+p^{\nu+1}b_{22} & *\\
* & * & *\end{array}\right)
$$
which is $\equiv$ mod $p^{\nu+1}$ to
$$\det(A)+p^{\nu+1}W(x,s^p)$$
where
$$
W(x,s^p)=\det\left(
\begin{array}{lll}
b_{11} & a_{12} & *\\
b_{21} & a_{22} & *\\
* & * & *\end{array}\right)
+\det\left(
\begin{array}{lll}
b_{11} & a_{12} & *\\
b_{21} & a_{22} & *\\
* & * & *\end{array}\right)+...
$$
So, for property 2), it is enough to find $Z(x)$ such that
\begin{equation}
\label{UW}U(x,s^p)+W(x,s^p)\equiv 0\ \ \ \text{mod}\ \ p.\end{equation}
On the other hand we have $A\equiv s\cdot 1_n-x^{(p)}$ mod $p$ so
$W(x,s^p)$ is $\equiv$ mod $p$ to
$$-z_1\det\left(\begin{array}{rrr}-x_{11} & -x_{12} & *\\
-x_{21} & s-x_{22} & *\\
* & * & *\end{array}\right)^p-z_2\det\left(\begin{array}{rrr}s-x_{11} & -x_{12} & *\\
-x_{21} & -x_{22} & *\\
* & * & *\end{array}\right)^p-...$$
Now  \ref{UW} is satisfied if the following system can be solved:
\begin{equation}
\label{bear}
\sum_{j=1}^n {\mathcal P}_{ij}(x)^p z_j= (-1)^iU_i(x),\ \ \ i=0,...,n-1.\end{equation}
  This system has  determinant  $\equiv  D^{**}(x)^p$ mod $p$ so it has a (unique) solution $(z_1,...,z_n)$ in the ring $R[x,D^{**}(x)^{-1}]\h$. By Cramer's rule and \ref{**P}, \ref{**U}, we have that $d\star z_j=z_j$ for all $j$ and $d\in T$. On the other hand by applying $w_{\sigma}$ to \ref{bear} we get
  $$\sum_{j=1}^n {\mathcal P}_{i\sigma^{-1}(j)}(x)^p (w_{\sigma}\star z_j)= (-1)^iU_i(x)$$
  hence $w_{\sigma}\star z_j=z_{\sigma^{-1}(j)}$, hence $w_{\sigma}\star Z=w_{\sigma}Zw_{\sigma}^{-1}$, hence $$w_{\sigma}\star \Lambda_{\nu+1}(x)=w_{\sigma}\Lambda_{\nu+1}(x)w_{\sigma}^{-1}.$$ This ends our construction and hence ends our proof.
\qed

\begin{remark}
\label{death}
By assertion 3 in Proposition \ref{crazy} and Lemma \ref{incaolema} it follows that ${\mathcal P}_i$ are prime integrals for the system of arithmetic differential equations 
$$\d x_{ij}-\Delta_{ij}^{**(\alpha)}(x),$$
where $\Delta_{ij}^{**(\alpha)}(x)=\d_{G^{**}}(x_{ij})$; i.e., if $\Delta^{**(\alpha)}=(\Delta^{**(\alpha)}_{ij})$ then for any solution $u\in G^{**}$ of the equation
\begin{equation}
\label{doggy}
\d u=\Delta^{**(\alpha)}(u)
\end{equation}
we have $\d({\mathcal P}_i(u))=0$ for all $i$. On the other hand we claim that Equation \ref{doggy} is equivalent to
\begin{equation}
\label{kitty}
l\d u=-_{\d}(u\star_{\d} \alpha(u))+_{\d}\alpha(u),
\end{equation}
where $l\d:G^{**}\ra {\mathfrak g}{\mathfrak l}_n$ is defined by
$$l\d u:=\frac{1}{p}(\phi(u)\Phi^{**}(u)^{-1}-1,\ \ \ \Phi^{**}(x):=\phi_{G^{**}}(x).$$
(This makes assertion 3 in our Proposition \ref{crazy} an analogue of the corresponding claim about the usual differential  Equation \ref{moon}.) 
To check the claim note that Equation \ref{kitty} is equivalent to 
$$1+pl\d u=(1+p\phi(u) \alpha(u)\phi(u)^{-1})^{-1}(1+p\alpha(u))$$
hence with
$$\Phi^{**}(u)=(1+p\alpha(u))^{-1}\phi(u)(1+p\alpha(u))=\epsilon(u)^{-1}\phi(u)\epsilon(u),$$
hence with Equation \ref{doggy}, and our claim is proved.
\end{remark}

\end{document}